# VERTEX-REINFORCED RANDOM WALK ON $Z$ EVENTUALLY GETS STUCK ON FIVE POINTS[1]

### By Pierre Tarrès

*CNRS, Université Paul Sabatier and Université de Neuchâtel*

Vertex-reinforced random walk (VRRW), defined by Pemantle in 1988, is a random process that takes values in the vertex set of a graph $G$, which is more likely to visit vertices it has visited before. Pemantle and Volkov considered the case when the underlying graph is the one-dimensional integer lattice $\mathbb{Z}$. They proved that the range is almost surely finite and that with positive probability the range contains exactly five points. They conjectured that this second event holds with probability 1. The proof of this conjecture is the main purpose of this paper.

**1. General introduction.** Let $(\Omega, \mathcal{F}, \mathbb{P})$ be a probability space. Let $G$ be a locally finite graph, let $\sim$ be its neighbor relationship and let $V(G)$ be its vertex set. Let $(X_n)_{n \in \mathbb{N}}$ be a process that takes values in $V(G)$. Let $\mathbb{F} = (\mathcal{F}_n)_{n \in \mathbb{N}}$ denote the filtration generated by the process [i.e., $\mathcal{F}_n = \sigma(X_0, \ldots, X_n)$ for all $n \in \mathbb{N}$] and let $\mathcal{F}_\infty = \sigma(\mathcal{F}_n, n \geq 0)$.

For any $v \in V(G)$, let $Z_n(v)$ be the number of times plus 1 that the process visits site $v$ up through time $n \in \mathbb{N} \cup \{\infty\}$, that is,

$$Z_n(v) = 1 + \sum_{i=0}^{n} \mathbb{1}_{\{X_i = v\}}.$$

Then $(X_n)_{n \in \mathbb{N}}$ is called vertex-reinforced random walk (VRRW) with starting point $v_0 \in V(G)$ if $X_0 = v_0$ and for all $n \in \mathbb{N}$,

$$\mathbb{P}(X_{n+1} = x | \mathcal{F}_n) = \mathbb{1}_{\{x \sim X_n\}} \frac{Z_n(x)}{\sum_{w \sim X_n} Z_n(w)}.$$

Received July 2002; revised April 2004.

[1]Supported in part by the Swiss National Science Foundation Grant 200021-1036251/1.

In other words, moves are restricted to the edges of $G$, with the probability of a move to a neighbor $x$ being proportional to the augmented occupation $Z_n(x)$ of $x$ at that time.

VRRWs were introduced in 1988 by Pemantle [7] in the spirit of the seminal work by Coppersmith and Diaconis [4], who defined the notion of edge-reinforced random walks, which have at each step a probability to move along an edge proportional to the number of times plus 1 that the process has visited this edge. Reinforced processes are useful in models involving self-organization and learning behavior; they can also describe spatial monopolistic competition in economics. For more details on applications and known results in connection with these models, refer to the articles by Pemantle and Volkov [8, 9].

VRRWs on finite complete graphs, with reinforcements weighted by factors associated to each edge of the graph, have been studied by Pemantle [8] and Benaïm [1]. Pemantle and Volkov obtained results in 1997 on reinforced random walks on $\mathbb{Z}$ [9], which are described in the following text. More recently, Volkov [13] generalized some of these results and proved that, on a fairly broad class of locally finite graphs (containing the graphs of bounded degree), the VRRW has finite range with positive probability. The remainder of this paper is devoted to VRRWs on $\mathbb{Z}$.

Define the two random sets

$$R := \{v \in \mathbb{Z} / \exists\, n \in \mathbb{N} \text{ s.t. } X_n = v\},$$

$$R' := \{v \in \mathbb{Z} / X_n = v \text{ infinitely often}\}$$

and, given $k \in \mathbb{Z}$ and $\alpha \in (0, 1)$, define the six events:

1. $\{R' = \{k-2, k-1, k, k+1, k+2\}\}$;
2. $\{\ln Z_n(k-2)/\ln n \to \alpha\}$;
3. $\{\ln Z_n(k+2)/\ln n \to 1 - \alpha\}$;
4. $\{Z_n(k-1)/n \to \alpha/2\}$;
5. $\{Z_n(k+1)/n \to (1-\alpha)/2\}$;
6. $\{Z_n(k)/n \to 1/2\}$.

Let $|\cdot|$ be the cardinality of a set. Pemantle and Volkov [9] proved the following results.

THEOREM 1.1.    *One has* $\mathbb{P}(|R| < \infty) = 1$ *and* $\mathbb{P}(|R| = 5) > 0$.

THEOREM 1.2.    *One has* $\mathbb{P}(|R'| \leq 4) = 0$.

THEOREM 1.3.    *For any open set* $I \subset [0, 1]$ *and any integer* $k \in \mathbb{Z}$ *there exists, with positive probability,* $\alpha \in I$ *such that events 1–6 occur.*

Pemantle and Volkov also proposed the following conjecture.

CONJECTURE 1. *There exist almost surely $k \in \mathbb{Z}$ and $\alpha \in (0,1)$ such that events 1–6 occur.*

The main purpose of the present article is to prove this conjecture. In fact, we prove the following result, which is slightly more accurate. Given $C_1$, $C_2 \in (0, \infty)$ and $k \in \mathbb{Z}$, define the two events:

2′. $Z_n(k-2)/n^\alpha \to C_1$;
3′. $Z_n(k+2)/n^{1-\alpha} \to C_2$.

THEOREM 1.4. *There exist almost surely $k \in \mathbb{Z}$, $\alpha \in (0,1)$ and $C_1$, $C_2 \in (0, \infty)$ such that events 1, 2′, 3′ and 4–6 occur.*

In our proof, we make use of Theorem 1.1 of Pemantle and Volkov ([9]; stated above). The heuristic developed by these authors on the comparison of VRRW to Pólya and Friedman urn models [9, 13] also has been very useful and is partly related to the results claimed in Section 3.1. Although we do not use it explicitly, the heuristic of a result from Benaïm about convergence with positive probability toward an attractor ([2], Chapter 7) has been very useful in Lemmas 2.4, 2.9 and 2.11.

## 2. Introduction to the ideas of the proof.

2.1. *Notation.* Let $\mathbb{R}^* = \mathbb{R} \setminus \{0\}$, $\mathbb{N}^* = \mathbb{N} \setminus \{0\}$, $\mathbb{R}_+^* = \mathbb{R}_+ \setminus \{0\}$ and $\mathbb{Q}_+^* = \mathbb{Q}_+ \setminus \{0\}$. Given a random sequence $(\gamma_n)$ of $\mathbb{F}$-adapted nondecreasing stopping times ($\forall q \in \mathbb{N}$, $\{\gamma_n \leq q\} \in \mathcal{F}_q$), let $\mathbb{F}_{(\gamma_n)_{n \in \mathbb{N}}} = (\mathcal{F}_{\gamma_n})_{n \in \mathbb{N}}$ denote the filtration defined as follows: For all $n \in \mathbb{N}$, $A \in \mathcal{F}_{\gamma_n} \iff \forall q \in \mathbb{N} \cup \{\infty\}$, $A \cap \{\gamma_n \leq q\} \in \mathcal{F}_q$. The equalities and inclusions between probability events are understood to hold almost surely. Given $x$, $y \in \mathbb{R}$, we use alternately the notation $x \wedge y$ and $\min(x, y)$ [resp. $x \vee y$ and $\max(x, y)$] for the minimum [resp. the maximum] of $x$ and $y$. We write $x = \square(y)$ iff $|x| \leq y$. We let $x^+ = \max(x, 0)$ and $x^- = \max(-x, 0)$. Let $\mathrm{Cst}(a_1, a_2, \ldots, a_p)$ denote a positive constant depending only on $a_1$, $a_2, \ldots, a_p$ and let Cst denote a universal positive constant. We say for simplicity that a property holds for $x < \mathrm{Cst}(a_1, \ldots, a_p)$ [resp. for $x > \mathrm{Cst}(a_1, \ldots, a_p)$] when there exists a constant $c$, which depends only on $a_1, \ldots, a_p$ so that this property holds for $x < c$ (resp. for $x > c$).

Let $(u_n)_{n \in \mathbb{N}}$ and $(v_n)_{n \in \mathbb{N}}$ be two sequences taking values in $\mathbb{R}$. We write $u_n = O(v_n)$ [resp. $u_n = o(v_n)$] when there exists an a.s. finite random variable $C$ [resp. a random sequence $(C_n)_{n \in \mathbb{N}}$ converging to 0 a.s.] such that, for all $n \in \mathbb{N}$, $u_n \leq Cv_n$ [resp. $u_n \leq C_n v_n$].

We write $u_n \asymp v_n$ iff either $\limsup |u_n| < \infty$ and $\limsup |v_n| < \infty$ or $u_n/v_n \to 1$, and write $u_n \preceq v_n$ iff, for all $\varepsilon > 0$, there exists $k_0 \in \mathbb{N}$ such that, for all $n \geq k \geq k_0$,

$$u_n - u_k \leq (1 + \varepsilon)(v_n - v_k) + \varepsilon.$$

Note that, if $u_n$ and $v_n$ are random variables, $k_0$ is a priori a random variable.

Similarly, given $a \in \mathbb{R}$ and another sequence $(w_n)_{n \in \mathbb{N}}$ that takes values in $\mathbb{R}$, we write $u_n \preceq_{w_n \geq a} v_n$ when, for all $\varepsilon > 0$, there exists $k_0 \in \mathbb{N}$ such that, for all $n \geq k \geq k_0$,

$$u_n - u_k \leq (1 + \varepsilon)(v_n - v_k) + \varepsilon \qquad \text{if } w_m \geq a \text{ for all } m \in [k, n].$$

We write $u_n \equiv v_n$ iff $\lim(u_n - v_n)$ exists a.s. and is finite, and write $u_n \doteq v_n$ iff there exists a random $k_0 \in \mathbb{N}$ such that for all $n \geq k_0$, $u_n - v_n = u_{k_0} - v_{k_0}$. In particular, we write $u_n \doteq v_n + o(w_n)$ iff there exists $\alpha \in \mathbb{R}$ (a priori random) such that $u_n - v_n = \alpha + o(w_n)$. Given $u, v \in \mathbb{R}_+^* \cup \{\infty\}$, we write $u \approx v$ iff either $u = v = \infty$ or $\max(u, v) < \infty$, and write $u \succeq v$ iff either $u = \infty$ or $\max(u, v) < \infty$.

We let $\mathbb{E}[\cdot]$ and $\mathbb{V}[\cdot]$ be the expectation and the variance of a random variable. If $\mathcal{G}$ is a sub-$\sigma$-field of $\mathcal{F}$, we let $\mathbb{E}[\cdot | \mathcal{G}]$ and $\mathbb{V}[\cdot | \mathcal{G}]$ be the expectation and variance conditionally to $\mathcal{G}$.

2.2. *Sketch of the proof.* Let us begin with some background on the study of VRRWs. First recall that VRRW are non-Markovian processes. Define, for all $n \in \mathbb{N}$, the vector of occupation densities of the random walk at time $n$ as

$$V(n) = \left( \frac{Z_n(v)}{n} \right)_{v \in V(G)}.$$

The works of Pemantle [8] and Benaïm [1] provide some methods to compare the behavior of $V(n)$ with solutions of ordinary differential equations when the graph is complete (i.e., any two vertices and adjacent). The heuristics of these results is as follows.

Let $L \gg 1$. For all $n \in \mathbb{N}$, try to compare $V(n + L)$ to $V(n)$. If $n \gg L$, then the VRRW between these times behaves as though $V$ were constant and, hence, approximates a Markov chain which we call $M(V(n))$. Let $\pi(V(n))$ be the invariant measure of $M(V(n))$. If $L$ is assumed to be large enough, then the occupation measure between these times will be close to $\pi(V(n))$. This means that, approximately,

(1) $$(n + L)V(n + L) = nV(n) + L\pi(V(n));$$

hence

(2) $$V(n + L) - V(n) = (L/n)(\pi(V(n)) - V(n)).$$

Passing to a continuous time limit gives

(3) $$\frac{d}{dt}V(t) = \frac{1}{t}(\pi(V(t)) - V(t)).$$

Up to an exponential time change, $V$ should behave like an integral curve for the vector field $\pi - I$.

If the graph is not assumed to be complete, then the relaxation time of the Markov chain $M(V(n))$ depends on $V(n)$, (1) and (2) do not make sense in general, and it is difficult to deduce some results from this heuristics. In the critical case where this relaxation time is on the order of $n$, it may occur that the random walk gets stuck with high probability in a proper subset of the graph.

Our work relies on the principle that, on $\mathbb{Z}$, when this relaxation time is large, there are some seldom visited vertices between some often visited vertices. In this case, the behavior of the occupation densities of the random walk can be studied nearly independently to the left and to the right of each seldom visited vertex.

This notion of asymptotically seldom visited vertex $x \in \mathbb{Z}$ corresponds, in the following notation, to event $\Upsilon(x)$ defined below. For $x \in \mathbb{Z}$ and $n \in \mathbb{N}$, denote

$$Z_n^{\pm}(x) := \sum_{k=1}^{n} \mathbb{1}_{\{X_{k-1}=x, X_k=x\pm 1\}},$$

$$Y_n(x) := \sum_{k=1}^{n} \mathbb{1}_{\{X_{k-1}=x\}} \frac{1}{Z_{k-1}(x-1) + Z_{k-1}(x+1)},$$

$$Y_n^{\pm}(x) := \sum_{k=1}^{n} \mathbb{1}_{\{X_{k-1}=x, X_k=x\pm 1\}} \frac{1}{Z_{k-1}(x\pm 1)},$$

$$\alpha_n^{\pm}(x) := \frac{Z_n(x\pm 1)}{Z_n(x-1) + Z_n(x+1)}, \qquad \beta_n^{\pm}(x) := \frac{Z_n(x\pm 1)}{Z_n(x)},$$

$$\widetilde{Y}_n^{\pm}(x) := \sum_{k=1}^{n} \mathbb{1}_{\{X_{k-1}=x\pm 1, X_k=x\}} \frac{1}{Z_{k-1}(x\pm 1)},$$

$$\bar{Y}_n^{\pm}(x) := \sum_{k=1}^{n} \frac{\mathbb{1}_{\{X_k=x\}}}{Z_k(x)} \alpha_k^{\pm}(x).$$

Since, for any fixed $x \in \mathbb{Z}$, the sequences $Z_n(x)$, $Y_n(x)$, $Y_n^{\pm}(x)$, $\widetilde{Y}_n^{\pm}(x)$ and $\bar{Y}_n^{\pm}(x)$ are monotone nondecreasing in $n$, it makes sense to denote

$$Z_\infty(x) := \lim_{n\to\infty} Z_n(x), \qquad Z_\infty^{\pm}(x) := \lim_{n\to\infty} Z_n^{\pm}(x),$$

$$Y_\infty(x) := \lim_{n\to\infty} Y_n(x), \qquad Y_\infty^{\pm}(x) := \lim_{n\to\infty} Y_n^{\pm}(x),$$

$$\widetilde{Y}_\infty^{\pm}(x) := \lim_{n\to\infty} \widetilde{Y}_n^{\pm}(x), \qquad \bar{Y}_\infty^{\pm}(x) := \lim_{n\to\infty} \bar{Y}_n^{\pm}(x).$$

Let us define the probability events

$$\Upsilon(x) := \{Y_\infty(x) < \infty\},$$

$$\Upsilon^-(x) := \{Y_\infty^-(x) < \infty\}, \qquad \Upsilon^+(x) := \{Y_\infty^+(x) < \infty\}.$$

Let us enumerate a few properties about these events $\Upsilon(x)$, $x \in \mathbb{Z}$. First, for all $x \in \mathbb{Z}$, $\Upsilon(x)$ coincides a.s. with the set $\Upsilon^+(x)$ on which there are a small number of visits from $x$ to $x+1$ and, by symmetry, with the set $\Upsilon^-(x)$. This property is stated in the following lemma.

LEMMA 2.1.   *For all $x \in \mathbb{Z}$,*

$$\Upsilon(x) = \Upsilon^+(x) = \Upsilon^-(x).$$

PROOF.    Using the conditional Borel–Cantelli lemma [Lemma A.1(i)],

$$\begin{aligned}
\Upsilon^\pm(x) &= \left\{ \sum_{k=1}^\infty \mathbb{1}_{\{X_{k-1}=x, X_k=x\pm1\}} \frac{1}{Z_{k-1}(x\pm1)} < \infty \right\} \\
&= \left\{ \sum_{k=1}^\infty \mathbb{E}\left[ \frac{\mathbb{1}_{\{X_{k-1}=x, X_k=x\pm1\}}}{Z_{k-1}(x\pm1)} \Big| \mathcal{F}_{k-1} \right] < \infty \right\} \\
&= \left\{ \sum_{k=1}^\infty \mathbb{1}_{\{X_{k-1}=x\}} \frac{1}{Z_{k-1}(x-1)+Z_{k-1}(x+1)} < \infty \right\} = \Upsilon(x). \quad \square
\end{aligned}$$

Second, there are at most two consecutive infinitely often visited sites $x \in \mathbb{Z}$ on which $\Upsilon(x)$ holds, as implied by Lemma 2.2.

LEMMA 2.2.   *For all $x \in \mathbb{Z}$,*

$$\Upsilon(x-1) \cap \Upsilon(x+1) = \{Z_\infty(x) < \infty\}.$$

PROOF.    Indeed, by Lemma 2.1,

$$\begin{aligned}
\Upsilon(x-1) &\cap \Upsilon(x+1) \\
&= \Upsilon^+(x-1) \cap \Upsilon^-(x+1) \\
&\subset \left\{ \ln Z_\infty(x) \approx \sum_{k=1}^\infty \frac{\mathbb{1}_{\{X_k=x\}}}{Z_{k-1}(x)} = Y_\infty^+(x-1) + Y_\infty^-(x-1) < \infty \right\} \\
&= \{Z_\infty(x) < \infty\}.
\end{aligned}$$

The reverse inclusion is straightforward.   $\square$

Third, if $\Upsilon(x)$ holds, we can give some information on the behavior of $\alpha_n^-(x+2)$ (since there are a small number of visits from $x$ to $x+1$) as stated in Corollary 3.1(ii). Note that the entire Corollary 3.1 is stated (and proved) in Section 3.1.

COROLLARY 3.1(ii).  *For all $x \in \mathbb{Z}$,*

$$\Upsilon(x) \subset \left\{ \exists \alpha_\infty^\mp(x \pm 2) := \lim_{n \to \infty} \alpha_n^\mp(x \pm 2) \in [0, 1) \right\}.$$

Fourth, we can claim a kind of propagation rule on seldom visited sites as given by the following proposition.

PROPOSITION 2.1.  *For all $x \in \mathbb{Z}$,*

$$\Upsilon(x) \subset \Upsilon(x + 1) \cup \Upsilon(x + 4).$$

This result is closely related to the dynamics inherent to the random walk. We cannot directly use the methods of comparison with the dynamical system, since there is no tool that gives a control on the behavior of the random walk on more than a few vertices.

The heuristic of the proposition is that there is a kind of competition between the numbers of visits to points $x + 1$ and $x + 4$. Its proof is divided into two cases. If $\alpha_\infty^-(x + 2)$ is positive, then $x + 4$ loses, which implies that $\Upsilon(x + 4)$ holds. On the other hand, if $\alpha_\infty^-(x + 2)$ is equal to zero, then $x + 1$ loses, which implies that $\Upsilon(x + 1)$ holds. These two results are implied, respectively, by Corollary 3.1(iii) and Lemma 2.3:

COROLLARY 3.1(iii).  *For all $x \in \mathbb{Z}$,*

$$\Upsilon(x) \cap \{\alpha_\infty^-(x + 2) > 0\} \subset \Upsilon(x + 4).$$

LEMMA 2.3.  *For all $x \in \mathbb{Z}$,*

$$\Upsilon(x) \cap \{\alpha_\infty^-(x + 2) = 0\} \subset \Upsilon(x + 1).$$

Recall that Corollary 3.1 is proved in Section 3.1. Lemma 2.3 is equivalent to the statement that the random set

$$\Upsilon_0(x) = \Upsilon(x) \cap \{\alpha_\infty^-(x + 2) = 0\} \cap \Upsilon(x + 1)^c$$

is of probability 0. Before proving this lemma, we prove in Lemma 2.4 (proved in Section 3.3) that $\Upsilon_0(x)$ is a.s. a subset of $\Upsilon_0'(x)$ (defined hereafter), on which we have a rough control on the behavior of the random walk on sites $x$ to $x + 5$. Then Lemma 2.5 (stated hereafter and proved in Section 5.1) completes the proof.

Let $e := \exp(1)$. Let, for all $x \in \mathbb{Z}$,

$$\Upsilon_0'(x) = \left\{ \limsup \frac{Z_n(x + 4)}{Z_n(x + 1)} \le e \right\} \cap \left\{ \limsup \left( \sup_{k \ge n} \frac{\alpha_k^-(x + 2)}{\alpha_n^-(x + 2)} \right) \le 1 \right\}$$

$$\cap \left\{ \lim \frac{\ln Z_n(x+1)}{\ln Z_n(x+2)} = \lim \frac{\ln Z_n(x+4)}{\ln Z_n(x+2)} = 1 \right\}$$

$$\cap \left\{ Z_\infty(x) = Z_\infty(x+4) = \infty \right\}$$

$$\cap \left\{ \lim \frac{Z_n(x+3)}{Z_n(x+2)} = 1 \right\}$$

$$\cap \left\{ \limsup \frac{Z_n(x+5)}{Z_n(x+3)} \le 1 \right\} \cap \left\{ \limsup \frac{Z_n(x)}{Z_n(x+2)} \le 1 \right\}.$$

LEMMA 2.4. *For all $x \in \mathbb{Z}$, $\Upsilon_0(x) \subset \Upsilon'_0(x)$.*

LEMMA 2.5. *For all $x \in \mathbb{Z}$, $\mathbb{P}(\Upsilon_0(x) \cap \Upsilon'_0(x)) = 0$.*

The case of $\Upsilon_0(x) \cap \Upsilon'_0(x)$ considered in Lemma 2.5 corresponds to an unstable set in the dynamical systems setting. To prove the nonconvergence to this set without a complete control on the behavior of the empirical density of occupation, we use a partial order on a certain class of random walks on $\mathbb{Z}$ and prove an appropriate result in some unstable situations (Section 4).

Let us now go back to the description of seldom visited sites. Theorem 1.1 implies that there exists a.s. a leftmost infinitely visited site $x_0$. By definition, $Z_\infty(x_0 - 1) < \infty$, which implies that $\Upsilon(x_0 - 1)$ and $\Upsilon(x_0) = \Upsilon^-(x_0)$ hold (using Lemma 2.1). Accordingly, Proposition 2.1 and Lemma 2.2 lead us to a pavement of the set of infinitely often visited vertices (which is connected) by sites on which $\Upsilon(x)$ holds.

More precisely, let us denote, for any finite sequence $(x_i)_{1 \le i \le n}$ taking values in $\mathbb{Z}$, the event

$$\Upsilon((x_i)_{1 \le i \le n}) = \bigcap_{1 \le i \le n} \Upsilon(x_i).$$

Let us define the events

$$\Omega(x) = \{x = \inf R'\},$$

$$\Omega_0(x) = \Omega(x) \cap \{Z_\infty(x+5) < \infty\},$$

$$\Omega_1(x) = \Upsilon(x, x+4, x+8) \cap \{Z_\infty(x+1) = Z_\infty(x+7) = \infty\},$$

$$\Omega_2(x) = \Upsilon(x-1, x, x+4, x+5, x+9, x+10)$$
$$\cap \{Z_\infty(x+1) = Z_\infty(x+8) = \infty\}.$$

We can state the following lemma.

LEMMA 2.6. *For all $x \in \mathbb{Z}$,*

$$\Omega(x) \subset \Omega_0(x) \cup \Omega_1(x) \cup \Omega_1(x+5) \cup \Omega_2(x).$$

PROOF. First, for all $y \in \mathbb{Z}$,

$$(4) \qquad \Upsilon(y-1, y) \cap \{Z_\infty(y) = \infty\} \subset \Upsilon(y+4),$$

since, by Proposition 2.1,

$$\Upsilon(y) \subset \Upsilon(y+1) \cup \Upsilon(y+4)$$

and, by Lemma 2.2,

$$\Upsilon(y-1) \cap \Upsilon(y+1) \subset \{Z_\infty(y) < \infty\}.$$

This implies

$$\Omega(x) \cap \Omega_0(x)^c$$
$$\subset \Upsilon(x-1, x) \cap \{Z_\infty(x) = Z_\infty(x+5) = \infty\}$$
$$\subset \Upsilon(x-1, x, x+4) \cap \{Z_\infty(x) = Z_\infty(x+5) = \infty\}$$
$$\subset (\Upsilon(x-1, x, x+4, x+5) \cup \Upsilon(x, x+4, x+8))$$
$$\cap \{Z_\infty(x) = Z_\infty(x+7) = \infty\},$$

where we use (4) with $y := x$ in the second inclusion and use Proposition 2.1 with $x := x + 4$ in the third inclusion $[Z_\infty(x+7) = \infty$ follows from the convergence of $\alpha_n^-(x+6)$ on $\Upsilon(x+4)$, by Corollary 3.1(ii), together with $Z_\infty(x+5) = \infty]$.

Now

$$\Upsilon(x-1, x, x+4, x+5) \cap \{Z_\infty(x+5) = \infty\}$$
$$\subset \Upsilon(x-1, x, x+4, x+5, x+9)$$
$$\subset \Upsilon(x-1, x, x+4, x+5, x+9, x+10)$$
$$\cup \Upsilon(x-1, x, x+4, x+5, x+9, x+13),$$

where we use (4) with $y := x + 5$ in the first inclusion and use Proposition 2.1 with $x := x + 9$ in the second inclusion.

Putting together these two equations, we obtain

$$\Omega(x) \cap \Omega_0(x)^c \cap \Omega_1(x)^c$$
$$\subset \Upsilon(x-1, x, x+4) \cap \{Z_\infty(x+1) = Z_\infty(x+7) = \infty\}$$
$$\cap (\Upsilon(x+5, x+9, x+10) \cup \Upsilon(x+5, x+9, x+13))$$
$$\subset \Omega_2(x) \cup \Omega_1(x+5),$$

where we note in the second inclusion that $Z_\infty(x+8) = \infty$ if $\Upsilon(x+8)$ does not hold and, similarly, $Z_\infty(x+10) = Z_\infty(x+12) = \infty$ [since $\alpha_\infty^-(x+11) \in [0,1)$] on $\Upsilon(x+9)$ if $\Upsilon(x+10)$ does not hold. $\square$

Now, for all $x \in \mathbb{Z}$, $\Omega_1(x)$ and $\Omega_2(x)$ are of probability 0, as stated in Lemmas 2.7 and 2.8. These results complete the proof of the conjecture.

LEMMA 2.7.  *For all $x \in \mathbb{Z}$, $\mathbb{P}(\Omega_1(x)) = 0$.*

LEMMA 2.8.  *For all $x \in \mathbb{Z}$, $\mathbb{P}(\Omega_2(x)) = 0$.*

Let us explain in a few words the proofs of these lemmas. Lemma 2.7 relies on the fact that there is a kind of competition between the numbers of visits to points $x+1$, $x+2$ and $x+3$ on the left-hand side of $x+4$, and $x+5$, $x+6$ and $x+7$ on the right-hand side of $x+4$. We first prove the following lemma.

LEMMA 2.9.  *For all $x \in \mathbb{Z}$, $\Omega_1(x) \subset \{\lim Z_n(x+6)/Z_n(x+2) = 1\}$.*

The heuristic of Lemma 2.9 is that if $Z_n(x+6)/Z_n(x+2)$ did not converge to 1, then it would converge to 0 or to $\infty$, and that these convergences would be so fast that $Z_\infty(x+6) < \infty$ in the first case and $Z_\infty(x+2) < \infty$ in the second case.

The proof of Lemma 2.7 therefore reduces to the study of the unstable case $Z_n(x+6)/Z_n(x+2) \to 1$. The methods used for this proof in Section 5.2 rely, similarly as in the proof of Lemma 2.5, on the tools introduced in Section 4.

The proof of Lemma 2.8 has roughly the same heuristic as Lemma 2.7, but we have to face the problem explained at the beginning of this section, that is, we have to discriminate between the case $\{\alpha_\infty^-(x+7) > \alpha_\infty^-(x+2)\}$, where the random walk regularly visits the set $\{x, \ldots, x+9\}$, and the case $\{\alpha_\infty^-(x+7) \leq \alpha_\infty^-(x+2)\}$, where the random walk eventually gets stuck in a strict subset [i.e., $Z_\infty(x+4) < \infty$ or $Z_\infty(x+5) < \infty$]. This study corresponds to Lemma 2.10, stated subsequently and proved in Section 3.5.

LEMMA 2.10.  *For all $x \in \mathbb{Z}$,*
$$\Omega_2(x) \subset \{\alpha_\infty^-(x+7) > \alpha_\infty^-(x+2)\}.$$

Next, we prove in Lemma 2.11 that $Z_n(x+7)/Z_n(x+2) \to 1$ on $\Omega_2(x)$ and we finish the proof of Lemma 2.8 in Section 5.3, using again the methods introduced in Section 4.

LEMMA 2.11.  *For all $x \in \mathbb{Z}$, $\Omega_2(x) \subset \{\lim Z_n(x+7)/Z_n(x+2) = 1\}$.*

2.3. *Outline of contents.*  Section 3 gives some preliminary results, based on martingales techniques. This section is divided into six parts. In Section 3.1, we prove results related to the Pólya and Friedman urn models. In Section 3.2 we give a comparison tool (Lemma 3.1) and prove an estimate that gives conditions for a site to be finitely often visited (Lemma 3.2). Finally, we deduce Lemmas 2.4, 2.9, 2.10 and 2.11, respectively, in Sections 3.3, 3.4, 3.5 and 3.6. In Section 4, we prove a result of nonconvergence in unstable situations, using a partial order on a certain class of random walks on $\mathbb{Z}$. This result is useful to the proofs of Lemmas 2.5, 2.7 and 2.8. In Section 5, we apply this result to the proofs of Lemmas 2.5, 2.7 and 2.8 (resp. in Sections 5.1, 5.2 and 5.3). In the Appendix, we state some general martingale results and, in particular, recall a generalized version of the conditional Borel–Cantelli lemma.

## 3. Preliminary results.

3.1. *Martingale results.* The following Proposition 3.1 and its Corollaries 3.1 and 3.2 provide us with some local properties of the VRRW, relating the quantities $Y_n(x)$, $Y_n^\pm(x)$, $\widetilde{Y}_n^\pm(x)$ and $\overline{Y}_n^\pm(x)$ defined in Section 2.2. These results enable us to describe the behavior of the random walk on the first few points following $x$ when one event like $\Upsilon_0(x)$, $\Omega_1(x)$ or $\Omega_2(x)$ holds.

Let us focus on the two key properties of this part, namely Proposition 3.1(a) and Corollary 3.2(i), which are related to the Pólya and Friedman urn models. Note that a detailed survey on the relationships between these urn models and reinforcement processes can be found in [9].

Proposition 3.1(a) can be explained by studying first the case of a VRRW on three consecutive points $\{x-1, x, x+1\}$. Under this assumption, the walk is half of the time in sites $x-1$ or $x+1$, and comes back to $x$ at the next step; the other half of the time, the walk is in site $x$ and moves to $x \pm 1$ with a probability equal to the number of times $x \pm 1$ has been visited up through time $n$ divided by the total number of visits to $x-1$ and $x+1$ [with the convention that the sites $x \pm 1$ have been visited $Z_0(x \pm 1)$ at time 0].

This construction is equivalent to a Pólya urn model with two colors $x-1$ and $x+1$, with $Z_n(x-1)$ and $Z_n(x+1)$ balls of colors $x-1$ and $x+1$ at time $n$. Indeed, this corresponds to the process of picking, half of the time, a ball at random in the urn and replacing it with a ball of the same color.

A classical result claims that the proportion of balls of color $x-1$ converges toward a random $\alpha \in (0,1)$. The random variable $\alpha$ has a beta distribution of parameters $Z_0(x-1)$ and $Z_0(x+1)$ (see, e.g., [5], Vol. 2, Chapter VII), but this result is difficult to use in our context, where, in the general case, we have to deal with visits from $x+2$ to $x+1$ and from $x-2$ to $x-1$. Observe that this convergence can be proved by Proposition 3.1(a). Indeed, it implies that $Y_n^+(x) - Y_n^-(x)$ converges, and we deduce from the convergence of $Y_n^\pm(x) - \ln Z_n(x \pm 1)$ (approximation of log by the harmonic series) that $\ln(Z_n(x+1)/Z_n(x-1))$ converges.

Let us now return to the study of the VRRW on $\mathbb{Z}$. The equivalence between the weighted numbers of visits from $x$ to $x+1$ and from $x$ to $x-1$, $Y_n^+(x)$ and $Y_n^-(x)$, claimed in Proposition 3.1(a), enables us to estimate in Corollary 3.1(i) the variation of $\ln Z_n(x+1)/Z_n(x-1)$, with respect to $Y_n^-(x+2)$ and $Y_n^+(x-2)$. Corollary 3.1(ii)–(iv) is a direct consequence of this claim.

Let us now explain the heuristic of Corollary 3.2(i). Let us consider the case where, given $x \in \mathbb{Z}$, the event

$$R' = \{x-2, x-1, x, x+1, x+2\}$$

holds. Then $\Upsilon^-(x-2) = \Upsilon(x-2)$ and $\Upsilon(x+2)$ hold. This implies by Corollary 3.1(ii) that $\alpha_n^-(x)$ converges to $\alpha_\infty^-(x) \in (0,1)$ (see also Remark 3.1). Let us study the behavior of the random walk on the border point $x-2$.

Let us denote by $t_n$ the $n$th visit time to site $x-1$. We again observe an urn model with two colors $x-2$ and $x$, $Z_{t_n}(x-2)$ and $Z_{t_n}(x)$ being the numbers of balls of color $x-2$ and $x$ at the $n$th iteration. Indeed, at time $t_n$, we move to $x-2$ with probability $\alpha_{t_n}^-(x-1)$; this operation is equivalent to picking a ball at random in the urn, and similarly for $x$. Now, if we move to $x-2$, we come back to $x-1$ at the next step (unless we move to $x-3$, which occurs only finitely often). On the other hand, if we move to $x$, the expected number of visits to $x$ before returning to $x-1$ is on the order of $1/\alpha_\infty^-(x)$ [expectation of a geometric random variable with success probability $\alpha_\infty^-(x)$].

In the urn model, this means that if we pick a ball of color $x-2$, we replace it with a ball of the same color, and that if we pick a ball of color $x$, we replace it, on average, by $1/\alpha_\infty^-(x)$ balls of color $x$. This corresponds to a generalized Pólya–Friedman urn model (see, e.g., [9]), and a classical result claims that this implies $Z_{t_n}(x) \approx Z_{t_n}(x-2)^{1/\alpha_\infty^-(x)}$. We obtain the same conclusion by Corollary 3.2(i) and (iv):

$$\ln Z_n(x-2) \equiv \bar{Y}_n^-(x) \equiv \sum_{k=1}^n \frac{\mathbb{1}_{\{X_k=x\}}}{Z_k(x)} \alpha_k^-(x) \equiv \alpha_\infty^-(x) \ln Z_n(x).$$

This gives an intuition for this Corollary 3.2, which is needed here instead of generalized Pólya–Friedman urn results since the corresponding martingale technique is more adaptable to the case of visits from $x-3$ to $x-2$.

Note that these methods provide the asymptotic behavior of the VRRW on $\mathbb{Z}$, conditional on the event that we eventually get stuck on five points. Indeed, we obtain that events 2–6 (defined in Section 1) hold, conditional on to event 1. Another proof of this result was given by Bienvenüe in his Ph.D. dissertation [3], using ideas related to the construction of continuous reinforced random walks (see [10]).

PROPOSITION 3.1. *For all $x \in \mathbb{Z}$ and $\nu < 1/2$,*

(a) $Y_n^\pm(x) \doteq Y_n(x) + o(Z_n(x \pm 1)^{-\nu})$,
(b) $\ln Z_n(x) \doteq Y_n^+(x-1) + Y_n^-(x+1) + O(Z_n(x)^{-1})$,

(c) $Y_n^\pm(x) \doteq \widetilde{Y}_n^\pm(x) - \mathbb{1}_{\{\pm X_n \leq \pm x\}}/(Z_{n-1}(x \pm 1)) \doteq \bar{Y}_n^\mp(x \pm 1) + o(Z_n(x \pm 1)^{-\nu})$.

PROOF. Let us first prove statement (a). Given $\nu_0 < 1/2$, we apply Lemma A.1(iii) with

$$\Gamma_k = \{X_{k-1} = x, X_k = x \pm 1\}, \qquad \xi_k = 1/Z_k(x \pm 1), \qquad \beta_k = Z_k(x \pm 1)^{2\nu_0},$$

to conclude that

$$(5) \qquad\qquad Y_n^\pm(x) - Y_n(x) \doteq O(Z_n(x \pm 1)^{-\nu_0}).$$

Indeed, using the notation of this lemma,

$$\sum_{k=0}^\infty \beta_k \delta_k = \sum_{k=0}^\infty \alpha_k^-(x) \alpha_k^+(x) \frac{\mathbb{1}_{\{X_k = x\}}}{Z_k(x \pm 1)^{2(1-\nu_0)}}$$

$$\leq \sum_{k=0}^\infty \alpha_k^\pm(x) \frac{\mathbb{1}_{\{X_k = x\}}}{Z_k(x \pm 1)^{2(1-\nu_0)}} \approx \sum_{k=0}^\infty \frac{\mathbb{1}_{\{X_k = x, X_{k+1} = x \pm 1\}}}{Z_k(x \pm 1)^{2(1-\nu_0)}} < \infty,$$

the last equivalence being a consequence of the conditional Borel–Cantelli lemma, Lemma A.1(i). Therefore, the conditions of Lemma A.1(iii) are statisfied and (5) holds.

Statement (a) follows directly if $Z_\infty(x \pm 1) = \infty$, by choosing $\nu_0 > \nu$. Otherwise, by Lemma A.1(i), $Y_\infty(x) \asymp Y_\infty^\pm(x) < \infty$, which also enables us to conclude (a).

Statement (b) follows from

$$Y_n^+(x-1) + Y_n^-(x+1) = \sum_{k=1}^n \frac{\mathbb{1}_{\{X_k = x\}}}{Z_{k-1}(x)} = \sum_{j=1+\mathbb{1}_{\{X_0 = x\}}}^{Z_n(x)-1} \frac{1}{j}$$

$$= \ln Z_n(x) + \mathrm{Cst}(x, v_0) + \square(Z_n(x)^{-1})$$

when $Z_n(x) \geq \mathrm{Cst}$.

Let us now prove statement (c) for $Y_n^+(x)$; the proof for $Y_n^-(x)$ is similar. For all $n \in \mathbb{N}^*$, let $u_n$ (resp. $v_n$) be the time of the $n$th visit from $x$ to $x+1$ (resp. from $x+1$ to $x$), that is,

$$u_n = \inf\{k \in \mathbb{N}^*/Z_k^+(x) = n\}, \qquad v_n = \inf\{k \in \mathbb{N}^*/Z_k^-(x+1) = n\}.$$

Recall that the number of visits $Z_k^\pm(x)$ from $x$ to $x \pm 1$ at time $k$ is defined in Section 2.1.

Assume, for instance, $v_1 < u_1$ (the other case is similar). Then, for all $n \in \mathbb{N}^*$, $u_n < v_{n+1} < u_{n+1}$ (with the convention that $\infty < \infty$) and $Z_{v_n-1}(x +$

$1) = Z_{u_{n}-1}(x+1)$. Therefore,

$$Y_{u_n}^+(x) - Y_{u_1}^+(x) = \sum_{k=u_1+1}^{u_n} \mathbb{1}_{\{X_{k-1}=x, X_k=x+1\}} \frac{1}{Z_{k-1}(x+1)}$$

$$= \sum_{j=2}^{n} \frac{1}{Z_{u_j-1}(x+1)}$$

$$= \sum_{j=2}^{n} \frac{1}{Z_{v_j-1}(x+1)} = \sum_{k=v_2}^{v_n} \mathbb{1}_{\{X_{k-1}=x+1, X_k=x\}} \frac{1}{Z_{k-1}(x+1)}$$

$$= \sum_{k=u_1+1}^{u_n} \mathbb{1}_{\{X_{k-1}=x+1, X_k=x\}} \frac{1}{Z_{k-1}(x+1)} = \widetilde{Y}_{u_n}^+(x) - \widetilde{Y}_{u_1}^+(x).$$

This gives the first equivalence of (c) when $u_n < \infty$ for all $n \in \mathbb{N}^*$. Otherwise, $Y_k^+(x) - \widetilde{Y}_k^+(x)$ is constant for large enough $k \in \mathbb{N}$, which also gives the equivalence.

Let us prove the second equivalence of (c). If $Z_\infty(x \pm 1) < \infty$, then $\bar{Y}_\infty^\pm(x) \asymp \widetilde{Y}_\infty^\pm(x) < \infty$ by Lemma A.1(i), which enables us to conclude. Otherwise, apply Lemma A.1(iii) and use its notation, with

$$\Gamma_k = \{X_{k-1} = x \pm 1, X_k = x\}, \qquad \xi_k = 1/Z_k(x \pm 1), \qquad \beta_k = Z_k(x \pm 1)^{2\nu},$$

where we note that

$$\sum_{k=0}^{\infty} \beta_k \delta_k \leq \sum_{k=0}^{\infty} \frac{\mathbb{1}_{\{X_k=x \pm 1\}}}{Z_k(x \pm 1)^{2(1-\nu)}} < \infty. \qquad \square$$

COROLLARY 3.1. *For all* $x \in \mathbb{Z}$ *and* $\nu < 1/2$:

(i) $\ln(Z_n(x-1))/(Z_n(x+1)) \doteq Y_n^+(x-2) - Y_n^-(x+2) + o(Z_n(x-1)^{-\nu}) + o(Z_n(x+1)^{-\nu})$;

(ii) $\Upsilon(x) \subset \{\exists \alpha_\infty^\mp(x \pm 2) := \lim_{n \to \infty} \alpha_n^\mp(x \pm 2) \in [0,1)\}$;

(iii) $\Upsilon(x) \cap \{\alpha_\infty^\mp(x \pm 2) > 0\} \subset \Upsilon(x \pm 4)$;

(iv) $\Upsilon(x) \subset \{\exists \beta_\infty^\mp(x \pm 2) := \lim_{n \to \infty} \beta_n^\mp(x \pm 2) \in [0,\infty)\}$;

(v) $\Upsilon(x) \cap \{Z_\infty(x \pm 2) = \infty\} \subset \{\beta_\infty^\mp(x \pm 2) = \alpha_\infty^\mp(x \pm 2)\}$.

REMARK 3.1. Corollary 3.1 implies that a.s. on $\Upsilon(x, x+4)$, $\alpha_n^\pm(x+2)$ [resp. $\beta_n^\pm(x+2)$] converges to $\alpha_\infty^\pm(x+2) \in (0,1)$ [resp. $\beta_\infty^\pm(x+2) > 0$] and that $\beta_\infty^\pm(x+2) = \alpha_\infty^\pm(x+2)$ if, moreover, $Z_\infty(x+2) = \infty$. This follows from an application of statements (ii) and (iv)–(v) successively to $x$ with $-$ instead of $\mp$ and to $x+4$ with $+$ instead of $\mp$, and from $\alpha_\infty^+(x+2) + \alpha_\infty^-(x+2) = 1$.

Proof of Corollary 3.1. It follows from statement (b) of Proposition 3.1, applied to $x-1$ and $x+1$ that, for all $\nu < 1/2$,

$$Y_n^-(x) \doteq \ln Z_n(x-1) - Y_n^+(x-2) + o(Z_n(x-1)^{-\nu}),$$

$$Y_n^+(x) \doteq \ln Z_n(x+1) - Y_n^-(x+2) + o(Z_n(x+1)^{-\nu}).$$

These equivalences remain true in the cases $Z_\infty(x-1) < \infty$ and $Z_\infty(x+1) < \infty$.

It also follows from statement (a) of the proposition that, for all $\nu < 1/2$,

$$Y_n^-(x) \doteq Y_n^+(x) + o(Z_n(x-1)^{-\nu}) + o(Z_n(x+1)^{-\nu}),$$

which completes the proof of (i). Let us now prove (ii) and (iii) for $\alpha_n^-(x+2)$; the case of $\alpha_n^+(x+2)$ is similar. Apply (i) for $x+2$: On $\Upsilon(x)$,

$$\ln \frac{Z_n(x+1)}{Z_n(x+3)} \equiv Y_n^+(x) - Y_n^-(x+4) \equiv -Y_n^-(x+4)$$

and $Y_n^-(x+4)$ is nondecreasing in $n$, which completes the proof.

Let us now prove (iv) and (v) for $\beta_n^-(x+2)$; the proof for $\beta_n^+(x-2)$ is similar. Assume that $\Upsilon(x)$ holds: By statements (a) (applied to $x+2$) and (b) (applied to $x+1$) of the proposition,

$$\ln Z_n(x+1) \equiv Y_n^+(x) + Y_n^-(x+2)$$

$$\equiv Y_n^-(x+2) \equiv Y_n(x+2) = \sum_{k=0}^{n-1} \frac{\mathbb{1}_{\{X_k=x+2\}}}{Z_k(x+2)} \frac{\alpha_k^-(x+2)}{\beta_k^-(x+2)},$$

and, therefore,

$$(6) \qquad \ln \beta_n^-(x+2) = \ln \frac{Z_n(x+1)}{Z_n(x+2)} \equiv \sum_{k=0}^{n-1} \frac{\mathbb{1}_{\{X_k=x+2\}}}{Z_k(x+2)} \left( \frac{\alpha_k^-(x+2)}{\beta_k^-(x+2)} - 1 \right).$$

Now, $Z_k(x+2) \le Z_k(x+1) + Z_k(x+3)$ implies $\alpha_k^-(x+2) \le \beta_k^-(x+2)$. Hence the right-hand side of the equation is nonincreasing in $n$, which implies (iv). Let us further assume that $Z_\infty(x+2) = \infty$. If $\beta_\infty^-(x+2) > 0$ and $\beta_\infty^-(x+2) \ne \alpha_\infty^-(x+2)$, then $\alpha_k^-(x+2)/\beta_k^-(x+2) - 1$ converges to a negative real and (6) implies $\ln \beta_\infty^-(x+2) = -\infty$, so that $\beta_\infty^-(x+2) = 0$, which leads to a contradiction. If $\beta_\infty^-(x+2) = 0$, then $\alpha_\infty^-(x+2) \le \beta_\infty^-(x+2) = 0$. This completes the proof of (v). □

Corollary 3.2. *For all $x \in \mathbb{Z}$, $\gamma \in (0,1)$ and $\nu < 1/2$:*

(i) $\Upsilon(x-1) \subset \{\ln Z_n(x) \equiv Y_n^-(x+1) \doteq \bar{Y}_n^-(x+2) + o(Z_n(x)^{-\nu})\}$;

(ii) $\Upsilon(x-1) \cap \{\limsup \alpha_n^-(x+2) \le \gamma\} \subset \Upsilon(x-1) \cap \{\ln Z_n(x) \preceq \gamma \ln Z_n(x+2)\} \subset \Upsilon(x-1,x)$;

(iii) $\{\liminf \alpha_n^-(x+2) \geq \gamma\} \subset \{\ln Z_n(x) \succeq \gamma \ln Z_n(x+2)\}$;

(iv) $\Upsilon(x-1,x,x+4,x+5) \subset \{\exists \delta > 0 / \alpha_n^-(x+2) - \alpha_\infty^-(x+2) = o(Z_n(x+2)^{-\delta})\} \cap \{\ln Z_n(x) \equiv \alpha_\infty^-(x+2) \ln Z_n(x+2) \equiv \alpha_n^-(x+2) \ln Z_n(x+2)\}$.

PROOF. Let us first prove (i). Assume that $\Upsilon(x-1)$ holds and apply Proposition 3.1(a), (b) and (c): For all $\nu < 1/2$,

$$
\begin{aligned}
\ln Z_n(x) &\equiv Y_n^-(x+1) \doteq Y_n^+(x+1) + o(Z_n(x)^{-\nu}) + o(Z_n(x+2)^{-\nu}) \\
&\doteq \bar{Y}_n^-(x+2) + o(Z_n(x)^{-\nu}) + o(Z_n(x+2)^{-\nu}) \\
&\doteq \bar{Y}_n^-(x+2) + o(Z_n(x)^{-\nu}),
\end{aligned}
$$

(7)

where we use in the last equation that $\alpha_n^-(x+1)$ converges to $\alpha_\infty^-(x+1) \in [0,1)$ by Corollary 3.1(ii) and, therefore, that $o(Z_n(x+2)^{-\nu})$ is upper bounded by $o(Z_n(x)^{-\nu})$. The first inclusion of statement (ii) and statement (iii) follow directly.

Let us prove the second part of (ii). Assume $\Upsilon(x-1) \cap \{\ln Z_n(x) \preceq \gamma \ln Z_n(x+2)\}$ holds. We prove the stronger statement that there exists $\delta > 0$ such that (8) holds, which completes the proof by Lemma 2.1 [$\Upsilon(x) = \Upsilon^+(x)$] and also is useful in the proof of (iv). If $Z_\infty(x+1) < \infty$, then $Y_\infty^+(x) < \infty$, which proves the statement. Otherwise, for all $\varepsilon > 0$, for sufficiently large $k \in \mathbb{N}$, using $Z_k(x+1) \leq Z_k(x) + Z_k(x+2)$,

$$
\begin{aligned}
\alpha_k^-(x+1) &= \frac{Z_k(x)}{Z_k(x) + Z_k(x+2)} \\
&= o((Z_k(x) + Z_k(x+2))^{\gamma+\varepsilon-1}) = o(Z_k(x+1)^{\gamma+\varepsilon-1}),
\end{aligned}
$$

which implies, for all $\delta < \min(1/2, 1-\gamma)$ and $\varepsilon < 1 - \gamma - \delta$, using Proposition 3.1(c),

$$
\begin{aligned}
Y_n^+(x) &\doteq \bar{Y}_n^-(x+1) + o(Z_n(x+1)^{-\delta}) \\
&= \sum_{k=n}^\infty \frac{\mathbb{1}_{\{X_k=x+1\}}}{Z_k(x+1)} \alpha_k^-(x+1) + o(Z_n(x+1)^{-\delta}) \\
&\doteq o(Z_n(x+1)^{\gamma+\varepsilon-1}) + o(Z_n(x+1)^{-\delta}) \doteq o(Z_n(x+1)^{-\delta}).
\end{aligned}
$$

(8)

Note that, conversely, (8) always holds on $\Upsilon(x-1,x)$ for all $\delta < \alpha_\infty^+(x+2) \wedge 1/2$, since $\limsup \alpha_n^-(x+2) = \alpha_\infty^-(x+2) < 1$ on $\Upsilon(x)$ by Corollary 3.1(ii).

Let us now assume that $\Upsilon(x-1,x,x+4,x+5)$ holds and prove (iv). First observe that both $\alpha_\infty^-(x+2)$ and $\alpha_\infty^+(x+2)$ are strictly positive, using Remark 3.1. For all $\delta < \alpha_\infty^+(x+2) \wedge 1/2$, $Y_n^+(x) \doteq o(Z_n(x+1)^{-\delta})$ by (8), and for all $\delta < \alpha_\infty^-(x+2) \wedge 1/2$ symmetrically (with respect to $x+2$), $Y_n^-(x+4) \doteq o(Z_n(x+3)^{-\delta})$.

Accordingly, using Corollary 3.1(i) with $x := x + 2$, there exists $\delta > 0$ such that

$$\ln \frac{Z_n(x+1)}{Z_n(x+3)} \doteq Y_n^+(x) - Y_n^-(x+4) + o(Z_n(x+1)^{-\delta}) + o(Z_n(x+3)^{-\delta})$$

(9)
$$\doteq o(Z_n(x+1)^{-\delta}) + o(Z_n(x+3)^{-\delta}) = o(Z_n(x+2)^{-\delta}),$$

where the last equality uses the observation that $o(Z_n(x+1)^{-\delta})$ and $o(Z_n(x+3)^{-\delta})$ are upper bounded by $o(Z_n(x+2)^{-\delta})$, since $\beta_\infty^\pm(x+2) = \alpha_\infty^\pm(x+2) > 0$ by Remark 3.1.

Equation (9) implies

$$\alpha_n^-(x+2) - \alpha_\infty^-(x+2) = o(Z_n(x+2)^{-\delta}).$$

Using (7), we deduce that

$$\ln Z_n(x) \equiv \bar{Y}_n^-(x+2) = \sum_{k=1}^n \frac{\mathbb{1}_{\{X_k = x+2\}}}{Z_k(x+2)} \alpha_k^-(x+2) \equiv \alpha_\infty^-(x+2) \ln Z_n(x+2).$$

$\square$

3.2. *Comparison results.* The following lemma considers the case of a sequence $u_n$ repelled by $a$ on its right-hand side, where the repulsion depends on a function $f$ of $u_n$ and on another sequence $v_n$. It yields, when $u_n$ does not asymptotically remain in $(-\infty, a]$, an estimate of $u_n$ as $n$ goes off to infinity.

LEMMA 3.1. *Let $f : \mathbb{R} \to \mathbb{R}$ be a nondecreasing function, positive on $(a, \infty)$, and let $(u_n)_{n \in \mathbb{N}}$ and $(v_n)_{n \in \mathbb{N}}$ be sequences that take values, respectively, in $\mathbb{R}$ and $\mathbb{R}^+$. Then*

$$\left\{ u_n \underset{u_n \geq a}{\succeq} \sum_{k=0}^{n-1} f(u_k) v_k \right\} \cap \left\{ \limsup_{n \to \infty} u_n > a \right\} \subset \left\{ \liminf_{n \to \infty} \frac{u_n - a}{1 + \sum_{k=0}^{n-1} v_n} > 0 \right\}.$$

PROOF. Assume that $u_n \underset{u_n \geq a}{\succeq} \sum_{k=0}^{n-1} f(u_k) v_k$ and $\limsup_{n \to \infty} u_n > a$. Let $\varepsilon > 0$ be such that $\limsup u_n \geq a + 3\varepsilon$ ($\varepsilon$ exists by the second assumption). The first assumption implies there exists $k_0 \in \mathbb{N}$ such that, for all $n > k \geq k_0$, if $u_m \geq a$ for all $m \in [k, n]$,

(10) $$u_n \geq u_k + (1+\varepsilon)^{-1} \left( \sum_{j=k}^{n-1} f(u_j) v_j \right) - \varepsilon.$$

By definition, there exists $k_1 \geq k_0$ such that $u_{k_1} \geq a + 2\varepsilon$. We easily prove by induction, using (10) with $k := k_1$ that, for all $n \geq k_1$, $u_n \geq a + \varepsilon$. It follows from this claim that, for all $n > k_1$,

$$u_n - a \geq \varepsilon + (1 + \varepsilon)^{-1} f(a + \varepsilon) \sum_{j=k_1}^{n-1} v_j$$

$$\geq \min(\varepsilon, (1 + \varepsilon)^{-1} f(a + \varepsilon)) \left( 1 + \sum_{j=k_1}^{n-1} v_j \right),$$

which enables us to conclude the proof.  □

The comparison result stated in Lemma 3.2 gives us a tool which allows us to estimate the behavior of $Z_n(x+6)/Z_n(x+2)$ [resp. $Z_n(x+7)/Z_n(x+2)$] on $\Omega_1(x)$ [resp. on $\Omega_2(x)$]. In particular, part (ii), which provides a sufficient condition for a site to be visited finitely often, implies on $\Omega_1(x)$ that if $Z_n(x+6)/Z_n(x+2)$ does not converge to 1, then either $x+3$ or $x+5$ will be visited finitely often [and a similar result on $\Omega_2(x)$], which allows us to conclude that it is a.s. impossible. The assumption that $A(x)$ holds is technical and easy to check in the cases of application.

LEMMA 3.2.   *For all $x \in \mathbb{Z}$, define the stopping time*

$$T_n(x) = \inf\{m \geq n / X_m = x\}$$

*and the event*

$$A(x) = \Upsilon(x) \cap \left\{ \exists \delta > 0 \Big/ \frac{Z_n(x-2)}{Z_n(x-1)Z_n(x)} = o(\min(Z_n(x-1)^{-\delta}, Z_n(x)^{-\delta})), \right.$$

$$\sup_{n \in \mathbb{N}} \frac{Z_n(x-2)}{Z_n(x-1) + Z_n(x+1)} < \infty,$$

$$\left. \frac{Z_{T_n(x)}(x-2)}{Z_n(x-2)} - 1 = o(Z_n(x-2)^{-\delta}) \right\}.$$

*Then*

(i)  $A(x) \subset \left\{ \ln Z_n(x-1) \equiv \sum_{k=0}^{n-1} \frac{\mathbb{1}_{\{X_k=x\}}}{Z_k(x)} \frac{Z_k(x-2)}{Z_k(x-1) + Z_k(x+1)} \right\},$

(ii)  $A(x) \cap \left\{ \limsup \frac{\ln Z_n(x-2)}{\ln Z_n(x+1)} < 1 \right\} \subset \{Z_\infty(x-1) < \infty\}.$

PROOF.   Assume that $A(x)$ holds. Let us prove that

$$\ln Z_n(x-1) \equiv \sum_{k=1}^{n} \frac{\mathbb{1}_{\{X_{k-1}=x-1\}}}{Z_{k-1}(x-1)}$$

$$\equiv \sum_{k=1}^{n} \frac{\mathbb{1}_{\{X_{k-1}=x-1, X_k=x\}}}{Z_{k-1}(x-1)} \frac{Z_{k-1}(x-2) + Z_{k-1}(x)}{Z_{k-1}(x)}$$

$$\equiv \sum_{k=1}^{n} \frac{\mathbb{1}_{\{X_{k-1}=x-1, X_k=x\}}}{Z_{k-1}(x-1)} \frac{Z_{k-1}(x-2)}{Z_{k-1}(x)}$$

$$\equiv \sum_{k=1}^{n} \frac{\mathbb{1}_{\{X_{k-1}=x, X_k=x-1\}}}{Z_{k-1}(x)} \frac{Z_{k-1}(x-2)}{Z_{k-1}(x-1)}$$

$$\equiv \sum_{k=1}^{n} \frac{\mathbb{1}_{\{X_{k-1}=x\}}}{Z_{k-1}(x)} \frac{Z_{k-1}(x-2)}{Z_{k-1}(x-1) + Z_{k-1}(x+1)}.$$

Indeed, the second equivalence follows from Theorem A.1(i), with

$$M_n = \sum_{k=1}^{n} \frac{\mathbb{1}_{\{X_{k-1}=x-1, X_k=x\}}}{Z_{k-1}(x-1)} \frac{Z_{k-1}(x-2) + Z_{k-1}(x)}{Z_{k-1}(x)} - \sum_{k=1}^{n} \frac{\mathbb{1}_{\{X_{k-1}=x-1\}}}{Z_{k-1}(x-1)}.$$

Indeed, $(M_n)_{n \in \mathbb{N}^*}$ is a square integrable martingale and

$$\langle M \rangle_\infty \leq \sum_{k=1}^{\infty} \frac{\mathbb{1}_{\{X_{k-1}=x-1\}}}{Z_{k-1}(x-1)^2} \frac{Z_{k-1}(x-2) + Z_{k-1}(x)}{Z_{k-1}(x)}$$

$$\leq \sum_{k=1}^{\infty} \frac{\mathbb{1}_{\{X_{k-1}=x-1\}}}{Z_{k-1}(x-1)^2} + \sum_{k=1}^{\infty} \frac{\mathbb{1}_{\{X_{k-1}=x-1\}}}{Z_{k-1}(x-1)} \frac{Z_{k-1}(x-2)}{Z_{k-1}(x-1)Z_{k-1}(x)} < \infty,$$

since $A(x)$ holds. The third equivalence follows from the fact that $\Upsilon(x) = \{\widetilde{Y}_\infty^-(x) < \infty\}$ (by Proposition 3.1) holds. The fourth equivalence follows from an argument similar to the proof of Proposition 3.1(c), using the assumption

$$\frac{Z_{T_n(x)}(x-2)}{Z_n(x-2)} - 1 = o(Z_n(x-2)^{-\delta}),$$

and observing that, by Lemma A.1(i), for all $\delta > 0$, if $A(x)$ holds,

$$\sum_{k=1}^{n} \frac{\mathbb{1}_{\{X_{k-1}=x, X_k=x-1\}}}{Z_{k-1}(x)} \frac{Z_{k-1}(x-2)^{1-\delta}}{Z_{k-1}(x-1)}$$

$$\asymp \sum_{k=1}^{n} \frac{\mathbb{1}_{\{X_{k-1}=x\}}}{Z_{k-1}(x)} \frac{Z_{k-1}(x-2)^{1-\delta}}{Z_{k-1}(x-1) + Z_{k-1}(x+1)}$$

$$\preceq \sum_{k=1}^{n} \frac{\mathbb{1}_{\{X_{k-1}=x\}}}{Z_{k-1}(x)^{1+\delta}} \left( \frac{Z_{k-1}(x-2)}{Z_{k-1}(x-1) + Z_{k-1}(x+1)} \right)^{1-\delta} < \infty,$$

where the second part of the equation follows from $Z_{k-1}(x) \leq Z_{k-1}(x-1) + Z_{k-1}(x+1)$.

To prove the fifth equivalence, we observe that the process

$$R_n = \sum_{k=1}^{n} \frac{\mathbb{1}_{\{X_{k-1}=x, X_k=x-1\}}}{Z_{k-1}(x)} \frac{Z_{k-1}(x-2)}{Z_{k-1}(x-1)}$$

$$- \sum_{k=1}^{n} \frac{\mathbb{1}_{\{X_{k-1}=x\}}}{Z_{k-1}(x)} \frac{Z_{k-1}(x-2)}{Z_{k-1}(x-1) + Z_{k-1}(x+1)}$$

is a martingale and that

$$\langle R \rangle_\infty \leq \sum_{k=1}^{\infty} \frac{\mathbb{1}_{\{X_{k-1}=x\}} \alpha_{k-1}^-(x)}{Z_{k-1}(x)^2} \frac{Z_{k-1}(x-2)^2}{Z_{k-1}(x-1)^2}$$

$$\leq \sum_{k=1}^{\infty} \frac{\mathbb{1}_{\{X_{k-1}=x\}}}{Z_{k-1}(x)} \frac{Z_{k-1}(x-2)}{Z_{k-1}(x)Z_{k-1}(x-1)} \frac{Z_{k-1}(x-2)}{Z_{k-1}(x-1) + Z_{k-1}(x+1)} < \infty$$

if $A(x)$ holds, and we conclude by Theorem A.1(i). Statement (ii) is a direct consequence of (i).  $\square$

3.3. *Proof of Lemma* 2.4: $\Upsilon_0(x) \subset \Upsilon_0'(x)$.  We suppose $x := 0$ for simplicity (the problem is translation-invariant, since the initial point $v_0$ of the VRRW is arbitrary). The inclusion

$$\Upsilon_0(0) \subset \left\{ \limsup \left( \sup_{k \geq n} \alpha_k^-(2)/\alpha_n^-(2) \right) \leq 1 \right\}$$

follows directly from Corollary 3.1(i) applied to site 2, that is,

$$\ln \frac{Z_n(1)}{Z_n(3)} \equiv Y_n^+(0) - Y_n^-(4) \equiv -Y_n^-(4).$$

We first prove the inclusion

(11)                    $$\Upsilon_0(0) \subset \{\limsup Z_n(4)/Z_n(1) \leq e\}.$$

Let us assume that $\Upsilon_0(0)$ holds. It follows from Corollary 3.1(iv)–(v) that $\beta_n^-(2) \to 0$. Fix $\varepsilon > 0$ and $k_0 \in \mathbb{N}$, and assume that for all $n \geq k_0$, $\beta_n^-(2) \leq \varepsilon$ ($\varepsilon > 0$ is chosen in the proof). Let $\mu > e$ and assume, given $p \geq k_0$, that $Z_p(4) \geq \mu Z_p(1)$. Let $(H_n)$ denote the property

$$\forall k \in [p, n], \qquad \frac{Z_k(4)}{Z_k(2)} \geq \beta_p^-(2) = \frac{Z_p(1)}{Z_p(2)} \qquad (H_n).$$

We prove that if $p$ has been chosen large enough, then for all $n \geq p$, $(H_{n+1})$ holds when $(H_n)$ holds. This implies that $(H_n)$ holds for all $n \geq p$ and, therefore, that $\limsup \alpha_n^-(3) < 1$, and subsequently by Corollary 3.2(ii) that $\Upsilon(1)$ holds, which leads to a contradiction by definition of $\Upsilon_0(0)$.

Let

$$\tilde{\alpha}_p^+ = \frac{\beta_p^-(2)}{1 + \beta_p^-(2)} \geq (1+\varepsilon)^{-1}\beta_p^-(2).$$

This proof is based on the following two inequalities, obtained on one hand by Corollary 3.1(i) (applied to $x := 3$) and Proposition 3.1(c) [$Y_n^+(1) \equiv \bar{Y}_n^-(2)$] and on the other hand by Corollary 3.2(i) applied to site 1: If $p$ has been chosen large enough, then for all $n \geq p$,

$$\ln \frac{Z_n(4)}{Z_n(2)} \geq \ln \frac{Z_p(4)}{Z_p(2)} - (\bar{Y}_n^-(2) - \bar{Y}_p^-(2)) - \varepsilon, \tag{12}$$

$$\ln Z_n(1) \leq \ln Z_p(1) + \bar{Y}_n^-(3) - \bar{Y}_p^-(3) + \varepsilon. \tag{13}$$

Recall that

$$\bar{Y}_n^-(3) - \bar{Y}_p^-(3) = \sum_{j=p+1}^{n} \frac{\mathbb{1}_{\{X_j=3\}}}{Z_j(3)} \alpha_j^-(3).$$

We use the following heuristic: As long as $(H_n)$ holds, $\alpha_n^-(3)$ remains far enough from 1, which implies by (13) that $Z_n(1)$ grows slowly in comparison with $Z_n(3)$, which implies that $\bar{Y}_n^-(2) - \bar{Y}_p^-(2)$ remains small, and subsequently by (12) that $\alpha_n^-(3)$ remains far enough from 1.

If $(H_n)$ holds, then for all $j \in [p, n]$, $\alpha_j^-(3) \leq 1 - \tilde{\alpha}_p^+$, which implies

$$\bar{Y}_j^-(3) - \bar{Y}_p^-(3) \leq (1 - \tilde{\alpha}_p^+) \ln \frac{Z_j(3)}{Z_p(3)} + \varepsilon,$$

which implies by (13)

$$\frac{Z_j(1)}{Z_p(1)} \leq e^{2\varepsilon} \left(\frac{Z_j(3)}{Z_p(3)}\right)^{1-\tilde{\alpha}_p^+},$$

and, therefore,

$$\alpha_j^-(2) = \frac{Z_j(1)}{Z_j(1) + Z_j(3)}$$

$$\leq e^{2\varepsilon} \frac{Z_p(1)}{Z_p(3)^{1-\tilde{\alpha}_p^+}} \frac{1}{(Z_j(1) + Z_j(3))^{\tilde{\alpha}_p^+}} \leq e^{2\varepsilon} \frac{Z_p(1)}{Z_p(3)^{1-\tilde{\alpha}_p^+}} \frac{1}{Z_j(2)^{\tilde{\alpha}_p^+}},$$

where we use $Z_j(2) \leq Z_j(1) + Z_j(3)$ in the second inequality. This implies

$$\bar{Y}_n^-(2) - \bar{Y}_p^-(2) \leq e^{2\varepsilon} \frac{Z_p(1)}{Z_p(3)^{1-\tilde{\alpha}_p^+}} \sum_{j=p+1}^{n} \frac{\mathbb{1}_{\{X_j=2\}}}{Z_j(2)^{1+\tilde{\alpha}_p^+}}$$

$$\leq e^{3\varepsilon} \frac{Z_p(1)}{Z_p(3)^{1-\tilde{\alpha}_p^+}} \frac{1}{\tilde{\alpha}_p^+} \frac{1}{Z_p(2)^{\tilde{\alpha}_p^+}}$$

$$\leq \frac{e^{4\varepsilon}}{\tilde{\alpha}_p^+} \frac{Z_p(1)}{Z_p(2)} = e^{4\varepsilon} \frac{\beta_p^-(2)}{\tilde{\alpha}_p^+} \leq (1+\varepsilon)e^{4\varepsilon},$$

where the third inequality follows from $Z_p(2) \leq Z_p(1) + Z_p(3) \leq e^\varepsilon Z_p(3)$ for large enough $p$ [recall that $\alpha_n^-(2) \to 0$]. Therefore, using (12),

$$\frac{Z_{n+1}(4)}{Z_{n+1}(2)} \geq \frac{Z_p(4)}{Z_p(2)} \exp(-\varepsilon - (1+\varepsilon)e^{4\varepsilon}) \geq \frac{Z_p(1)}{Z_p(2)}$$

if $\varepsilon < \mathrm{Cst}(\mu)$, which completes the proof of (11).

The fact that $\lim Z_n(3)/Z_n(2) = 1$ a.s. on $\Upsilon_0(0)$ follows from (11) and $\alpha_\infty^-(2) = 0$, using $Z_n(2) \leq Z_n(1) + Z_n(3)$, and $Z_n(3) \leq Z_n(2) + Z_n(4)$.

Corollary 3.2(ii) and (iii) implies on one hand

$$\Upsilon_0(0) \subset \Upsilon(0) \cap \{\lim \alpha_n^-(3) = 1\} \subset \left\{ \lim \frac{\ln Z_n(1)}{\ln Z_n(2)} = 1 \right\} \tag{14}$$

and on the other hand

$$\Upsilon_0(0) \subset \{\liminf \alpha_n^+(2) = 1\} \subset \left\{ \liminf \frac{\ln Z_n(4)}{\ln Z_n(2)} = 1 \right\}, \tag{15}$$

which gives the third part of the inclusion.

Let us now prove that the $\limsup$ of $Z_n(5)/Z_n(3)$ is less than or equal to 1. Assume $\Upsilon_0(0)$ holds. Let

$$u_n = \ln \frac{Z_n(5)}{Z_n(3)}.$$

By Corollary 3.1(i),

$$u_n = \ln \frac{Z_n(5)}{Z_n(3)} \equiv Y_n^-(6) - Y_n^+(2).$$

For all $a > 0$, using Proposition 3.1(c) and the result given by (11),

$$Y_n^-(6) \equiv \widetilde{Y}_n^-(6) = \sum_{k=1}^n \frac{\mathbb{1}_{\{X_{k-1}=5, X_k=6\}}}{Z_{k-1}(5)} \underset{u_n \geq a}{\succeq} \sum_{k=1}^n \frac{\mathbb{1}_{\{X_{k-1}=5, X_k=6\}}}{Z_{k-1}(6)} \tag{16}$$
$$= Y_n^+(5) \equiv \bar{Y}_n^+(4).$$

The $\succeq$ inequality comes from $Z_n(6) \geq Z_n(5) - Z_n(4) \underset{u_n \geq 0}{\succeq} Z_n(5)$, since $Z_n(4)/Z_n(5) = e^{-u_n} Z_n(4)/Z_n(3) \underset{u_n \geq 0}{\to} 0$, as a consequence of (11) and $\alpha_\infty^-(2) = 0$.

We prove similarly that $Y_n^+(2) \underset{u_n \geq a}{\preceq} \bar{Y}_n^-(4)$, which implies, together with (16),

$$u_n \underset{u_n \geq a}{\succeq} \bar{Y}_n^+(4) - \bar{Y}_n^-(4) = \sum_{k=1}^n \frac{\mathbb{1}_{\{X_n=4\}}}{Z_n(4)} \frac{Z_n(5) - Z_n(3)}{Z_n(3) + Z_n(5)}$$
$$\equiv \sum_{k=0}^{n-1} \frac{\mathbb{1}_{\{X_n=4\}}}{Z_n(4)} \frac{Z_n(5) - Z_n(3)}{Z_n(3) + Z_n(5)}.$$

Let us apply Lemma 3.1 with $f(x) = 1 - e^{-x}$ and

$$v_n = \frac{\mathbb{1}_{\{X_n=4\}}}{Z_n(4)} \frac{Z_n(5)}{Z_n(3) + Z_n(5)} \underset{u_n \geq 0}{\geq} \frac{\mathbb{1}_{\{X_n=4\}}}{2Z_n(4)}.$$

Note that $\sum_{k=0}^{n-1} v_k \succeq_{u_n \geq 0} \ln Z_n(4)/2$.

Using (15), we obtain that for all $a > 0$, $\limsup u_n > a$ implies

$$\liminf \frac{\ln Z_n(5)}{\ln Z_n(4)} = \liminf \frac{\ln Z_n(5)}{\ln Z_n(3)} = 1 + \liminf \frac{\ln(Z_n(5)/Z_n(3))}{\ln Z_n(4)} > 1$$

and, accordingly, that $\Upsilon(4)$ holds. Remark 3.1 implies $\alpha_\infty^-(2) > 0$, which leads to a contradiction on $\Upsilon_0(0)$. The proof concerning $Z_n(0)/Z_n(2)$ is similar.

The statement $Z_\infty(0) = Z_\infty(4) = \infty$ follows from $Z_\infty(0) \succeq Y_\infty^-(1) \approx Y_\infty(1) = \infty$ a.s. on $\Upsilon_0(0) \subset \Upsilon(1)^c$ [by Proposition 3.1(a)], which implies $Z_\infty(4) = \infty$ by the other statements of this lemma.

3.4. *Proof of Lemma* 2.9. We suppose $x := 0$ for simplicity. Let us prove that, on $\Omega_1(0)$, the $\limsup$ of $Z_n(6)/Z_n(2)$ is less than or equal to 1. The symmetrical statement (with respect to site 4) completes the proof. Our goal is to describe the evolution of the quantity $Z_n(6)/Z_n(2)$; this description is obtained in (18). Assume in the sequel that $\Omega_1(0) \cap \Upsilon(0,4,8)$ holds.

First, we apply Lemma 3.2(i). Let us use its notation for $x := 4$ and prove that $A(4)$ holds. Using Remark 3.1 and Corollary 3.2(iii), we obtain that $\alpha_n^+(2)$ and $\beta_n^+(2)$ converge to $\alpha_\infty^+(2) > 0$ and that $\ln Z_n(4) \succeq \alpha_\infty^+(2) \ln Z_n(2)$. Hence, on one hand,

$$\limsup \frac{Z_n(2)}{Z_n(3) + Z_n(5)} \leq \lim \beta_n^+(2)^{-1} < \infty,$$

$$\frac{Z_n(2)}{Z_n(3) Z_n(4)} \asymp (\alpha_\infty^+(2))^{-1} Z_n(4)^{-1} = O(Z_n(4)^{-1}) = O(Z_n(3)^{-\alpha_\infty^+(2)}).$$

On the other hand, let us prove that there exists a.s. $\delta > 0$ such that

$$(17) \qquad \frac{Z_{T_n(4)}(2)}{Z_n(2)} - 1 = o(Z_n(2)^{-\delta}).$$

Indeed, for all $n \in \mathbb{N}^*$, let $t_n$ be the $n$th visit time to site 4. For all $m \in \mathbb{N}$ and $a > 0$, let $T^{m,a} := \inf\{n \geq m \text{ s.t. } \alpha_n^+(2) \leq a\}$. There exists a.s. $a > 0$ and $m \in \mathbb{N}$ such that $T^{m,a} = \infty$.

Given $a, \varepsilon > 0$ and $n \in \mathbb{N}$, let

$$\Gamma_{n+1} := \{Z_{t_{n+1}}(2) - Z_{t_n}(2) > Z_{t_n}(2) n^{\varepsilon-1}\} \cup \{T^{m,a} \leq t_{n+1}\}.$$

Let $n \in \mathbb{N}$ be such that $t_n \geq m$ and assume $n \geq \text{Cst}$. Given $t \geq t_n$ such that $t < T^{m,a}$, $X_t = 2$ and $Z_t(2) \leq Z_{t_n}(2)(1 + n^{\varepsilon-1})$, the probability to reach site

4 in two steps starting from site 2 at time $t$ is greater than $an/(2Z_{t_n}(2))$ and, therefore,

$$\mathbb{P}(\Gamma_{n+1}^c | \mathcal{F}_{t_n}) \leq \left(1 - \frac{an}{2Z_{t_n}(2)}\right)^{Z_{t_n}(2)n^{\varepsilon-1}/2} \leq \exp\left(-\frac{an^\varepsilon}{4}\right).$$

Accordingly,

$$\sum_{n \in \mathbb{N}^*} \mathbb{P}(\Gamma_n^c) < \infty \qquad \text{a.s.}$$

and the Borel–Cantelli lemma implies that $\Gamma_n^c$ occurs only finitely often. This gives (17), using $\ln Z_n(4) \succeq \alpha_\infty^+(2) \ln Z_n(2)$.

Therefore, Lemma 3.2(i) implies, together with $\beta_n^+(2) \to \alpha_\infty^+(2) > 0$,

$$\ln Z_n(2) \equiv \ln Z_n(3) \equiv \sum_{k=0}^{n-1} \frac{\mathbb{1}_{\{X_k=4\}}}{Z_k(4)} \frac{Z_k(2)}{Z_k(3) + Z_k(5)}.$$

The situation being symmetrical with respect to site 4, we have a similar estimate for $\ln Z_n(6)$. Hence,

$$(18) \qquad \ln \frac{Z_n(6)}{Z_n(2)} \equiv \sum_{k=0}^{n-1} \frac{\mathbb{1}_{\{X_k=4\}}}{Z_k(4)} \frac{Z_k(6) - Z_k(2)}{Z_k(3) + Z_k(5)}.$$

Let us apply Lemma 3.1 with $a = 0$, $u_n = \ln(Z_n(6)/Z_n(2))$, $f(x) = 1 - e^{-x}$ and

$$v_n = \frac{\mathbb{1}_{\{X_n=4\}}}{Z_n(4)} \frac{Z_n(6)}{Z_n(3) + Z_n(5)} \underset{u_n \geq 0}{\geq} (\beta_n^+(2) + \beta_n^-(6))^{-1} \frac{\mathbb{1}_{\{X_n=4\}}}{2Z_n(4)},$$

using that $Z_n(3) + Z_n(5) = \beta_n^+(2)Z_n(2) + \beta_n^-(6)Z_n(6)$. We obtain, if $\limsup Z_n(6)/Z_n(2) > 1$,

$$\liminf \frac{\ln Z_n(6)}{\ln Z_n(2)} = 1 + \liminf \frac{\ln(Z_n(6)/Z_n(2))}{\ln Z_n(4)} \frac{\ln Z_n(4)}{\ln Z_n(2)} > 1,$$

using the estimate $\ln Z_n(4) \succeq \alpha_\infty^+(2) \ln Z_n(2)$.

Now, Lemma 3.2(ii) completes the proof. Indeed,

$$\Omega_1(0) \cap \left\{\limsup \frac{Z_n(6)}{Z_n(2)} > 1\right\}$$

$$\subset \Omega_1(0) \cap A(4) \cap \left\{\limsup \frac{\ln Z_n(2)}{\ln Z_n(5)} < 1\right\}$$

$$\subset \Omega_1(0) \cap \{Z_\infty(3) < \infty\} = \varnothing \qquad \text{a.s.},$$

where we use in the first inclusion that $\beta_n^-(6) \to \alpha_\infty^-(6) > 0$, by Remark 3.1.

3.5. *Proof of Lemma 2.10.* Suppose $x := 0$ for simplicity. We use Lemma A.2 in the Appendix. Let us introduce a some notation first. Let, for all $n \in \mathbb{N}^*$,

$$T_n := \inf\{k \in \mathbb{N} \text{ s.t. } Z_k(3) = n \text{ or } Z_k(6) = n\}$$

and let $\mathbb{G} := (\mathcal{G}_n)_{n \in \mathbb{N}^*} := (\mathcal{F}_{T_n})_{n \in \mathbb{N}^*}$. We easily prove by induction that, for all $n \geq 2$, $X_{T_n} \in \{3, 6\}$ and $Z_{T_n}(X_{T_n}) = n = \max(Z_{T_n}(3), Z_{T_n}(6))$. For all $n \geq 2$, let us define $\bar{X}_{T_n} := 5$ if $X_{T_n} = 3$ and $:= 4$ if $X_{T_n} = 6$. Let $\Gamma_0 = \Gamma_1 = \Gamma_2 := \varnothing$ and, for all $n \geq 2$,

$$\Gamma_{n+1} = \{T_n < \infty\} \cap \{X_{T_n+2} = \bar{X}_{T_n}\}.$$

Assume that $\Omega_2(0) \cap \{\alpha_\infty^-(7) \leq \alpha_\infty^-(2)\}$ holds. Let us apply Lemma A.2 to prove that $\Gamma_n$ holds only finitely often a.s., which implies that $Z_\infty(4) \wedge Z_\infty(5) < \infty$ or $Z_\infty(3) \wedge Z_\infty(6) < \infty$ (if $\exists n \in \mathbb{N}$ s.t. $T_n = \infty$) and, therefore, enables us to conclude.

We settle upon the notation of Lemma A.2. Let us choose the sequence $(\gamma_n)$ that satisfies the upper bound of $\mathbb{P}(\Gamma_{n+1}|\mathcal{F}_{t_n})$. Given $\delta > 0$ and $m \in \mathbb{N}$, let

$$A_{\delta,m} := \{\forall n \geq m, |\alpha_\infty^-(2) - \alpha_n^-(2)| \leq Z_n(3)^{-\delta} \text{ and } |\alpha_\infty^+(7) - \alpha_n^+(7)| \leq Z_n(6)^{-\delta}\}.$$

By Remark 3.1 $[\beta_\infty^+(2) = \alpha_\infty^+(2) > 0$ and $\beta_\infty^-(6) = \alpha_\infty^-(6) > 0]$ and Corollary 3.2(iv), there exists a.s. $\delta > 0$, $m \in \mathbb{N}$ such that $A_{\delta,m}$ holds. We fix $\delta > 0$, $m \in \mathbb{N}$ and suppose $A_{\delta,m}$ holds. We choose $\gamma_n := (\alpha_{T_n}^-(2) - n^{-\delta})\mathbb{1}_{X_{T_n}=3} + (\alpha_{T_n}^+(7) - n^{-\delta})\mathbb{1}_{X_{T_n}=6}$. Note that $\gamma_n \leq \alpha_\infty^-(2)\mathbb{1}_{X_{T_n}=3} + \alpha_\infty^+(7)\mathbb{1}_{X_{T_n}=6}$ by definition of $A_{\delta,m}$.

By Remark 3.1 and Corollary 3.2(iv), there exists a.s. $h > 0$ such that, for all $n \in \mathbb{N}$, $Z_n(4) \leq hZ_n(3)^{\alpha_\infty^+(2)}$ and $Z_n(5) \leq hZ_n(6)^{\alpha_\infty^-(7)}$. Then, if $n \geq m$ and $X_{T_n} = 3$, using $Z_{T_n}(6) \leq \tau_n$ and $Z_{T_n}(3) = n$,

$$\mathbb{P}(\Gamma_{n+1}|\mathcal{F}_{T_n}) = \alpha_{T_n}^+(3)\alpha_{T_n}^+(4) \leq \beta_{T_n}^+(3)\alpha_{T_n}^+(4) \leq Z_{T_n}(4)Z_{T_n}(5)/Z_{T_n}(3)^2$$

$$\leq h^2 Z_{T_n}(6)^{\alpha_\infty^-(7)} Z_{T_n}(3)^{\alpha_\infty^+(2)-2} \leq h^2 \tau_n^{\alpha_\infty^-(2)}/n^{1+\alpha_\infty^-(2)}$$

$$\leq h^2 \tau_n^{\gamma_{\tau_n}}/n^{1+\gamma_{\tau_n}},$$

using $\alpha_\infty^-(7) \leq \alpha_\infty^-(2)$ in the fourth inequality and using, in the last inequality, that $\gamma_{\tau_n} \leq \alpha_\infty^-(2)$ and that $x \mapsto \tau_n^x/n^{1+x}$ is nonincreasing on $\mathbb{R}_+$ (since $\tau_n \leq n$). The estimate of $\mathbb{P}(\Gamma_{n+1}|\mathcal{F}_{T_n})$ is very similar when $X_{T_n} = 6$, which enables us to conclude.

3.6. *Proof of Lemma 2.11.* We suppose $x := 0$ for simplicity. Let us prove that, on $\Omega_2(0)$, the $\limsup$ of $Z_n(7)/Z_n(2)$ is less than or equal to 1. The symmetrical statement (with respect to the number 4.5) completes the proof. Assume subsequently that $\Omega_2(0)$ holds. Let $u_n = \ln(Z_n(7)/Z_n(2))$. Similarly as in the proof of Lemma 2.9, we try to describe the evolution of the quantity $u_n$; this description is obtained in (20).

Let us begin with some elementary properties. By Remark 3.1, $\alpha_\infty^\pm(2) = \beta_\infty^\pm(2) \in (0,1)$ and $\alpha_\infty^\pm(7) = \beta_\infty^\pm(7) \in (0,1)$. By Corollary 3.2(iv), there exist a.s. $\gamma_\infty^1, \gamma_\infty^2 > 0$ such that $Z_n(4) \asymp \gamma_\infty^1 Z_n(2)^{\alpha_\infty^+(2)}$ and $Z_n(5) \asymp \gamma_\infty^2 Z_n(7)^{\alpha_\infty^-(7)}$.

Now, we can adapt the proof of Lemma 3.2(i) for $x := 5$ to show that

$$\ln Z_n(4) \underset{u_n \geq 0}{\preceq} \sum_{k=1}^n \frac{\mathbb{1}_{\{X_k=5\}}}{Z_k(5)} \frac{Z_k(3)}{Z_k(4) + Z_k(6)}. \tag{19}$$

Indeed, there exists a.s. $\delta > 0$ such that

$$\frac{Z_n(3)}{Z_n(4)Z_n(5)} \asymp \frac{\beta_\infty^+(2) Z_n(2)}{\gamma_\infty^1 \gamma_\infty^2 Z_n(2)^{\alpha_\infty^+(2)} Z_n(7)^{\alpha_\infty^-(7)}}$$

$$\underset{u_n \geq 0}{=} O(Z_n(7)^{-\delta_\infty}) \underset{u_n \geq 0}{=} O((Z_n(4) \vee Z_n(5))^{-\delta}),$$

letting $\delta_\infty := \alpha_\infty^-(7) - \alpha_\infty^-(2)$ ($> 0$ by Lemma 2.10) and, on the other hand,

$$\frac{Z_n(3)}{Z_n(4) + Z_n(6)} \leq \frac{Z_n(3)}{Z_n(6)} \asymp \frac{\beta_\infty^+(2)}{\beta_\infty^-(7)} e^{-u_n} = \frac{\alpha_\infty^+(2)}{\alpha_\infty^-(7)} e^{-u_n},$$

$$\frac{Z_{T_n(5)}(3)}{Z_n(3)} - 1 = o(Z_n(3)^{-\delta}).$$

The last equality comes from the fact that when $u_n \geq 0$, the probability to go from 3 to 5 is greater than a term on the order of $Z_n(3)^{\varepsilon-1}$ for $\varepsilon > 0$, using Lemma 2.10, and its proof is very similar to the proof of (17) in Section 3.4.

Inequality (19) follows, which implies with the upper bound of $Z_n(3)/(Z_n(4) + Z_n(6))$ that

$$\ln Z_n(4) \equiv \alpha_\infty^+(2) \ln Z_n(2) \underset{u_n \geq 0}{\preceq} \sum_{k=1}^n \frac{\mathbb{1}_{\{X_k=5\}}}{Z_k(5)} \frac{\alpha_\infty^+(2) Z_k(2)}{\alpha_\infty^-(7) Z_k(7)}.$$

Hence

$$\ln \frac{Z_n(7)}{Z_n(2)} \underset{u_n \geq 0}{\succeq} \sum_{k=0}^{n-1} \frac{\mathbb{1}_{\{X_k=5\}}}{\alpha_\infty^-(7) Z_k(5)}\left(1 - \frac{Z_k(2)}{Z_k(7)}\right) \tag{20}$$

since, by Corollary 3.2(iv),

$$\ln Z_n(7) \equiv \frac{\ln Z_n(5)}{\alpha_\infty^-(7)} \equiv \sum_{k=0}^{n-1} \frac{\mathbb{1}_{\{X_k=5\}}}{\alpha_\infty^-(7) Z_k(5)}.$$

Now, Lemma 3.1 with $a = 0$, $u_n$ defined below,

$$v_n = \frac{\mathbb{1}_{\{X_k=5\}}}{\alpha_\infty^-(7) Z_k(5)}, \qquad f(x) = 1 - e^{-x},$$

implies that $\liminf u_n / \ln Z_n(5) > 0$ and, therefore, $\limsup \ln Z_n(3) / \ln Z_n(6) < 1$, since

$$\liminf \frac{\ln Z_n(6)}{\ln Z_n(3)} = \liminf \frac{\ln Z_n(7)}{\ln Z_n(2)} = 1 + \liminf \frac{u_n}{\ln Z_n(5)} \frac{\ln Z_n(5)}{\ln Z_n(2)} > 1.$$

Hence we can conclude, again by inequality (19) [as in Lemma 3.2(ii) for $x := 5$], that $Z_\infty(4) < \infty$, which is a.s. impossible [on $\Omega_2(0)$].

## 4. Nonconvergence toward unstable situations.

4.1. *Introduction.* The aim of this section is to provide a result that ensures nonconvergence in the unstable situations that correspond to Lemmas 2.5, 2.7 and 2.8, which are proved, respectively, in Sections 5.1, 5.2 and 5.3. This result makes use of the particular structure of reinforced random walks to overcome the fact that, in general, we can only obtain partial information on the behavior of the random walk. Indeed, as explained in the Introduction, it is not, in general, possible to describe the behavior of the density of occupation of the random walk by the differential equation (3), which would enable us to interpret these unstable situations by unstable sets of the corresponding dynamical system and, therefore, allow us to use the classical results of nonconvergence toward these sets. For this reason, we provide a result that requires only an equation of evolution of the considered unstable quantity. This information is sufficient when the evolution involved is in some sense compatible with a partial order constructed on a certain class of random walks on $\mathbb{Z}$. More precisely, we study the behavior of a $(\mathcal{G}_n)_{n\in\mathbb{N}}$-adapted sequence $(z_n)_{n\in\mathbb{N}}$ that takes values in $\mathbb{R}$ and we try to prove that its behavior around 0 is unstable, so that convergence to 0 is a.s. impossible.

The evolution of $(z_n)_{n\in\mathbb{N}}$ is given by an equation of the form

$$(21) \qquad z_{n+1} - z_n = y_n + \varepsilon_{n+1} + r_n,$$

where $(y_n)_{n\in\mathbb{N}}$, $(\varepsilon_n)_{n\in\mathbb{N}^*}$ and $(r_n)_{n\in\mathbb{N}}$ are $(\mathcal{G}_n)_{n\in\mathbb{N}}$-adapted and

$$\mathbb{E}(\varepsilon_{n+1}|\mathcal{G}_n) = 0.$$

Let us, for instance, consider the case discussed in Lemma 2.5: $(z_n)_{n\in\mathbb{N}}$ and $(y_n)_{n\in\mathbb{N}}$ correspond to

$$z_n = \ln \frac{Z_{t_n}(3)}{Z_{t_n}(2)} \quad \text{and} \quad y_n = \frac{R_{t_n}}{Z_{t_n}(2) Z_{t_n}(3)},$$

where

$$t_n = \inf\{m \in \mathbb{N} / Z_m^+(2) \geq n\},$$
$$R_n = Z_n(4) + Z_n(2) - (Z_n(1) + Z_n(3)),$$

and

$$\mathbb{E}(\varepsilon_{n+1}^2 | \mathcal{G}_n) \asymp \alpha_{t_n}^-(2)/n^2, \qquad |r_n| = O(1/n^{2-\varepsilon}) \qquad \text{for all } \varepsilon > 0.$$

If $z_n$ and $y_n$ were of the same sign, we would be able to conclude (see, e.g., [8], [11] and [12], Chapter 3) that the unstable point $z = 0$ is a.s. avoided, namely that

$$(22) \qquad\qquad \mathbb{P}\left(\lim_{n\to\infty} z_n = 0\right) = 0.$$

This is not the case here and, in fact, $y_n$ does not depend only on $z_n$. However, we can observe that the term $R_{t_n}$ increases only with visits from 5 to 4 and decreases with visits from 0 to 1. Indeed, it is easy to prove by induction that, for all $n \in \mathbb{N}$,

$$(23) \quad R_n = Z_n^-(5) - Z_n^+(0) + (\mathbb{1}_{\{X_n=2 \text{ or } X_n \geq 4\}} - \mathbb{1}_{\{X_n \leq 1 \text{ or } X_n=3\}})/2 + \text{Cst}(v_0),$$

so that

$$(24) \qquad\qquad R_{t_n} = Z_{t_n}^+(4) - Z_{t_n}^-(1) + \text{Cst}(v_0),$$

using that $Z_{t_n}^-(5) = Z_{t_n}^+(4) + \text{Cst}(v_0)$ and $Z_{t_n}^+(0) = Z_{t_n}^-(1) + \text{Cst}(v_0)$.

Heuristically, when $z_n$ tends to increase (resp. to decrease), the random walk tends to go more to the right (resp. to the left), which implies that $y_n$ also tends to increase. The precise tool behind these remarks is the definition of a partial order on the random walks. Lemma 4.1 claims the following result.

Assume we deal with two random walks $\mathcal{M}$ and $\mathcal{M}'$ such that at each point $j \in \mathbb{Z}$, for the same number of visits to $j$, if $\mathcal{M}'$ has more visited $j+1$ than $\mathcal{M}$ and less visited $j-1$, then $\mathcal{M}'$ has a greater probability than $\mathcal{M}$ to go right. Then we can couple $\mathcal{M}$ and $\mathcal{M}'$ so that for all $j \in \mathbb{Z}$, for the same number of visits to $j$, $\mathcal{M}'$ has more visited $j+1$ than $\mathcal{M}$ and less visited $j-1$. In this case, we write that $\mathcal{M}' \gg \mathcal{M}$.

It is easy to prove that, given two random walks $\mathcal{M}$ and $\mathcal{M}'$ such that $\mathcal{M}' \gg \mathcal{M}$, if we keep the same notation for $\mathcal{M}$ and add a superscript prime for $\mathcal{M}'$, then, for all $n \in \mathbb{N}$,

$$R'_{t'_n} \geq R_{t_n}.$$

Having put down this partial order on random walks on $\mathbb{Z}$, we observe in the considered cases that a significant part of the noise inherent in the behavior of $z_n$ is generated by the uncertainty on the visits from a certain vertex $v$ to $v-1$ (in the case of Lemma 2.5, $v := 2$). This leads us to define, concurrent to the VRRW called $\mathcal{M}$, a random walk $\mathcal{M}'$ as follows. Starting from all points except from $v$, $\mathcal{M}'$ has the same conditional probabilities as $\mathcal{M}$. From $v$ the conditional probability to visit $v-1$ is the probability designed for $\mathcal{M}$ minus a term on the order of the standard deviation of this probability on a large time interval. This new random walk is constructed in Definition 4.12.

Lemma 4.1 implies that we can couple $\mathcal{M}$ and $\mathcal{M}'$ so that $\mathcal{M}' \gg \mathcal{M}$. This property has the consequence that, roughly, $y'_n$ is greater than $y_n$; more precisely, Assumption (H3) of Proposition 4.1 is satisfied. Furthermore, the different moving probabilities from $v$ imply here that $z_n$ in $z'_n$ undergoes a drift toward the right significant enough to cover the noise, which corresponds to Assumptions (H1) and (H2) [$W(k)$ is of the order of the standard deviation of $z_n$ starting at time $t_k$]. The probabilities of a same group of paths for $\mathcal{M}$ and for $\mathcal{M}'$ being of the same order (stated in Lemma 4.2), these properties imply that the conditional probability not to converge to 0 is always greater than a positive constant. This enables us to conclude that $z_n$ a.s. does not converge to 0.

The section is divided as follows. In Section 4.2, we introduce some notation, and state and prove a coupling result for nearest-neighbor random walks on $\mathbb{Z}$. In Section 4.3, we state the nonconvergence result Proposition 4.1, which is applied in Sections 5.1–5.3. Proposition 4.1 is proved in Section 4.4.

### 4.2. *Notation and a coupling result.*

DEFINITION 4.1. Given $k \in \mathbb{N} \cup \{\infty\}$ and $\mathbf{v} \in \mathbb{Z}^{k+1}$, for all $n \leq k$, we let $\mathbf{v}_n$ be the $(n+1)$th coordinate of $\mathbf{v}$. We say that $\mathbf{v}$ is a $k$ path (or a path, when there is no ambiguity) on $\mathbb{Z}$ iff, for all $0 \leq n \leq k-1$, there exists $\varepsilon_n \in \{-1, 1\}$ such that $\mathbf{v}_{n+1} - \mathbf{v}_n = \varepsilon_n$. Let $\mathcal{P}_k$ be the set of $k$ paths. Let $\mathcal{P} := \mathcal{P}_\infty$. Given $i \in \mathbb{N}^*$, $j \in \mathbb{Z}$ and $\mathbf{v} \in \mathcal{P}_k$, we let $n_{i,j}(\mathbf{v})$ be the time the sequence $(\mathbf{v}_n)_{0 \leq n \leq k}$ makes its $i$th visit to site $j$, with the convention that $n_{i,j}(\mathbf{v}) = \infty$ if $j$ is visited less than $i$ times.

We make use of the notation introduced in Sections 1 and 2.2, that is, $Z_n(x)(\mathbf{v})$, $Y_n(x)(\mathbf{v}), \ldots$, for all $n \leq k$, replacing the underlying $(X_j)_{j \in \mathbb{N}}$ in these definitions by $(\mathbf{v}_j)_{j \in \mathbb{N}}$.

DEFINITION 4.2. For all $k \in \mathbb{N}$, let $\mathcal{T}_k$ be the smallest $\sigma$-field on $\mathcal{P}$ that contains the cylinders $C_\mathbf{v} = \{\mathbf{w} \in \mathcal{P}/\mathbf{w}_0 = \mathbf{v}_0, \ldots, \mathbf{w}_k = \mathbf{v}_k\}$, $\mathbf{v} \in \mathcal{P}_k$. Let $\mathcal{T} := \vee_{k \in \mathbb{N}} \mathcal{T}_k$. Let us define the filtration $\mathbb{T} := (\mathcal{T}_k)_{k \in \mathbb{N}}$.

DEFINITION 4.3. On a probability space $(\Omega, \mathcal{F}, \mathbb{P})$, we call a random walk a process $(X_k)_{k \in \mathbb{N}}$ taking values in $\mathbb{Z}$, starting from a fixed point $X_0 := v_0 \in \mathbb{Z}$, satisfying the following conditions: for a.e. $\omega \in \Omega$, $(X_i(\omega))_{i \in \mathbb{N}} \in \mathcal{P}$ and, for all $k \in \mathbb{N}^*$ and $\mathbf{v} \in \mathcal{P}_k$ such that $\mathbf{v}_0 = v_0$, $\mathbb{P}((X_i)_{0 \leq i \leq k} = \mathbf{v}) > 0$.

DEFINITION 4.4. Let $\mathcal{M} := (X_k)_{k \in \mathbb{N}}$ be a random walk on a probability space $(\Omega, \mathcal{F}, \mathbb{P})$. Let $\mathcal{I}_{\mathcal{M}} : (\Omega, \mathcal{F}) \to (\mathcal{P}, \mathcal{T})$ be the measurable function $\omega \in \Omega \mapsto (X_i(\omega))_{i \in \mathbb{N}}$. Note that $\mathcal{I}_{\mathcal{M}}$ defines naturally a probability on $\mathcal{P}$ by, for all $\mathcal{C} \in \mathcal{T}$, $\mathbb{P}_{\mathcal{M}}(\mathcal{C}) := \mathbb{P}(\mathcal{I}_{\mathcal{M}}^{-1}(\mathcal{C}))$. For all $n \in \mathbb{N}$, let $\mathcal{E}_n^{\mathcal{M}} := \mathcal{I}_{\mathcal{M}}^{-1}(\mathcal{T}_n) = \sigma(X_0, \ldots, X_n)$. Let us define the filtration $\mathbb{E}^{\mathcal{M}} := (\mathcal{E}_n^{\mathcal{M}})_{n \in \mathbb{N}}$. Given a $\mathcal{T}$-measurable random variable $u$, let us define the $\mathcal{F}$-measurable random variable $u^{\mathcal{M}} := u \circ \mathcal{I}_{\mathcal{M}}$.

Note that, if $T$ is a $\mathbb{T}$ stopping time, then $T^{\mathcal{M}}$ is a $\mathbb{E}^{\mathcal{M}}$ stopping time. If $(t_n)_{n \in \mathbb{N}}$ is a nondecreasing sequence of $\mathbb{T}$ stopping times and if $T$ [resp. $(a_n)_{n \in \mathbb{N}}$] is a $(\mathcal{T}_{t_n})_{n \in \mathbb{N}}$ stopping time (resp. adapted process), then $T^{\mathcal{M}}$ [resp. $(a_n^{\mathcal{M}})_{n \in \mathbb{N}}$] is a $(\mathcal{E}_{t_n^{\mathcal{M}}}^{\mathcal{M}})_{n \in \mathbb{N}}$ stopping time (resp. adapted process).

DEFINITION 4.5. Let $\mathcal{M} := (X_k)_{k \in \mathbb{N}}$ be a random walk on a probability space $(\Omega, \mathcal{F}, \mathbb{P})$, starting from $v_0 \in \mathbb{Z}$. For all $k \in \mathbb{N} \cup \{\infty\}$, $\mathbf{v} \in \mathcal{P}_k$ such that $\mathbf{v}_0 = v_0$ and $n \leq k$, let $q_{\mathcal{M}}(\mathbf{v}, n) := \mathbb{P}_{\mathcal{M}}(\mathcal{C}_{(\mathbf{v}_0, \ldots, \mathbf{v}_n, \mathbf{v}_n+1)})/\mathbb{P}_{\mathcal{M}}(\mathcal{C}_{(\mathbf{v}_0, \ldots, \mathbf{v}_n)})$ be the conditional probability to go to the right at time $n$, knowing $X_0 = \mathbf{v}_0, \ldots, X_n = \mathbf{v}_n$. For all $\mathbf{v} \in \mathcal{P}$, $i \in \mathbb{N}^*$, $j \in \mathbb{Z}$ such that $n_{i,j}(\mathbf{v}) < \infty$, let $p_{i,j}^{\mathcal{M}}(\mathbf{v}) := q_{\mathcal{M}}(\mathbf{v}, n_{i,j}(\mathbf{v}))$ be the conditional probability to go to the right just after the $i$th visit to site $j$, knowing $X_0 = \mathbf{v}_0, \ldots$ and $X_{n_{i,j}(\mathbf{v})} = \mathbf{v}_{n_{i,j}(\mathbf{v})}$.

Subsequently, we fix the probability space $(\Omega, \mathcal{F}, \mathbb{P})$, on which we take i.i.d. uniform $[0,1]$ random variables $(\omega_{i,j})_{i \in \mathbb{N}^*, j \in \mathbb{Z}}$.

DEFINITION 4.6. We construct (and settle) the random walks on $(\Omega, \mathcal{F}, \mathbb{P})$ by the following method. Given the initial point $v_0$ and the conditional probabilities of move $q_{\mathcal{M}}(\cdot, \cdot)$ of a random walk $\mathcal{M}$, we let $\mathcal{M} := (X_k)_{k \in \mathbb{N}}$ on $\Omega$ be as follows: $X_0 := v_0$ and, for all $n \in \mathbb{N}$, given $(X_0, \ldots, X_n)$,

$$X_{n+1} = \begin{cases} X_n + 1, & \text{if } \omega_{Z_n(X_n)(X_0, \ldots, X_n)-1, X_n} \leq q_{\mathcal{M}}((X_0, \ldots, X_n), n), \\ X_n - 1, & \text{otherwise.} \end{cases}$$

DEFINITION 4.7. Let $\mathcal{M} := (X_k)_{k \in \mathbb{N}}$ be a random walk [on the probability space $(\Omega, \mathcal{F}, \mathbb{P})$]. Let, for all $n \in \mathbb{N}$, $\mathcal{F}_n^{\mathcal{M}} := \sigma(\omega_{i,j})_{(i,j) \in \mathbb{N}^* \times \mathbb{Z}/n_{i,j}(\mathcal{I}_{\mathcal{M}}(\omega)) \leq n} = \sigma(\{\omega_{i,j} \in I\} \cap \{n_{i,j}(\mathcal{I}_{\mathcal{M}}(\omega)) \leq n\}, i \in \mathbb{N}^*, j \in \mathbb{Z}, I \subset [0,1]$ interval). Note that $\mathcal{E}_n^{\mathcal{M}} \subset \mathcal{F}_n^{\mathcal{M}}$. Let us define the filtration $\mathbb{F}^{\mathcal{M}} := (\mathcal{F}_n^{\mathcal{M}})_{n \in \mathbb{N}}$.

DEFINITION 4.8. Given $\mathbf{v}, \mathbf{v}' \in \mathcal{P}$, let us define, for all $i \in \mathbb{N}^*$ and $j \in \mathbb{Z}$, the property $E_{i,j}(\mathbf{v}, \mathbf{v}')$ as

$$Z_{n_{i,j}(\mathbf{v}')}(j+1)(\mathbf{v}') \geq Z_{n_{i,j}(\mathbf{v})}(j+1)(\mathbf{v}) \quad \text{and}$$

$$Z_{n_{i,j}(\mathbf{v}')}(j-1)(\mathbf{v}') \leq Z_{n_{i,j}(\mathbf{v})}(j-1)(\mathbf{v})$$

with the convention that $E_{i,j}(\mathbf{v}, \mathbf{v}')$ holds whenever $n_{i,j}(\mathbf{v}) = \infty$ or $n_{i,j}(\mathbf{v}') = \infty$.

DEFINITION 4.9. Let $\mathcal{M}$ and $\mathcal{M}'$ be two random walks on $(\Omega, \mathcal{F}, \mathbb{P})$. Let $\mathcal{M}' \gg \mathcal{M}$ denote the following property: for a.e. $\omega \in \Omega$, $E_{i,j}(\mathcal{I}_{\mathcal{M}}(\omega), \mathcal{I}_{\mathcal{M}'}(\omega))$ holds for all $i \in \mathbb{N}^*$ and $j \in \mathbb{Z}$.

Thus $\mathcal{M}' \gg \mathcal{M}$ means that, for the same number $i$ of visits to $j$, $\mathcal{M}'$ has visited site $j + 1$ (right hand from $j$) more often than $\mathcal{M}$ and has visited $j - 1$ less often than $\mathcal{M}$.

LEMMA 4.1. *Let $\mathcal{M}$ and $\mathcal{M}'$ denote two random walks on $\mathbb{Z}$ [on the probability space $(\Omega, \mathcal{F}, \mathbb{P})$], starting from the same point $X_0 = X_0' = v_0$. Suppose that for all $\mathbf{v}$, $\mathbf{v}' \in \mathcal{P}$, for all $i \in \mathbb{N}^*$, $j \in \mathbb{Z}$, $p_{i,j}^{\mathcal{M}'}(\mathbf{v}') \geq p_{i,j}^{\mathcal{M}}(\mathbf{v})$ whenever $E_{i,j}(\mathbf{v}, \mathbf{v}')$ holds and $\max(n_{i,j}(\mathbf{v}), n_{i,j}(\mathbf{v}')) < \infty$. Then $\mathcal{M}' \gg \mathcal{M}$.*

PROOF. Consider an arbitrary element $\omega \in \Omega$. Let $\mathbf{v} := \mathcal{I}_{\mathcal{M}}(\omega)$ and $\mathbf{v}' := \mathcal{I}_{\mathcal{M}'}(\omega)$, $p_{i,j} := p_{i,j}^{\mathcal{M}}(\mathbf{v})$, $p_{i,j}' := p_{i,j}^{\mathcal{M}'}(\mathbf{v}')$ and $E_{i,j} := E_{i,j}(\mathbf{v}, \mathbf{v}')$. We want to prove that, for all $i \in \mathbb{N}^*$, $j \in \mathbb{Z}$, $E_{i,j}$ holds. Observe the following facts:

- One has $\mathbf{v}_0 = v$.
- For all $p \in \mathbb{N}$, there exists $i \in \mathbb{N}^*$ and $j \in \mathbb{Z}$ such that $p = n_{i,j}(\mathbf{v})$:
  - if $\omega_{i,j} \leq p_{i,j}$, then $\mathbf{v}_{p+1} = \mathbf{v}_p + 1$;
  - if $\omega_{i,j} > p_{i,j}$, then $\mathbf{v}_{p+1} = \mathbf{v}_p - 1$.

The same remark holds, with $\mathbf{v}'$ and $p_{i,j}'$ instead of $\mathbf{v}$ and $p_{i,j}$.

Let us introduce the property

$$P_k = \{\forall i \in \mathbb{N}^*, j \in \mathbb{Z}, \text{ s.t. } n_{i,j}(\mathbf{v}) \leq k \text{ and } n_{i,j}(\mathbf{v}') \leq k, E_{i,j} \text{ holds}\}.$$

Let us prove by induction on $k$ that $P_k$ holds for all $k \in \mathbb{N}$. Note that $P_0$ follows from $X_0 = X_0' = v_0$. Suppose $P_{k-1}$ holds. We want to deduce $P_k$, which is different from $P_{k-1}$ if there exists $(i,j) \in \mathbb{N}^* \times \mathbb{Z}$ such that $[n_{i,j}(\mathbf{v}) = k$ and $n_{i,j}(\mathbf{v}') \leq k]$ or $[n_{i,j}(\mathbf{v}') = k$ and $n_{i,j}(\mathbf{v}) \leq k]$. Select such a couple $(i,j)$.

If $i = 1$, then suppose, for instance, that $j > \mathbf{v}_0 = v_0$ (the case $j < v_0$ is analogous, and $j = v_0$ is obvious). Then $Z_{n_{i,j}(\mathbf{v}')}(j+1)(\mathbf{v}') = Z_{n_{i,j}(\mathbf{v})}(j+1)(\mathbf{v}) = 1$ and we aim to prove that $Z_{n_{i,j}(\mathbf{v}')}(j-1)(\mathbf{v}') \leq Z_{n_{i,j}(\mathbf{v})}(j-1)(\mathbf{v})$. Suppose the contrary $Z_{n_{i,j}(\mathbf{v}')}(j-1)(\mathbf{v}') > Z_{n_{i,j}(\mathbf{v})}(j-1)(\mathbf{v})$. Let $a = Z_{n_{i,j}(\mathbf{v})}(j-1)(\mathbf{v}) - 1$. Since $P_{k-1}$ holds, $E_{a,j-1}$ holds and, therefore, $p_{a,j-1}' \geq p_{a,j-1}$. Now $\mathbf{v}_{n_{a,j-1}+1} = \mathbf{v}_{n_{a,j-1}} + 1$ since $a = Z_{n_{i,j}(\mathbf{v})}(j-1)(\mathbf{v}) - 1$; this implies $\mathbf{v}_{n_{a,j-1}+1}' = \mathbf{v}_{n_{a,j-1}}' + 1$ and leads to a contradiction.

If $i > 1$, take $\mathbf{v}$ at time $n_{i-1,j}(\mathbf{v})$ and take $\mathbf{v}'$ at time $n_{i-1,j}(\mathbf{v}')$. We make use of the notation $\mathcal{M}$ or $\mathcal{M}' \to l$ or $r$ to indicate that $\mathbf{v}$ (resp. $\mathbf{v}'$) goes to the left or to the right at this time $n_{i-1,j}(\mathbf{v})$ [resp. $n_{i-1,j}(\mathbf{v}')$]. Since $E_{i-1,j}$ is satisfied, $p_{i-1,j}' \geq p_{i-1,j}$ and it is impossible that $\mathcal{M} \to r$ and $\mathcal{M}' \to l$. Hence, there are three cases:

- $\mathcal{M} \to l$ and $\mathcal{M}' \to r$, and the conclusion follows;

- $\mathcal{M} \to r$ and $\mathcal{M}' \to r$, and the conclusion follows from a proof similar to the case $i = 1$.
- $\mathcal{M} \to l$ and $\mathcal{M}' \to l$, and the conclusion follows from an analogous argument.

$\square$

DEFINITION 4.10. Let $\mathcal{M} = (X_n)_{n \in \mathbb{N}}$ be the VRRW on $\mathbb{Z}$ [on $(\Omega, \mathcal{F}, \mathbb{P})$ according to Definition 4.6] defined in the Introduction, that is, defined by $X_0 := v_0$ and the transition probabilities, for all $\mathbf{v} \in \mathcal{P}$ and $n \in \mathbb{N}$,

$$q_{\mathcal{M}}(\mathbf{v}, n) := \alpha_n^+(\mathbf{v}_n)(\mathbf{v}) = 1 - \alpha_n^-(\mathbf{v}_n)(\mathbf{v}).$$

DEFINITION 4.11. For all $x \in \mathbb{Z}$, $k \in \mathbb{N}$ and $M > 1$, let us define the $\mathbb{T}$ stopping time $U_{x,k,M}$ by, for all $\mathbf{v} \in \mathcal{P}$,

$$U_{x,k,M} := \inf\{n \geq k \text{ s.t. } Z_n(x)(\mathbf{v}) > M Z_k(x)(\mathbf{v}) \text{ or } \alpha_n^-(x)(\mathbf{v}) > M \alpha_k^-(x)(\mathbf{v})\}.$$

DEFINITION 4.12. Let $x \in \mathbb{Z}$ and $M \in \mathbb{R}_+^*$. Let us define, for all $\mathbb{T}$ stopping times $k$ and $V$, and $g \in \mathbb{R}_+^*$, the modified VRRW $\mathcal{M}'_{k,V,g}$ [on $(\Omega, \mathcal{F}, \mathbb{P})$ according to Definition 4.6] by $X_0' := v_0$ and the transition probabilities, for all $\mathbf{v} \in \mathcal{P}$ and $n \in \mathbb{N}$,

$$q_{\mathcal{M}'_{k,V,g}}(\mathbf{v}, n)$$
$$:= 1 - \alpha_n^-(\mathbf{v}_n)(\mathbf{v})(1 - \gamma_k(\mathbf{v}) \mathbb{1}_{\{\mathbf{v}_n = x \text{ and } n \in [k, U_{x,k,M} \wedge V]\} \cap \{\gamma_k(\mathbf{v}) < 1, \alpha_k^-(x)(\mathbf{v}) < 1/4M\}}),$$

where $\gamma_k(\mathbf{v}) := g / \sqrt{Z_k(x)(\mathbf{v}) \alpha_k^-(x)(\mathbf{v})}$.

Observe that, for all $k \in \mathbb{N}$ and $g \in \mathbb{R}_+^*$, $\mathcal{M}'_{k,V,g} \gg \mathcal{M}$ by Lemma 4.1.

4.3. *Nonconvergence result.* Let $c$, $d$, $M \in \mathbb{R}_+^*$ and $x \in \mathbb{Z}$. Let $(t_n)_{n \in \mathbb{N}}$ be an increasing sequence of $\mathbb{T}$ stopping times. Let $(y_n)_{n \in \mathbb{N}}$, $(z_n)_{n \in \mathbb{N}}$ and $(W(n))_{n \in \mathbb{N}}$ be $(\mathcal{T}_{t_n})_{n \in \mathbb{N}}$-adapted random processes taking values, respectively, in $\mathbb{R}$, in $\mathbb{R}$ and in $\mathbb{R}_+^*$, and let $T$ be a $(\mathcal{T}_{t_n})_{n \in \mathbb{N}}$ stopping time such that $n < T$ implies $t_n < \infty$. Let $\mathcal{M}$ be the random walk in Definition 4.10. By a slight abuse of notation, for all $n \in \mathbb{N}$, we let $t_n := t_n^{\mathcal{M}}$, $T := T^{\mathcal{M}}$, $y_n := y_n^{\mathcal{M}}$, $z_n := z_n^{\mathcal{M}}$ and $\mathcal{F}_n := \mathcal{F}_n^{\mathcal{M}}$.

ASSUMPTION (H1). There exists $n_0 \in \mathbb{N}$ (deterministic) and $(\mathcal{F}_{t_n})$-adapted processes $(\varepsilon_n)_{n \in \mathbb{N}^*}$ and $(r_n)_{n \in \mathbb{N}}$ such that, for all $n \geq n_0$,

$$(25) \qquad z_{n+1} - z_n = y_n + \varepsilon_{n+1} + r_n \qquad \text{if } n + 1 < T$$

and, if $n < T$,

$$\mathbb{E}(\varepsilon_{n+1}|\mathcal{F}_{t_n}) = 0, \qquad \mathbb{E}\left(\sum_{j=n}^{\infty}\varepsilon_{j+1}^2\Big|\mathcal{F}_{t_n}\right) \leq d^2 W(n)^2,$$

(26)

$$\sum_{j=n}^{\infty}|r_j| \leq dW(n).$$

Let us use again the notation of Definition 4.12 ($x \in \mathbb{Z}$ and $M \in \mathbb{R}_+^*$ are already fixed). For all $k$, $n \in \mathbb{N}$ and $g \in \mathbb{R}_+^*$, let $t'_{n,k,g} := t_n^{\mathcal{M}'_{t_k,t_{2k},g}}$, $T' := T^{\mathcal{M}'_{t_k,t_{2k},g}}$, $y'_{n,k,g} := y_n^{\mathcal{M}'_{t_k,t_{2k},g}}$, $z'_{n,k,g} := z_n^{\mathcal{M}'_{t_k,t_{2k},g}}$ and $\mathcal{F}'_n := \mathcal{F}_n^{\mathcal{M}'_{t_k,t_{2k},g}}$. Note that $\mathcal{F}'_{t'_{k,k,g}} = \mathcal{F}_{t_k}$, since the random walks $\mathcal{M}$ and $\mathcal{M}'_{t_k,t_{2k},g}$ are the same up to time $t_k$.

ASSUMPTION (H2). For all $g > 0$, there exists $n_0$ (deterministic, but dependent on $g$) and $(\mathcal{F}'_{t'_{n,k,g}})$-adapted processes $(\varepsilon'_n)_{n\in\mathbb{N}^*}$ and $(r'_n)_{n\in\mathbb{N}}$ such that, for all $k \geq n_0$ and $n \geq n_0$,

$$z'_{n+1,k,g} - z'_{n,k,g} \geq y'_{n,k,g} + \frac{cgW(k)}{k}\mathbb{1}_{n\in[k,2k)} + \varepsilon'_{n+1} + r'_n \qquad \text{if } n+1 < T'$$

(27)

and, if $n < T'$,

$$\mathbb{E}(\varepsilon'_{n+1}|\mathcal{F}'_{t'_{n,k,g}}) = 0, \qquad \mathbb{E}\left(\sum_{j=n}^{\infty}\varepsilon'^2_{j+1}\Big|\mathcal{F}'_{t'_{n,k,g}}\right) \leq d^2 W(n)^2,$$

(28)

$$\sum_{j=n}^{\infty}|r'_j| \leq dW(n).$$

Note that $T'$, $\mathcal{F}'_n$, $\varepsilon'_n$ and $r'_n$ also depend on $k$ and $g$, but we omit it for simplicity.

ASSUMPTION (H3). For all $g \in \mathbb{R}_+^*$, there exists $n_0 \in \mathbb{N}$ such that, for all $n \geq k \geq n_0$,

(29) $$y'_{n,k,g} - y_n \geq -d|y_n|[|z'_{n,k,g} - z_n| + W(k)].$$

PROPOSITION 4.1. *Let $c$, $d$, $M \in \mathbb{R}_+^*$ and let $x \in \mathbb{Z}$. Let $(t_n)_{n\in\mathbb{N}}$ be an increasing sequence of $\mathbb{T}$ stopping times. Let $(y_n)_{n\in\mathbb{N}}$, $(z_n)_{n\in\mathbb{N}}$ and $(W(n))_{n\in\mathbb{N}}$ be $(\mathcal{T}_{t_n})_{n\in\mathbb{N}}$-adapted random processes that take values, respectively, in $\mathbb{R}$, in* $\mathbb{R}$ and in $\mathbb{R}_+^*$, and let $T$ be a $(\mathcal{T}_{t_n})_{n \in \mathbb{N}}$ stopping time such that $n < T$ implies $t_n < \infty$. Suppose Assumptions (H1)–(H3) hold. Then

$$(30) \qquad \mathbb{P}\left[\left\{\lim_{n \to \infty} z_n = 0\right\} \cap \left\{\sum_{n < T} |y_n| < \infty\right\} \cap \{T = \infty\}\right] = 0.$$

NOTATION.   In the remainder of the paper (except Lemma 4.2 and its proof), we let $\mathcal{M}' := \mathcal{M}'_{t_k, t_{2k}, g}$, forgetting the dependence on $k$ and $g$. For any $\mathbb{T}$-measurable random variable $u$, we write $u$ (resp. $u'$) instead of $u^{\mathcal{M}}$ (resp. $u^{\mathcal{M}'}$). In particular, we use notation $\alpha_n^-(x)$ [resp. $\alpha_n'^-(x)$] for $\alpha_n^-(x)^{\mathcal{M}}$ [resp. $\alpha_n^-(x)^{\mathcal{M}'}$], $Z_n^+(x)$ [resp. $Z_n'^+(x)$] for $Z_n^+(x)^{\mathcal{M}}$ [resp. $Z_n^+(x)^{\mathcal{M}'}$] and so on. The notation $t'_{n,k,g}$, $y'_{n,k,g}$ and $z'_{n,k,g}$ is used in this current section to emphasize the link with variables $k$ and $g$ [especially in (27)], but is replaced subsequently by $t'_n$, $y'_n$ and $z'_n$.

REMARK 4.1.   Assume that $y_n$ can be written, when $n < T$, as

$$y_n = \frac{R_{t_n}}{S_{t_n}}$$

[where $(R_n)_{n \in \mathbb{N}}$ and $(S_n)_{n \in \mathbb{N}}$ are $\mathbb{T}$-adapted processes, taking values, respectively, in $\mathbb{R}$ and $\mathbb{R}_+^*$] and that, when $n < T \wedge T'$,

$$R'_{t'_n} \geq R_{t_n}.$$

Then a sufficient condition for Assumption (H3) is that, when $n \in [n_0, T \wedge T')$,

$$\left|\frac{S_{t_n}}{S'_{t'_n}} - 1\right| \leq d[|z'_n - z_n| + W(m)].$$

Indeed,

$$y'_n = \frac{R'_{t'_n}}{S'_{t'_n}} \geq \frac{R_{t_n}}{S'_{t'_n}} = y_n \frac{S_{t_n}}{S'_{t'_n}},$$

which implies

$$y'_n - y_n \geq y_n\left(\frac{S_{t_n}}{S'_{t'_n}} - 1\right).$$

REMARK 4.2.   In the conclusion of Proposition 4.1 [equation (30)], $\sum |y_n| < \infty$ can be replaced by $\sum y_n^\pm < \infty$ or by $\sum |y_n| \mathbb{1}_{y_n z_n \leq 0} < \infty$. Indeed, assume (H1) holds and $z_n \to 0$. Let us prove that $\sum y_n^+ < \infty$ implies $\sum y_n^- < \infty$ (and thus $\sum |y_n| < \infty$); the proof of the converse implication is very analogous. Summing (25) from $n := m$ to $p$, if $m \geq n_0$,

$$z_p - z_m = \sum_{k=m}^{p-1} (y_k^+ - y_k^-) + \sum_{k=m}^{p-1} (\varepsilon_{k+1} + r_k)$$

and, therefore,

$$\sum_{k=m}^{\infty} y_k^- \leq z_m + \sum_{k=m}^{\infty} y_k^+ + \sup_{p \geq m} \left| \sum_{k=m}^{p-1} \varepsilon_{k+1} \right| + \sum_{k=m}^{\infty} |r_k| < \infty$$

using (26) and Doob's inequality $[W(m) < \infty]$.

Similarly, $\sum |y_n| \mathbb{1}_{y_n z_n \leq 0} < \infty$ implies $\sum |y_n| \mathbb{1}_{y_n z_n \geq 0} < \infty$ (and thus $\sum |y_n| < \infty$), using

$$|z_{n+1}| - |z_n| \geq |y_n| \mathbb{1}_{y_n z_n \geq 0} - |y_n| \mathbb{1}_{y_n z_n \leq 0} + \operatorname{sign}(z_n)(\varepsilon_{n+1} + r_n),$$

where, for all $x \in \mathbb{R}$, $\operatorname{sign}(x)$ denotes the sign of $x$.

4.4. *Proof of Proposition* 4.1. The following lemma is useful in the proof of Proposition 4.1.

LEMMA 4.2. *Let $x \in \mathbb{Z}$ and $M \in \mathbb{R}_+^*$. For all $\mathbb{T}$ stopping times $k$ and $V$, and $g \in \mathbb{R}_+^*$, let $\mathcal{M}$ and $\mathcal{M}'_{k,V,g}$ be the two random walks on $(\Omega, \mathcal{F}, \mathbb{P})$ introduced in Definitions* 4.6, 4.10 *and* 4.12. *Then, for all $\eta > 0$ and $\mathcal{C} \in \mathcal{T}$,*

$$\mathbb{P}_{\mathcal{M}'_{k,V,g}}(\mathcal{C}|\mathcal{T}_k) \leq \operatorname{Cst}(g, M, \eta) \mathbb{P}_{\mathcal{M}}(\mathcal{C}|\mathcal{T}_k) + \eta.$$

PROOF. For simplicity, let $\mathcal{M}' := \mathcal{M}'_{k,V,g} = (X'_n)_{n \in \mathbb{N}}$ and, for any $\mathbb{T}$-measurable r.v. $u$, $u := u^{\mathcal{M}}$ and $u' := u^{\mathcal{M}'}$. Given $n \in \mathbb{N}$ and $\mathbf{v} \in \mathcal{P}$, let $I_{\mathbf{v}}(n) := \mathbb{P}_{\mathcal{M}'_{k,V,g}}(\mathcal{C}_{(\mathbf{v}_0,\ldots,\mathbf{v}_n)}) / \mathbb{P}_{\mathcal{M}}(\mathcal{C}_{(\mathbf{v}_0,\ldots,\mathbf{v}_n)})$. If $\gamma_k(\mathbf{v}) \geq 1$ or $\alpha_k^-(x)(\mathbf{v}) \geq 1/4M$, then $I_{\mathbf{v}}(n) = 1$. Otherwise,

$$I_{\mathbf{v}}(n) = \prod_{j=k+1}^{n \wedge V \wedge U_{x,k,M}(\mathbf{v})} \left[ \left(1 - \gamma_k(\mathbf{v}) \mathbb{1}_{\{\mathbf{v}_{j-1} = x, \mathbf{v}_j = x-1\}}\right) \right.$$

$$\left. \times \left(1 + \frac{\gamma_k(\mathbf{v}) \mathbb{1}_{\{\mathbf{v}_{j-1} = x, \mathbf{v}_j = x+1\}} \alpha_{j-1}^-(x)(\mathbf{v})}{\alpha_{j-1}^+(x)(\mathbf{v})}\right) \right].$$

Let us upper bound $\ln I_{\mathbf{v}}(n)$, using $\ln(1+x) \leq x$ for all $x > -1$. We obtain

$$\ln I_{\mathbf{v}}(n) \leq R_n(\mathbf{v}),$$

where $(R_n)_{n\geq k}$ is the $(\mathcal{T}_n)_{n\geq k}$-adapted process defined by $R_k = 0$ and, for $n \geq k+1$,

$$R_n(\mathbf{v}) = \gamma_k(\mathbf{v}) \sum_{j=k+1}^{n\wedge V \wedge U_{x,k,M}} \mathbb{1}_{\{\mathbf{v}_{j-1}=x\}} \left[ -\mathbb{1}_{\{\mathbf{v}_j = x-1\}} + \frac{\mathbb{1}_{\{\mathbf{v}_j = x+1\}}\alpha_{j-1}^{-}(x)(\mathbf{v})}{\alpha_{j-1}^{+}(x)(\mathbf{v})} \right].$$

Our goal is to overestimate $\sup_n R_n(\mathbf{v})$ $[\geq \sup_n \ln I_{\mathbf{v}}(n)]$ on a $\mathcal{T}$-measurable subset of $\mathcal{P}$ of large probability for $\mathcal{M}'$. To this end, we now analyze the behavior of the $(\mathcal{F}'_n)_{n\in\mathbb{N}}$-measurable process $R'_n := R_n^{\mathcal{M}'}$, depending on the random walk $\mathcal{M}'$ (recall that $\mathcal{F}'_n = \mathcal{F}_n^{\mathcal{M}'}$).

Let $(\widetilde{R}'_n)_{n\geq k}$ be the compensator of $(R_n)_{n\geq k}$, that is, the $(\mathcal{F}'_n)_{n\geq k}$-predictable process such that the process $\widehat{R}'_n := R'_n - \widetilde{R}'_n$ is a martingale. For all $n \geq k$,

$$\widetilde{R}'_{n+1} - \widetilde{R}'_n$$
$$= \gamma_k \mathbb{1}_{\{X'_n = x\}} \mathbb{1}_{\{n < V' \wedge U'_{x,k,M}\}} \left[ -\alpha_n'^{-}(x)(1-\gamma_k) + \frac{\alpha_n'^{+}(x) + \gamma_k \alpha_n'^{-}(x)}{\alpha_n'^{+}(x)} \alpha_n'^{-}(x) \right]$$
$$= \gamma_k^2 \mathbb{1}_{\{X'_n = x\}} \mathbb{1}_{\{n < V' \wedge U'_{x,k,M}\}} \alpha_n'^{-}(x) \left[ 1 + \frac{\alpha_n'^{-}(x)}{\alpha_n'^{+}(x)} \right].$$

Therefore, for all $n \geq k$, using $\alpha_j'^{-}(x) \leq M\alpha_k^{-}(x) \leq 1/4$ and $Z'_j(x) \leq MZ_k(x)$ for all $j < U'_{x,k,M}$,

$$(31) \qquad \widetilde{R}'_n \leq 2\gamma_k^2(M-1)Z_k(x)M\alpha_k^{-}(x) = 2g^2 M(M-1).$$

Now, for all $n \geq k$,

$$\mathbb{V}[\widehat{R}'_{n+1}|\mathcal{F}'_n] \leq \mathbb{E}[(R'_{n+1} - R'_n)^2|\mathcal{F}'_n] \leq \gamma_k^2 \mathbb{1}_{\{X'_n = x, n < U'_{x,k,M}\}} 2\alpha_n'^{-}(x)$$
$$\leq \mathbb{1}_{\{X'_n = x, n < U'_{x,k,M}\}} 2M\alpha_k^{-}(x)\gamma_k^2,$$

where $\mathbb{V}[\cdot|\mathcal{F}_n]$ is the variance conditional on $\mathcal{F}_n$.

Successively using the Bienaymé–Tchebycheff and Doob inequalities, for all $A \in \mathbb{R}_+^*$,

$$\mathbb{P}\left(\sup_{n\geq k}\widehat{R}'_n \geq A \Big| \mathcal{F}'_k\right) \leq \mathbb{E}\left(\sup_{n\geq k}\widehat{R}_n'^2 \Big| \mathcal{F}_k\right)\Big/A^2 \leq 4\mathbb{E}[\widehat{R}_{V'\wedge U'_{x,k,M}}'^2|\mathcal{F}_k]/A^2$$
$$\leq 8M\alpha_k^{-}(x)\gamma_k^2(M-1)Z_k(x)/A^2 = 8M(M-1)g^2/A^2.$$

Hence, by letting $\Gamma$ be the $\mathcal{T}$-measurable event

$$\Gamma := \left\{ \sup_{n\geq k} R_n(\mathbf{v}) \geq \sqrt{8M(M-1)g^2/\eta} + 2g^2 M(M-1) \right\},$$

we deduce, using $R'_n = \widetilde{R}'_n + \widehat{R}'_n$ and (31), that

$$\mathbb{P}_{\mathcal{M}'}(\Gamma|\mathcal{T}_k) \leq \eta.$$

On the other hand, for any path $\mathbf{v} \in \mathcal{P} \cap \Gamma^c$, for all $n \in \mathbb{N}$,

$$\mathbb{P}_{\mathcal{M}'}(\mathcal{C}_{(\mathbf{v}_0,\dots,\mathbf{v}_n)}|\mathcal{T}_k) \leq \operatorname{Cst}(g,M,\eta)\mathbb{P}_{\mathcal{M}}(\mathcal{C}_{(\mathbf{v}_0,\dots,\mathbf{v}_n)}|\mathcal{T}_k),$$

which enables us to conclude. $\square$

Let us now prove Proposition 4.1. Given $\varepsilon > 0$, which is chosen subsequently, let us define the stopping time $Z^m = \inf\{n \geq m/\sum_{i=m}^n |y_i| > \varepsilon\}$. Let, for all $m \in \mathbb{N}$,

$$\Gamma_m := \{\lim z_n \neq 0\} \cup \{T \wedge Z^m < \infty\}, \qquad \mathcal{C}_m := \mathcal{I}_{\mathcal{M}}^{-1}(\Gamma_m).$$

It suffices here to prove that there exists $g = \operatorname{Cst}(c,d,M)$ such that if $k \geq m \vee n_0$, then $\mathbb{P}[\mathcal{C}_m|\mathcal{F}_{t_k}] \geq \operatorname{Cst}(c,d,M)$. Indeed, by a standard martingale convergence theorem, $\mathbb{P}[\mathcal{C}_m|\mathcal{F}_{t_k}] = \mathbb{E}[\mathbb{1}_{\mathcal{C}_m}|\mathcal{F}_{t_k}] \to_{k\to\infty} \mathbb{1}_{\mathcal{C}_m}$ since $\mathcal{C}_m \in \mathcal{F}_{t_\infty} = \sigma(\bigcup \mathcal{F}_{t_n})$. Therefore, $\mathbb{1}_{\mathcal{C}_m} \geq \operatorname{Cst}(c,d,M) > 0$ a.s. and $\mathbb{P}(\mathcal{C}_m) = 1$ for all $m \in \mathbb{N}$. We conclude that $\mathbb{P}(\{\lim z_n \neq 0\} \cup \{\sum |y_n| = \infty\} \cup \{T < \infty\}) = \mathbb{P}(\bigcap_{m\in\mathbb{N}} \mathcal{C}_m) = 1$.

Let $g \in \mathbb{R}_+^*$ be fixed later. Let us introduce the random walk $\mathcal{M}'_{t_k,t_{2k},g}$ (see Definition 4.12) and use the notation introduced in Section 4.3. Let us consider, for all $k \in \mathbb{N}$,

$$\Delta = \left\{\sup_{n \geq k}\left|\sum_{i=k}^n \varepsilon_{i+1}\right| \leq 4\,dW(k), \sup_{n \geq k}\left|\sum_{i=k}^n \varepsilon'_{i+1}\right| \leq 4\,dW(k)\right\}.$$

Let us apply the Bienaymé–Tchebycheff and Doob inequalities, and use Assumptions (H1) and (H2): For all $k \geq n_0$,

$$\mathbb{P}(\Delta^c|\mathcal{F}_{t_k}) \leq \frac{\mathbb{E}(\sup_{n\geq k}(\sum_{i=k}^n \varepsilon_{i+1})^2|\mathcal{F}_{t_k})}{16\,d^2 W(k)^2} + \frac{\mathbb{E}(\sup_{n\geq k}(\sum_{i=k}^n \varepsilon'_{i+1})^2|\mathcal{F}_{t_k})}{16\,d^2 W(k)^2}$$
$$\leq \frac{8}{16} = \frac{1}{2}.$$

We want to prove that on $\Delta$, $T \wedge T' \wedge Z^m < \infty$ or $z_n \not\to 0$. From now on, we suppose that $\Delta$ holds and that $T \wedge T' \wedge Z^m = \infty$.

Observe that Assumptions (H1)–(H3) imply, for all $i \geq k$ and $n \geq i$,

$$(32) \qquad z'_n - z_n \geq z'_i - z_i - d\varepsilon \sup_{j\in[i,n-1]} |z'_j - z_j| - d(10+\varepsilon)W(k)$$

and that for all $n \geq 2k$ (using $z'_k = z_k$, the coupled random walks $\mathcal{M}$ and $\mathcal{M}'$ being identical up to time $t_k$),

$$(33) \qquad z'_n - z_n \geq [cg - d(10+\varepsilon)]W(k) - d\varepsilon \sup_{j\in[k,n]} |z'_j - z_j|.$$

Suppose $n \geq 2k$ and let

$$u(n) := \sup\left\{i \in [k,n] \text{ s.t. } |z'_i - z_i| = \sup_{j\in[k,n]} |z'_j - z_j|\right\}, \qquad \tau_n = z'_{u(n)} - z_{u(n)}.$$

On one hand, using (33),

$$|\tau_n| \geq \tau_n \geq [cg - d(10 + \varepsilon)]W(k) - d\varepsilon|\tau_n|.$$

We take $g := 24d/c$ and suppose $\varepsilon < \mathrm{Cst}(d)$; hence $|\tau_n| > 12\,dW(k)$. On the other hand, using (32) with $i := k$ and $n := u(n)$,

$$\tau_n \geq -d\varepsilon|\tau_n| - d(10 + \varepsilon)W(k).$$

Therefore $\tau_n < 0$ implies $|\tau_n|(1 - d\varepsilon) \leq d(10 + \varepsilon)W(k)$ and hence $|\tau_n| \leq 12\,dW(k)$ [if we suppose $\varepsilon < \mathrm{Cst}(d)$], which leads to a contradiction. Therefore $\tau_n \geq 12\,dW(k)$. For all $n \geq 2k$, (32) with $i := u(n)$ implies

$$z_n' - z_n \geq \tau_n(1 - d\varepsilon) - d(10 + \varepsilon)W(k) \geq dW(k) > 0$$

if we assume $\varepsilon < \mathrm{Cst}(d)$. Therefore,

$$\begin{aligned}
&\mathbb{P}_{\mathcal{M}}(\Gamma_m|\mathcal{T}_{t_k}) + \mathbb{P}_{\mathcal{M}'}(\Gamma_m|\mathcal{T}_{t_k}) \\
&= (\mathbb{P}(\mathcal{I}_{\mathcal{M}}^{-1}(\Gamma_m)|\mathcal{F}_{t_k}) + \mathbb{P}(\mathcal{I}_{\mathcal{M}'}^{-1}(\Gamma_m)|\mathcal{F}_{t_k})) \circ \mathcal{I}_{\mathcal{M}}^{-1} \\
&\geq \mathbb{P}\left(\lim_{n\to\infty} z_n \neq 0 \text{ or } \lim_{n\to\infty} z_n' \neq 0 \text{ or } T \wedge T' \wedge W < \infty|\mathcal{F}_{t_k}\right) \circ \mathcal{I}_{\mathcal{M}}^{-1} \\
&\geq \mathbb{P}(\Delta|\mathcal{F}_{t_k}) \circ \mathcal{I}_{\mathcal{M}}^{-1} \geq 1/2.
\end{aligned}$$

In the equality, we use that $\mathcal{F}_{t_k} = \mathcal{F}_{t_k'}'$ and, for all $\mathcal{C} \in \mathcal{P}$ and $n \in \mathbb{N}$, $\mathbb{P}_{\mathcal{M}}(\mathcal{C}|\mathcal{T}_{t_n}) = \mathbb{P}(\mathcal{I}_{\mathcal{M}}^{-1}(\mathcal{C})|\mathcal{F}_{t_n}) \circ \mathcal{I}_{\mathcal{M}}^{-1}$ [for all $\mathbf{v} \in \mathcal{P}$, $\mathbb{P}(\mathcal{I}_{\mathcal{M}}^{-1}(\mathcal{C})|\mathcal{F}_{t_n})$ is constant on $\mathcal{I}_{\mathcal{M}}^{-1}(\mathbf{v})$]; similarly $\mathbb{P}_{\mathcal{M}'}(\mathcal{C}|\mathcal{T}_{t_n}) = \mathbb{P}(\mathcal{I}_{\mathcal{M}'}^{-1}(\mathcal{C})|\mathcal{F}_{t_n'}) \circ \mathcal{I}_{\mathcal{M}'}^{-1}$. We also use that for all $\mathbf{v} \in \mathcal{P}$, $\mathcal{I}_{\mathcal{M}}^{-1}(\mathbf{v})$ and $\mathcal{I}_{\mathcal{M}'}^{-1}(\mathbf{v})$ have the same projections on the $t_k$ first coordinates (the coupled random walks $\mathcal{M}$ and $\mathcal{M}'$ are the same up to time $t_k$).

Now, Lemma 4.2 implies that for all $\eta > 0$,

$$\mathbb{P}_{\mathcal{M}'}(\Gamma_m|\mathcal{T}_{t_k}) \leq \mathrm{Cst}(g, M, \eta)\mathbb{P}_{\mathcal{M}}(\Gamma_m|\mathcal{T}_{t_k}) + \eta$$

and, therefore,

$$\mathbb{P}(\mathcal{C}_m|\mathcal{F}_{t_k}) = \mathbb{P}_{\mathcal{M}}(\Gamma_m|\mathcal{T}_{t_k}) \circ \mathcal{I}_{\mathcal{M}} \geq \frac{1/2 - \eta}{1 + \mathrm{Cst}(g, M, \eta)} = \mathrm{Cst}(c, d, M)$$

if we take $\eta = 1/4$ (recall that we have chosen $g = 24d/c$).

## 5. Proofs of Lemmas 2.5, 2.7 and 2.8.

### 5.1. *Proof of Lemma 2.5*.

5.1.1. *Notation.* It suffices to prove that $\mathbb{P}(\Upsilon_0(0) \cap \Upsilon_0'(0)) = 0$, because the problem is translation-invariant. We apply Proposition 4.1 to show it (as explained in Section 4.1) and use its notation. Let us first introduce our choice of the variables that appear in this lemma.

Let us define a sequence $(t_n)_{n \in \mathbb{N}}$ by

$$t_n := \inf\{m \in \mathbb{N} \text{ s.t. } Z_m^+(2) \geq n\}.$$

Let, for all $n \in \mathbb{N}$,

$$R_n = Z_n(2) + Z_n(4) - (Z_n(1) + Z_n(3)).$$

Let $(y_n)_{n \in \mathbb{N}}$ and $(z_n)_{n \in \mathbb{N}}$ be the $(\mathcal{T}_{t_n})$-adapted processes defined by

$$y_n := \frac{R_{t_n}}{Z_{t_n}(2) Z_{t_n}(3)}, \qquad z_n := \ln \frac{Z_{t_n}(3)}{Z_{t_n}(2)}$$

if $t_n < \infty$, and $y_n = z_n := 0$ else.

Given $\varepsilon > 0$, for all $m \in \mathbb{N}$, let $T_1^{m,\varepsilon}$ and $T_2^{m,\varepsilon}$ be the $\mathbb{F}_{(t_n)_{n \in \mathbb{N}}}$ stopping times defined by

$$T_1^{m,\varepsilon} := \inf\{n \geq m \text{ s.t. } Z_{t_n}(4) \geq e(1+\varepsilon) Z_{t_n}(1) \text{ or } \alpha_{t_n}^-(2) \geq (1+\varepsilon) \inf_{m \leq j \leq n} \alpha_{t_j}^-(2)$$

$$\text{or } Z_{t_n}(1) \wedge Z_{t_n}(4) \notin [(Z_{t_n}(2) \vee Z_{t_n}(3))^{1-\varepsilon}, \varepsilon Z_{t_n}(2) \wedge Z_{t_n}(3))]\},$$

$$T_2^{m,\varepsilon} := \inf\{n \geq m \text{ s.t. } t_n = \infty \text{ or } \exists y \in [1,4] \text{ s.t. } Z_{t_n}(y) - Z_{t_{n-1}}(y) \geq Z_{t_{n-1}}(y)^\varepsilon\}$$

and $T^{m,\varepsilon} := T_1^{m,\varepsilon} \wedge T_2^{m,\varepsilon}$.

LEMMA 5.1. *For all $\varepsilon > 0$,*

$$\Upsilon_0(0) \cap \Upsilon_0'(0) \subset \{\lim z_n = 0\} \cap \left\{\sum y_n^- < \infty\right\} \cap \left(\bigcup_{m \in \mathbb{N}} \{T^{m,\varepsilon} = \infty\}\right).$$

PROOF. Observe that a.s. on $\Upsilon_0(0) \cap \Upsilon_0'(0)$,

$$\sum_{n \in \mathbb{N}} y_n^- \preceq \sum_{n \in \mathbb{N}} \frac{Z_{t_n}(0)}{Z_{t_n}(2) Z_{t_n}(3)}$$

$$\approx \sum_{n \in \mathbb{N}} \frac{Z_{t_n}(0)}{Z_{t_n}(0) + Z_{t_n}(2)} \frac{Z_{t_n}(1)}{Z_{t_n}(1) + Z_{t_n}(3)} \frac{1}{Z_{t_n}(1)}$$

$$\preceq \sum_{n \in \mathbb{N}} \alpha_n^-(1) \frac{\mathbb{1}_{\{X_n = 2, X_{n+1} = 1\}}}{Z_n(1)}$$

$$\preceq \sum_{n \in \mathbb{N}} \alpha_n^-(1) \frac{\mathbb{1}_{\{X_n = 1\}}}{Z_n(1)} \approx \sum_{n \in \mathbb{N}} \frac{\mathbb{1}_{\{X_n = 1, X_{n+1} = 0\}}}{Z_n(1)}$$

$$\preceq \sum_{n \in \mathbb{N}} \frac{\mathbb{1}_{\{X_n = 0, X_{n+1} = 1\}}}{Z_n(1)} < \infty,$$

where we use (23) in the first relationship, the definition of $\Upsilon_0(0) \cap \Upsilon'_0(0)$ in the second and seventh relationships, Lemma A.1(i) in the third and fifth relationships, Proposition 3.1(c) in the sixth relationship and, finally, $\Upsilon(0) = \{Y^+_\infty(0) < \infty\}$ by Lemma 2.1.

The fact that there exists a.s. on $\Upsilon_0(0) \cap \Upsilon'_0(0)$, $m \in \mathbb{N}$ such that $T_2^{m,\varepsilon} = \infty$ can be proved as follows: For instance, for $Z_{t_n}(1) - Z_{t_{n-1}}(1)$, the probability to visit 3 starting from 1 in two steps is asymptotically greater than $1/4$, so we can use a method very similar to the proof of (17) to estimate the number of visits to 1 between times $t_{n-1}$ and $t_n$. The other points follow from the definitions.   □

Lemma 5.1 implies, together with Remark 4.2, that it suffices to apply Proposition 4.1, $m$ and $\varepsilon$ being fixed, with $(t_n)_{n\in\mathbb{N}}$, $(y_n)_{n\in\mathbb{N}}$ and $(z_n)_{n\in\mathbb{N}}$ as defined below, to conclude that $\mathbb{P}(\Upsilon_0(0) \cap \Upsilon'_0(0)) = 0$. Let us choose here the other constants that appear in the application of this lemma (the choice is justified afterward): $c := 1/64$, $d := 8\sqrt{1+e}$, $M := 4$, $x := 2$, $W(n) := \sqrt{\alpha^-_{t_n}(2)/n}$ and $T := T^{m,\varepsilon}$. We choose $\varepsilon$ in the following text. We check in Sections 5.1.2–5.1.4 that Assumptions (H1)–(H3) hold.

5.1.2. *Assumption* (H1) *of Proposition* 4.1 *holds.*   Let $n \in \mathbb{N}$ be such that $n < T$. We need to define a continuation of $z_{n+1}$ on $n+1 = T$, so that conditions (25) and (26) hold. Let us define

$$\widetilde{Z}_{t_{n+1}}(2) := \inf\{i \geq Z_{t_n}(2)+1 \text{ s.t. } \omega_{i-1,2} \leq \max(\alpha^+_{n_{i-1,2}}(2), \alpha^+_{t_n}(2) - Z_{t_n}(1)^{\varepsilon-1})\},$$

where $n_{i,j}$, $i \in \mathbb{N}^*$, $j \in \mathbb{Z}$, is the $i$th visit time to site $j$, as in Definition 4.1. We apply the convention that $\alpha^+_{n_{i-1,2}}(2) = 0$ whenever $n_{i-1,2} = \infty$.

Note that $\widetilde{Z}_{t_{n+1}}(2)$ is $\mathcal{F}_{t_{n+1}}$-measurable and that $\widetilde{Z}_{t_{n+1}}(2) = Z_{t_{n+1}}(2)$ when $n+1 < T$. Indeed, $Z_{t_n}(1) \leq Z_{t_n}(1) + Z_{t_n}(1)^\varepsilon$ implies, for all $k \in [t_n, t_{n+1}]$, $\alpha^-_{t_n}(2) \leq \alpha^+_{t_n}(2) + Z_{t_n}(1)^{\varepsilon-1}$ [hence $\alpha^+_k(2) \geq \alpha^+_{t_n}(2) - Z_{t_n}(1)^{\varepsilon-1}$].

Let $V_{n+1} := \widetilde{Z}_{t_{n+1}}(2) - Z_{t_n}(2)$. Then $V_{n+1}$ is lower and upper bounded by two geometric random variables with success probabilities $\alpha^+_{t_n}(2) - Z_{t_n}(1)^{\varepsilon-1}$ and $\alpha^+_{t_n}(2)$, and, therefore, assuming $\varepsilon < \text{Cst}$ and $n \geq \text{Cst}$ [note that $Z_{t_n}(2) \geq n$],

$$(34) \qquad \mathbb{E}(V_{n+1}|\mathcal{F}_{t_n}) \in [\alpha^+_{t_n}(2)^{-1}, (\alpha^+_{t_n}(2) - Z_{t_n}(1)^{\varepsilon-1})^{-1}]$$
$$\subset [\alpha^+_{t_n}(2)^{-1}, \alpha^+_{t_n}(2)^{-1} + 2Z_{t_n}(1)^{\varepsilon-1}]$$

and

$$(35) \qquad \mathbb{E}((V_{n+1}-1)^2|\mathcal{F}_{t_n}) \leq 3(\alpha^-_{t_n}(2) + Z_{t_n}(1)^{\varepsilon-1}) \leq 6\alpha^-_{t_n}(2).$$

Indeed, if $\zeta$ is a geometric r.v. with success probability $p$ and if $1 - p < \text{Cst}$,

$$\mathbb{E}(\zeta) = 1/p, \qquad \mathbb{E}((\zeta-1)^2) = (1-p)(2-p)/p^2 \leq 3(1-p).$$

Let us define $\widetilde{Z}_{t_{n+1}}(3)$ similarly so that if $W_{n+1} := \widetilde{Z}_{t_{n+1}}(3) - Z_{t_n}(3)$,

$$\mathbb{E}(W_{n+1}|\mathcal{F}_{t_n}) \in [\alpha_{t_n}^-(3)^{-1}, \alpha_{t_n}^-(3)^{-1} + 2Z_{t_n}(4)^{\varepsilon-1}]$$

and such that the estimate analogous to (35) holds for $W_{n+1}$. Let

$$\tilde{z}_{n+1} := \ln \frac{\widetilde{Z}_{t_{n+1}}(3)}{Z_{t_{n+1}}(2)}, \qquad r_n := \mathbb{E}(\tilde{z}_{n+1} - z_n|\mathcal{F}_{t_n}) - y_n,$$

$$\varepsilon_{n+1} := \tilde{z}_{n+1} - z_n - \mathbb{E}(\tilde{z}_{n+1} - z_n|\mathcal{F}_{t_n}).$$

Note that if $\varepsilon < \mathrm{Cst}$ and $n \geq \mathrm{Cst}(\varepsilon)$, using notation $\square$ in Section 2.1,

$$(36) \quad \begin{aligned} \mathbb{E}\left(\ln \frac{\widetilde{Z}_{t_{n+1}}(2)}{Z_{t_n}(2)}\Big|\mathcal{F}_{t_n}\right) &= \frac{\mathbb{E}(V_{n+1}|\mathcal{F}_{t_n})}{Z_{t_n}(2)} + \square\left(\frac{\mathbb{E}(V_{n+1}^2|\mathcal{F}_{t_n})}{2Z_{t_n}(2)^2}\right) \\ &= \frac{\alpha_{t_n}^+(2)^{-1}}{Z_{t_n}(2)} + \square\left(\frac{2Z_{t_n}(1)^{\varepsilon-1}}{n}\right) + \square\left(\frac{1}{n^2}\right) \\ &= \frac{\alpha_{t_n}^+(2)^{-1}}{Z_{t_n}(2)} + \square\left(\frac{1}{2n^{2-3\varepsilon}}\right), \end{aligned}$$

using $|\ln(1+x) - x| \leq x^2/2$ for $x \geq 0$ in the first equation, $\mathbb{E}(V_{n+1}^2|\mathcal{F}_{t_n}) \leq 2$, $Z_{t_n}(2) \geq n$ and (34) in the second equation, and $Z_{t_n}(1)^{\varepsilon-1} \leq n^{-(1-\varepsilon)^2} \leq 1/(4n^{1-3\varepsilon})$ in the third equation.

Using a similar estimate for $\mathbb{E}(\ln \widetilde{Z}_{t_{n+1}}(3)/Z_{t_n}(3)|\mathcal{F}_{t_n})$, we obtain

$$(37) \quad \mathbb{E}(\tilde{z}_{n+1} - z_n|\mathcal{F}_{t_n}) = y_n + \square(1/n^{2-3\varepsilon}).$$

Therefore, $r_n = \square(1/n^{2-3\varepsilon})$ is such that

$$\sum_{j=n}^{\infty} |r_j| = \square\left(\sum_{j=n}^{\infty} \frac{1}{j^{2-3\varepsilon}}\right) \leq \frac{1}{(1-3\varepsilon)(n-1)^{1-3\varepsilon}} \leq 2\sqrt{\frac{\alpha_{t_n}^-(2)}{n}} = 2W(n),$$

using $\alpha_{t_n}^-(2) \geq 1/n^{2\varepsilon}$, if $\varepsilon < \mathrm{Cst}$ and $n \geq \mathrm{Cst}(\varepsilon)$.

Let us now estimate $\mathbb{E}(\varepsilon_{n+1}^2|\mathcal{F}_{t_n})$:

$$\begin{aligned} &\mathbb{V}\left(\ln \frac{\widetilde{Z}_{t_{n+1}}(2)}{Z_{t_n}(2)}\Big|\mathcal{F}_{t_n}\right) \\ &\quad \leq \mathbb{E}\left(\left[\ln\left(1 + \frac{V_{n+1}}{Z_{t_n}(2)}\right) - \ln\left(1 + \frac{1}{Z_{t_n}(2)}\right)\right]^2\Big|\mathcal{F}_{t_n}\right) \\ &\quad \leq \mathbb{E}\left(\left[\ln\left(1 + \frac{V_{n+1}-1}{Z_{t_n}(2)+1}\right)\right]^2\Big|\mathcal{F}_{t_n}\right) \leq \frac{\mathbb{E}[(V_{n+1}-1)^2|\mathcal{F}_{t_n}]}{n^2} \leq \frac{6\alpha_{t_n}^-(2)}{n^2}, \end{aligned}$$

using (35) (note that the second inequality is an equality). Using a similar estimate of conditional variance of $\log \widetilde{Z}_{t_{n+1}}(2)/Z_{t_n}(2)$, we obtain, if $\varepsilon < \mathrm{Cst}$ and $n \geq k \geq \mathrm{Cst}(\varepsilon)$,

$$(38) \qquad \mathbb{E}(\varepsilon_{n+1}^2 | \mathcal{F}_{t_n}) \leq \frac{12(\alpha_{t_n}^-(2) + \alpha_{t_n}^+(3))}{n^2}$$

$$\leq \frac{24\alpha_{t_n}^-(2)(1+e)}{n^2} \leq \frac{48(1+\varepsilon)(1+e)\alpha_{t_k}^-(2)}{n^2},$$

using $\alpha_{t_n}^+(3) \leq 2e\alpha_{t_n}^-(2)$, since $n < T_1^{m,\varepsilon}$. Indeed, $Z_{t_n}(2) \in [(1+\varepsilon)^{-1}Z_{t_n}(3), (1+\varepsilon)Z_{t_n}(3)]$ $[Z_{t_n}(3) \leq Z_{t_n}(2) + Z_{t_n}(4) \leq (1+\varepsilon)Z_{t_n}(2)$ and the other inequality is similar].

We conclude from (37) and (38) that (26) holds and, therefore, that Assumption (H1) is satisfied when $\varepsilon < \mathrm{Cst}$ and $n \geq \mathrm{Cst}(\varepsilon)$.

5.1.3. *Assumption* (H2) *of Proposition* 4.1 *holds.* Note that $\gamma_{t_k} := g/\sqrt{Z_{t_k}(2)\alpha_{t_k}^-(2)} < 1$ if $k \geq \mathrm{Cst}(g)$ and $\varepsilon < \mathrm{Cst}$. Let us prove that $2k < T$ implies $t_{2k} < U_{2,t_k,M}$ so that in Definition 4.12 for the transition probabilities of $\mathcal{M}' := \mathcal{M}'_{t_k, t_{2k}, g}$ (in the statement of Proposition 4.1), $t_{2k} \wedge U_{2,t_k,M} = t_{2k}$. This implies that if $2k \leq T'$, the probabilities of moving from site 2 to site 1 for $\mathcal{M}'$ are the VRRW probabilities multiplied by the factor $1 - \gamma_{t_k}$ from time $t_k$ to time $t'_{2k}$. First, $\alpha_{t_k}^-(2) \leq (1+\varepsilon)\alpha_{t_k}^-(2) \leq 2\alpha_{t_k}^-(2)$ as long as $n < T$ by definition of $T$, if $\varepsilon \leq 1$. Second, if $n \leq 2k < T$, then $Z_{t_n}(2) \leq Z_{t_n}^+(2) + Z_{t_n}(1) \leq Z_{t_n}^+(2) + \varepsilon Z_{t_n}(2)$, so that $Z_{t_n}(2) \leq (1-\varepsilon)^{-1}Z_{t_n}^+(2) = (1-\varepsilon)^{-1}n \leq 4k$ if $\varepsilon \leq 1/2$.

Let $n \in \mathbb{N}$ be such that $n < T'$. We define

$$\widetilde{Z}'_{t'_{n+1}}(2) := \inf\{i \geq Z'_{t'_n}(2) + 1 \text{ s.t.}$$

$$1 - \omega_{i-1,2} \geq (1 - \gamma_{t_k}\mathbb{1}_{n \in [k,2k)})\min(\alpha_{n'_{i-1,2}}^{'-}(2), \alpha_{t'_n}^{'-}(2) + Z'_{t'_n}(1)^{\varepsilon-1})\},$$

with the convention that $\alpha_{n'_{i-1,2}}^{'-}(2) = 1$ whenever $n'_{i-1,2} = \infty$.

Similarly as in Section 5.1.2, $\widetilde{Z}'_{t'_{n+1}}(2)$ is $\mathcal{F}'_{t'_{n+1}}$-measurable and $\widetilde{Z}'_{t'_{n+1}}(2) = Z'_{t'_{n+1}}(2)$ when $n+1 < T'$. We define $\widetilde{Z}'_{t'_{n+1}}(3)$ in the same way (but the factor $1 - \gamma_{t_k}$ does not appear, since the probabilities of moving are only changed starting from 2, according to Definition 4.12).

Instead of (36), we obtain

$$\mathbb{E}\left(\ln\frac{\widetilde{Z}'_{t'_{n+1}}(2)}{Z'_{t'_n}(2)}\Big| \mathcal{F}'_{t'_n}\right) = \frac{[1 - (1 - \gamma_{t_k})\mathbb{1}_{n \in [k,2k)}\alpha_{t'_n}^{'-}(2)]^{-1}}{Z'_{t'_n}(2)} + \square\left(\frac{1}{2n^{2-3\varepsilon}}\right)$$

and when $n \in [k, 2k \wedge T')$,

$$\frac{[1 - (1 - \gamma_{t_k})\alpha_{t'_n}^{'-}(2)]^{-1}}{Z'_{t'_n}(2)} - \frac{\alpha_{t'_n}^{'+}(2)^{-1}}{Z'_{t'_n}(2)}$$

$$= \frac{\alpha_{t'_n}'^+(2)^{-1}}{Z'_{t'_n}(2)} \left[ \left( 1 + \frac{\gamma_{t_k}\alpha_{t'_n}'^-(2)}{\alpha_{t'_n}'^+(2)} \right)^{-1} - 1 \right]$$

$$\leq -\frac{\gamma_{t_k}\alpha_{t'_n}'^-(2)}{2Z'_{t'_n}(2)\alpha_{t'_n}'^+(2)^2} \leq -\frac{\gamma_{t_k}\alpha_{t_k}^-(2)}{32k} \leq -\frac{gW(k)}{64k},$$

using in the first inequality that $(1+x)^{-1} \leq 1 - x/2$ for $x \in [0,1]$ and using in the second inequality that $Z'_{t'_n}(2) \leq 4k$ (see the first paragraph of this section, and similarly $Z'_{t'_n}(3) \leq 4k$, which implies $\alpha_{t'_n}'^-(2) \geq \alpha_{t_k}^-(2)/4$ [recall that $\alpha_{t_k}^-(2) = \alpha_{t_k}'^-(2)$]. The rest of the proof is analogous to the proof of Assumption (H1), so that we obtain (27) and (28).

5.1.4. *Assumption* (H3) *of Proposition* 4.1 *holds.* Let us use Remark 4.1, with $R_n := R_n$ and $S_n := Z_n(2)Z_n(3)$. First observe that $R'_{t'_n} \geq R_{t_n}$, since $\mathcal{M}' \gg \mathcal{M}$. This comes from the fact that $\mathcal{M}'$ tends to go more to the right than $\mathcal{M}$, so that $R_{t_n}$, which increases only with visits from 5 to 4 and decreases with visits from 0 to 1, is larger for $\mathcal{M}'$ than for $\mathcal{M}$. More precisely, $\mathcal{M}' \gg \mathcal{M}$ implies on one hand that $Z'^+_{t'_n}(4) \geq Z^+_{t_n}(4)$ and $Z'^-_{t'_n}(1) \leq Z^-_{t_n}(1)$, and on the other hand by (24) that $R_t \leq R'_{t'}$ for all $t, t' \in \mathbb{N}$ such that $X_t = X'_{t'}$, $Z'^+_{t'}(4) \geq Z^+_t(4)$ and $Z'^-_{t'}(1) \leq Z^-_t(1)$. Second, if $\varepsilon < \mathrm{Cst}$,

$$\left| \frac{Z_{t_n}(2)Z_{t_n}(3)}{Z'_{t'_n}(2)Z'_{t'_n}(3)} - 1 \right| \leq 2 \left| \ln \frac{Z_{t_n}(2)Z_{t_n}(3)}{Z'_{t'_n}(2)Z'_{t'_n}(3)} \right|$$

$$\leq 2 \left( \ln \frac{Z'_{t'_n}(3)}{Z_{t_n}(3)} + \ln \frac{Z_{t_n}(2)}{Z'_{t'_n}(2)} \right) = 2|z'_n - z_n|.$$

The first inequality follows from $|x| \leq 2|\ln(1+x)|$ for $|x| < 1$: Indeed, $Z_{t_n}(2)$, $Z_{t_n}(3)$, $Z'_{t'_n}(2)$ and $Z'_{t'_n}(3)$ are close to $n$ (e.g., $Z_{t_n}(2) \in [Z^+_{t_n}(2), Z^+_{t_n}(2) + \varepsilon Z_{t_n}(2)]$ and $Z_{t_n}(2) = n$; see Section 5.1.3), so that $Z'_{t'_n}(2)/Z_{t_n}(2)$ and $Z'_{t'_n}(3)/Z_{t_n}(3)$ are close to 1. The second inequality is a consequence of $Z'_{t'_n}(3) \geq Z_{t_n}(3)$ and $Z'_{t'_n}(2) \leq Z_{t_n}(2)$, which follows again from $\mathcal{M}' \gg \mathcal{M}$.

5.2. *Proof of Lemma* 2.7. It suffices to prove that $\mathbb{P}(\Omega_1(0)) = 0$, because the problem is translation-invariant. We prove a preliminary result in Section 5.2.1 and then show the result in Section 5.2.2, using Proposition 4.1.

5.2.1. *A preliminary result.* The following lemma implies that a.s. on $\Omega_1(0)$, $Z_n(0) \vee Z_n(4) \vee Z_n(8)/Z_n(2) \wedge Z_n(6)$ tends to 0.

LEMMA 5.2. *One has*

$$\Omega_1(0) \subset \bigcap_{i=0,4,8} \left\{ \lim_{n\to\infty} \frac{Z_n(i)}{Z_n(2)} = 0 \right\}.$$

PROOF.    Let us prove, for instance, that $Z_n(0)/Z_n(2) \to 0$ on $\Omega_1(0)$; the proofs of the other statements are similar (using also Lemma 2.9). Observe that, by Proposition 3.1(c) and Remark 3.1,

$$\Omega_1(0) \subset \Upsilon(0,4) \cap \{Z_\infty(2) = \infty\} \subset \{\bar{Y}_\infty^-(1) < \infty\} \cap \{\beta_\infty^-(2) = \alpha_\infty^-(2) > 0\}.$$

Assume that $\Omega_1(0)$ holds. It suffices to prove $Z_n(0)/Z_n(1) \to 0$. There exist a.s. $a \in (0,1)$ and $p \in \mathbb{N}$ such that, for all $n \geq p$, $Z_n(1) \geq aZ_n(2)$. This implies, for all $n \geq p$,

$$\bar{Y}_\infty^-(1) - \bar{Y}_{n-1}^-(1) = \sum_{k \geq n} \frac{\mathbb{1}_{X_k=1}}{Z_k(1)} \frac{Z_k(0)}{Z_k(0) + Z_k(2)} \geq \sum_{k \geq n} \frac{\mathbb{1}_{X_k=1}}{Z_k(1)} \frac{Z_n(0)}{Z_n(0) + Z_k(2)}$$

$$\geq aZ_n(0) \sum_{k \geq n} \frac{\mathbb{1}_{X_k=1}}{(Z_k(1) + Z_n(0))^2} \geq a\frac{Z_n(0)}{Z_n(0) + Z_n(1) + 1}. \quad \square$$

5.2.2. *Application of Proposition* 4.1.    We apply Proposition 4.1 and use its notation to prove $\mathbb{P}(\Omega_1(0)) = 0$. Let us first introduce our choice of variables that appear in this lemma. Let us define a sequence $(t_n)_{n \in \mathbb{N}}$:

$$t_n := \inf\{m \in \mathbb{N} \text{ s.t. } Z_m(4) \geq n\}.$$

For $n \geq 2$, $t_n$ is the time of the $(n-1)$th visit to site 4.

For all $n \in \mathbb{N}$, let

$$R_n = Z_n(5) + Z_n(7) - (Z_n(1) + Z_n(3)).$$

Let $(y_n)_{n \in \mathbb{N}}$ and $(z_n)_{n \in \mathbb{N}}$ be the $(\mathcal{T}_{t_n})$-adapted processes defined by

$$y_n := \frac{R_{t_n}}{n(Z_{t_n}(3) + Z_{t_n}(5))}, \qquad z_n := \ln \frac{Z_{t_n}(6)}{Z_{t_n}(2)}$$

if $t_n < \infty$ and by $y_n = z_n := 0$ otherwise. Given $\varepsilon$, $a > 0$, for all $m \in \mathbb{N}$, let $T_1^{a,m,\varepsilon}$ and $T_2^{m,\varepsilon}$ be the $\mathbb{F}_{(t_n)_{n \in \mathbb{N}}}$ stopping times defined by

$$T_1^{a,m,\varepsilon} := \inf\Big\{n \geq m \text{ s.t. } Z_{t_n}(6)/Z_{t_n}(2) \notin [1-\varepsilon, 1+\varepsilon]$$

$$\text{or } \sup_{j \in \{0,4,8\}} Z_{t_n}(j) > \varepsilon Z_{t_n}(2) \text{ or } \alpha_{t_n}^+(2) \wedge \alpha_{t_n}^-(6) \leq a$$

$$\text{or } \inf_{j \in [m,n]} [\ln(Z_{t_n}(5)/Z_{t_j}(5)) - 2a^{-1}\ln(Z_{t_n}(4)/Z_{t_j}(4))] \geq \varepsilon\Big\},$$

$$T_2^{m,\varepsilon} := \inf\{n \geq m \text{ s.t. } t_n = \infty$$

$$\text{or } \exists y \in [1,7] \text{ s.t. } Z_{t_n}(y) - Z_{t_{n-1}}(y) \geq Z_{t_{n-1}}(y)n^{\varepsilon-1}\}$$

and $T^{a,m,\varepsilon} := T_1^{a,m,\varepsilon} \wedge T_2^{m,\varepsilon}$.

LEMMA 5.3. *For all $\varepsilon' > 0$,*

$$\Omega_1(0) \subset \{\lim z_n = 0\} \cap \left\{ \sum_{n \in \mathbb{N}} |y_n| \mathbb{1}_{y_n z_n \leq 0} < \infty \right\} \cap \left( \bigcup_{\substack{a > 0, \ m \in \mathbb{N} \\ \varepsilon < \varepsilon' \wedge 2^{-2a-1-3}}} \{T^{a,m,\varepsilon} = \infty\} \right).$$

PROOF. Let $\varepsilon > 0$ and suppose $\Omega_1(0)$ holds. The existence of $m \in \mathbb{N}$, $a > 0$ and $\varepsilon < \varepsilon' \wedge 2^{-2a-1-3}$ such that $T_1^{a,m,\varepsilon} = \infty$ follows from Lemmas 2.9 and 5.2, Remark 3.1 and Corollary 3.2(iii). The proof of the existence of $m \in \mathbb{N}$ such that $T_2^{m,\varepsilon} = \infty$ follows from an argument very similar to the proof of (17) [which gives the estimate of $Z_{t_n}(y) - Z_{t_{n-1}}(y)$ for $y = 2$].

To estimate the sum of $|y_n| \mathbb{1}_{y_n z_n \leq 0}$, observe that $y_n z_n \leq 0$ implies $R_{t_n}(Z_{t_n}(6) - Z_{t_n}(2)) \leq 0$; hence,

$$|R_{t_n}| \leq |R_{t_n} - (Z_{t_n}(6) - Z_{t_n}(2))| \leq Z_{t_n}(0) + Z_{t_n}(4) + Z_{t_n}(8).$$

Therefore, it suffices to prove that

$$\sum_{n \in \mathbb{N}} \frac{\mathbb{1}_{\{X_n = 4\}}}{Z_n(4)} \frac{Z_n(0)}{Z_n(3) + Z_n(5)} < \infty, \qquad \sum_{n \in \mathbb{N}} \frac{\mathbb{1}_{\{X_n = 4\}}}{Z_n(4)} \frac{Z_n(8)}{Z_n(3) + Z_n(5)} < \infty,$$

since the sum involving $Z_{t_n}(4)$ is obviously finite on $\Omega_1(0) \subset \{Y_\infty(4) < \infty\}$. Whereas $\Omega_1(0)$ is symmetric with respect to 4, let us prove the first inequality:

$$\sum_{n \in \mathbb{N}} \frac{\mathbb{1}_{\{X_n = 4\}}}{Z_n(4)} \frac{Z_n(0)}{Z_n(3) + Z_n(5)}$$

$$\approx \sum_{n \in \mathbb{N}} \frac{\mathbb{1}_{\{X_n = 4, X_{n+2} = 2\}}}{Z_n(4)} \frac{Z_n(0)}{Z_n(3)} \approx \sum_{n \in \mathbb{N}} \frac{\mathbb{1}_{\{X_n = 2, X_{n+2} = 4\}}}{Z_n(4)} \frac{Z_n(0)}{Z_n(3)}$$

$$\approx \sum_{n \in \mathbb{N}} \frac{\mathbb{1}_{\{X_n = 2\}}}{Z_n(2)} \frac{Z_n(0)}{Z_n(1) + Z_n(3)} \approx \sum_{n \in \mathbb{N}} \frac{\mathbb{1}_{\{X_n = 0\}}}{Z_n(0)} \frac{Z_n(0)}{Z_n(-1) + Z_n(1)}$$

$$= \sum_{n \in \mathbb{N}} \frac{\mathbb{1}_{\{X_n = 0\}}}{Z_n(-1) + Z_n(1)}.$$

In the second equivalence, we use an argument similar to the proof of Proposition 3.1(c). In the last equivalence, we use the same principle as in the previous arguments. □

Lemma 5.3 implies, together with Remark 4.2, that it suffices to apply Proposition 4.1, $a$, $m$ and $\varepsilon < \varepsilon' < \text{Cst}$ being fixed, with $(t_n)_{n \in \mathbb{N}}$, $(y_n)_{n \in \mathbb{N}}$ and $(z_n)_{n \in \mathbb{N}}$ as defined below, to conclude that $\mathbb{P}(\Omega_1(0)) = 0$. Let us choose here the other constants that appear in the application of this lemma: $c := 1/8$, $d := 4a^{-1}$, $M := 2^{2a^{-1}+1}$, $x := 5$, $W(n) := 1/\sqrt{n}$ and $T := T^{a,m,\varepsilon}$.

We leave to the reader the proof that Assumption (H1) holds, since it is very similar to the proof of the same fact in Section 5.1.2. Concerning the proof of Assumption (H2), let us show that $2k < T$ implies $t_{2k} < U_{5,t_k,M}$. It follows directly from $2k = Z_{t_{2k}}(4) = Z_{t_k}(4)$ that $\alpha_{t_n}^-(5) \le 2\alpha_{t_k}^-(5)$ for all $n \le 2k$. On the other hand, if $n \le 2k < T$, then $\ln(Z_{t_n}(5)/Z_{t_{2k}}(5)) \le 2a^{-1}\ln 2 + 1$ if $\varepsilon \le 1$, which implies $Z_{t_n}(5) \le 2^{2a^{-1}+1}Z_{t_k}(5)$. The rest of the proof that Assumption (H2) is fulfilled is left to the reader.

Let us check Assumption (H3), using Remark 4.1. First, $\mathcal{M}' \gg \mathcal{M}$ implies $R'_{t'_n} \ge R_{t_n}$. Second, we need an estimate of

$$\left| \frac{Z_{t_n}(3) + Z_{t_n}(3)}{Z'_{t'_n}(3) + Z'_{t'_n}(5)} - 1 \right| \le \left| \frac{Z_{t_n}(3) + Z_{t_n}(5) - (Z'_{t'_n}(3) + Z'_{t'_n}(5))}{Z'_{t'_n}(3) + Z'_{t'_n}(5)} \right|.$$

Note that

$$Z_{t_n}(3) + Z_{t_n}(5) = n + Z_{t_n}^-(3) + Z_{t_n}^+(5) + \mathrm{Cst}(v_0),$$

$$Z'_{t'_n}(3) + Z'_{t'_n}(5) = n + Z'^-_{t'_n}(3) + Z'^+_{t'_n}(5) + \mathrm{Cst}(v_0)$$

and

$$Z'^+_{t'_n}(5) - Z_{t_n}^+(5) \le Z'_{t'_n}(6) - Z_{t_n}(6), \qquad Z_{t_n}^-(3) - Z'^-_{t'_n}(3) \le Z_{t_n}(2) - Z'_{t'_n}(2).$$

Indeed, the first inequality follows, for instance, from

$$Z_{t_n}(6) = Z_{t_n}^+(5) + Z_{t_n}^+(6) + \mathrm{Cst}(v_0), \qquad Z'_{t'_n}(6) = Z'^+_{t'_n}(5) + Z'^+_{t'_n}(6) + \mathrm{Cst}(v_0),$$

which implies

$$Z'_{t'_n}(6) - Z_{t_n}(6) = Z'^+_{t'_n}(5) - Z_{t_n}^+(5) + Z'^+_{t'_n}(6) - Z_{t_n}^+(6) \ge Z'^+_{t'_n}(5) - Z_{t_n}^+(5),$$

where the last inequality follows from $Z'_{t'_n}(y) \ge Z_{t_n}(y)$ and $Z'^+_{t'_n}(y) \ge Z_{t_n}^+(y)$ for $y \ge 4$, and $Z'_{t'_n}(y) \le Z_{t_n}(y)$ and $Z'^-_{t'_n}(y) \le Z_{t_n}^-(y)$ for $y \le 4$, as a consequence of $\mathcal{M}' \gg \mathcal{M}$.

To summarize, as long as $n < T \wedge T'$, if $\varepsilon < \mathrm{Cst}$,

$$\left| \frac{Z_{t_n}(3) + Z_{t_n}(5)}{Z'_{t'_n}(3) + Z'_{t'_n}(5)} - 1 \right| \le \frac{Z'_{t'_n}(6) - Z_{t_n}(6)}{Z'_{t'_n}(3) + Z'_{t'_n}(5)} + \frac{Z_{t_n}(2) - Z'_{t'_n}(2)}{Z'_{t'_n}(3) + Z'_{t'_n}(5)}$$

$$\le 3a^{-1}\left( \frac{Z_{t_n}(2)}{Z'_{t'_n}(2)} - \frac{Z_{t_n}(6)}{Z'_{t'_n}(6)} \right) \le 4a^{-1}|z'_n - z_n|.$$

5.3. *Proof of Lemma* 2.8. It suffices to prove that $\mathbb{P}(\Omega_2(0)) = 0$, because the problem is translation-invariant. We prove a preliminary result in Section 5.3.1 and then prove the result in Section 5.3.2, using Proposition 4.1.

5.3.1. *A preliminary result.* Let, for all $n \in \mathbb{N}$,

$$\Delta_n := \frac{1}{Z_n^+(4)} \frac{Z_n(4)Z_n(5)}{Z_n(2)}.$$

LEMMA 5.4. *One has*

$$\Omega_2(0) \subset \left\{ \exists\, \Delta_\infty := \lim_{n \to \infty} \Delta_n \right\}.$$

PROOF. Suppose $\Omega_2(0)$ holds. Then $Z_n(7)/Z_n(2) \to 1$ by Lemma 2.11 and $\beta_n^+(2) \to \alpha_\infty^+(2) > 0$ by Remark 3.1. Using Corollary 3.2(iv), we obtain that $Z_n(5)/Z_n(2)^{\alpha_\infty^-(7)}$ and $Z_n(4)/Z_n(2)^{\alpha_\infty^+(2)}$ converge to a positive value. Note that these statements imply $Z_n(5)/Z_n(3) = \beta_n^+(2)^{-1} Z_n(5)/Z_n(2) \to 0$.

Therefore, there exist a.s. finite random variables $\gamma_\infty^i$, $i \in \{1,2,3,4,5\}$, such that

$$
\begin{aligned}
Z_{n+1}^+(4) &= \sum_{k=0}^n \mathbb{1}_{\{X_k=4, X_{k+1}=5\}} \\
&\asymp \sum_{k=0}^n \mathbb{1}_{\{X_k=4\}} \frac{Z_k(5)}{Z_k(3)+Z_k(5)} \asymp \sum_{k=0}^n \mathbb{1}_{\{X_k=4\}} \frac{Z_k(5)}{Z_k(3)} \\
&\asymp \gamma_\infty^1 \sum_{k=0}^n \mathbb{1}_{\{X_k=4\}} Z_k(2)^{\alpha_\infty^-(7)-1} \asymp \gamma_\infty^2 \sum_{k=0}^n \mathbb{1}_{\{X_k=4\}} Z_k(4)^{\alpha_\infty^-(2)^{-1}(\alpha_\infty^-(7)-1)} \\
&\asymp \gamma_\infty^3 Z_n(4)^{\alpha_\infty^+(2)^{-1}(\alpha_\infty^-(7)-1)} Z_n(4) \\
&\asymp \gamma_\infty^4 Z_n(2)^{\alpha_\infty^-(7)} \frac{Z_n(4)}{Z_n(2)} \asymp \gamma_\infty^5 \frac{Z_n(4)Z_n(5)}{Z_n(2)},
\end{aligned}
$$

where we use the conditional Borel–Cantelli lemma, Lemma A.1(i), in the first equivalence and use $\alpha_\infty^+(2)^{-1}(\alpha_\infty^-(7)-1) > -1$ [since $\alpha_\infty^-(2) < \alpha_\infty^-(7)$ by Lemma 2.10] in the fifth equivalence. $\quad\square$

5.3.2. *Application of Proposition 4.1.* We apply Proposition 4.1 and use its notation to prove $\mathbb{P}(\Omega_2(0)) = 0$. Let us first introduce our choice of variables that appear in this lemma. Let us define a sequence $(t_n)_{n \in \mathbb{N}}$ by

$$t_n = \inf\{m \in \mathbb{N}/Z_m^+(4) = n\}.$$

For all $n \in \mathbb{N}$, let

$$R_n := Z_n(6) + Z_n(8) - (Z_n(1) + Z_n(3)).$$

Let $(y_n)_{n \in \mathbb{N}}$ and $(z_n)_{n \in \mathbb{N}}$ be the $(\mathcal{T}_{t_n})$-adapted processes defined by

$$y_n := \frac{R_{t_n}}{Z_{t_n}(4)Z_{t_n}(5)}, \qquad z_n := \ln \frac{Z_{t_n}(7)}{Z_{t_n}(2)}$$

if $t_n < \infty$ and by $y_n = z_n := 0$ otherwise.

Given $\varepsilon$, $a$, $b > 0$, for all $m \in \mathbb{N}$, let $T_1^{a,b,m,\varepsilon}$ and $T_2^{m,\varepsilon}$ be the $\mathbb{F}_{(t_n)_{n \in \mathbb{N}}}$ stopping times defined by

$$T_1^{a,b,m,\varepsilon} := \inf\Big\{ n \geq m \text{ s.t. } Z_{t_n}(7)/Z_{t_n}(2) \notin [1 - \varepsilon, 1 + \varepsilon]$$

$$\text{or } \Delta_{t_n}/b \notin [1 - \varepsilon, 1 + \varepsilon] \text{ or } \alpha_{t_n}^+(2) \wedge \alpha_{t_n}^-(7) \leq a$$

$$\text{or } \sup_{m \leq k \leq n} \alpha_{t_n}^-(5)/\alpha_{t_k}^-(5) \geq 2$$

$$\text{or } \sup_{m \leq k \leq n} (Z_{t_n}(5)/n^{a-1})/(Z_{t_k}/k^{a-1}) \geq 2$$

$$\text{or } \sup_{v \in \{0,4,5,9\}} Z_{t_n}(v) > \varepsilon Z_{t_n}(2)$$

$$\text{or } \sup_{m \leq k \leq n} \sqrt{k}|(Y_{t_n}^+(6) - Y_{t_n}^-(6)) - (Y_{t_k}^+(6) - Y_{t_k}^-(6))| \geq 1 \Big\},$$

$$T_2^{m,\varepsilon} := \inf\{ n \geq m \text{ s.t. } t_n = \infty$$

$$\text{or } \exists\, y \in [1, 8] \text{ s.t. } Z_{t_n}(y) - Z_{t_{n-1}}(y) \geq Z_{t_{n-1}}(y) n^{\varepsilon - 1}\}$$

and $T^{a,b,m,\varepsilon} := T_1^{a,b,m,\varepsilon} \wedge T_2^{m,\varepsilon}$.

LEMMA 5.5.   *For all $\varepsilon' > 0$,*

$$\Omega_2(0) \subset \{\lim z_n = 0\} \cap \Big\{ \sum_{n \in \mathbb{N}} |y_n| \mathbb{1}_{y_n z_n \leq 0} < \infty \Big\} \cap \Big( \bigcup_{\substack{a,b > 0,\ m \in \mathbb{N} \\ \varepsilon < \varepsilon' \wedge a 2^{-2a-1-3}}} \{T^{a,b,m,\varepsilon} = \infty\} \Big).$$

PROOF.   Let $\varepsilon' > 0$ and suppose $\Omega_2(0)$ holds. Then $t_n < \infty$ for all $n \in \mathbb{N}$ and $z_n \to 0$. Let us prove the existence of $m \in \mathbb{N}$, $a$, $b > 0$ and $\varepsilon < \varepsilon' \wedge a 2^{-2a-1-3}$ such that $T_1^{a,b,m,\varepsilon} = \infty$. Lemma 2.11, Remark 3.1 and Lemma 5.4 imply that $Z_n(7)/Z_n(2) \to 1$, $\alpha_n^+(2) \to \alpha_\infty^+(2) \in (0, 1)$, $\alpha_n^-(7) \to \alpha_\infty^-(7) \in (0, 1)$ and $\Delta_n \to \Delta_\infty > 0$. Moreover, Corollary 3.2(iv) implies

$$\ln Z_n(4) \equiv \alpha_\infty^+(2) \ln Z_n(2), \qquad \ln Z_n(5) \equiv \alpha_\infty^-(7) \ln Z_n(7),$$

and similar estimates for $\ln Z_n(0)$ and $\ln Z_n(9)$, so that $\sup_{v \in \{0,4,5,9\}} Z_{t_n}(v)/Z_{t_n}(2) \to 0$ and $\sup_{n \geq k} \alpha_n^-(5)/\alpha_k^-(5) \to 1$ [there exists a.s. $\gamma_\infty > 0$ s.t. $\alpha_n^-(5) \asymp \gamma_\infty Z_n(2)^{\alpha_\infty^-(7)-1}$].

There exist random variables $\gamma_\infty^1$, $\gamma_\infty^2 > 0$ such that

$$Z_{t_n}(7) \asymp \gamma_\infty^1 n^{(\alpha_\infty^-(7) + \alpha_\infty^+(2) - 1)^{-1}}$$

[using $\Delta_{t_n} = Z_{t_n}(4)Z_{t_n}(5)/(nZ_{t_n}(2)) \to \Delta_\infty$] and

$$Z_{t_n}(5) \asymp \gamma_\infty^2 n^{\alpha_\infty^+(2)(\alpha_\infty^-(7)+\alpha_\infty^+(2)-1)^{-1}}.$$

Hence, there exists a.s. $a' > 0$ such that $\sup_{n \geq k}(Z_{t_n}(5)/n^{a^{-1}})/(Z_{t_k}(5)/k^{a^{-1}}) \to 1$ for all $a \in (0, a')$.

Next, observe that $\alpha_\infty^+(2)(\alpha_\infty^-(7) + \alpha_\infty^+(2) - 1)^{-1} > 1$ implies that there exists $\nu < 1/2$ such that $Z_{t_n}(5)^\nu/\sqrt{n} \to \infty$. Accordingly, using Proposition 3.1(a),

$$Y_{t_n}^+(6) \doteq Y_{t_n}^-(6) + o(Z_{t_n}(5)^{-\nu}) \doteq Y_{t_n}^-(6) + o(n^{-1/2}).$$

The proof of the existence of $m \in \mathbb{N}$ such that $T_2^{m,\varepsilon} = \infty$ relies on the observation that if $t > t_n$ is such that, for all $y \in [1, 9]$, $Z_t(y) - Z_{t_n}(y) \leq Z_{t_n}(y)n^{\varepsilon-1}$, then the probability to visit 4 at time $t + |X_t - 4|$ if $X_t \in [5, 9]$ (resp. to visit 5 at time $t + |X_t - 5|$ if $X_t \in [1, 4]$) is greater than a constant multiplied by $n/Z_{t_n}(X_t)$ (using $\Delta_{t_n} \to \Delta_\infty$). The result follows from an argument very similar to the proof of (17) in Section 3.4.

To estimate the sum of $|y_n|\mathbb{1}_{y_n z_n \leq 0}$, observe that $y_n z_n \leq 0$ implies $R_{t_n}(Z_{t_n}(7) - Z_{t_n}(2)) \leq 0$; hence,

$$|R_{t_n}| \leq |R_{t_n} - (Z_{t_n}(7) - Z_{t_n}(2))| \leq Z_{t_n}(0) + Z_{t_n}(4) + Z_{t_n}(5) + Z_{t_n}(9),$$

which enables us to conclude the proof using

$$\sum \frac{\sup_{v \in \{0,4,5,9\}} Z_{t_n}(v)}{Z_{t_n}(4)Z_{t_n}(5)} \preceq \sum \frac{\sup_{v \in \{0,4,5,9\}} Z_{t_n}(v)}{Z_{t_n}(2)} \frac{1}{n} < \infty. \qquad \square$$

Lemma 5.5 implies, together with Remark 4.2, that it suffices to apply Proposition 4.1, $a$, $b$, $m$ and $\varepsilon < \varepsilon' < \text{Cst}$ being fixed, with $(t_n)_{n \in \mathbb{N}}$, $(y_n)_{n \in \mathbb{N}}$ and $(z_n)_{n \in \mathbb{N}}$ as defined below, to conclude that $\mathbb{P}(\Omega_2(0)) = 0$. Let us choose here the other constants that appear in the application of this lemma: $c := \sqrt{a}/(4b^{3/2})$, $d := 20 \vee 4b^{-1}$, $M := 2^{a^{-1}+1}$, $x := 5$, $W(n) := 1/\sqrt{n}$ and $T := T^{a,b,m,\varepsilon}$.

We leave to the reader the proof that Assumption (H1) holds, since it is very similar to the proof of the same fact in Section 5.1.2.

Concerning the proof of Assumption (H2), let us show that $2k < T$ implies $t_{2k} < U_{5,t_k,M}$. For all $n \geq 2k < T$, $\alpha_{t_n}^-(5) \leq 2\alpha_{t_k}^-(5)$ and $Z_{t_n}(5) \leq 2(2k)^{a^{-1}} Z_{t_k}(5)/k^{a^{-1}} \leq 2^{a^{-1}+1}Z_{t_k}(5)$ by definition of $T$. The rest of the proof that Assumption (H2) is fulfilled is left to the reader.

Let us check Assumption (H3), using Remark 4.1 with $R_n := R_n$ and $S_n := Z_n(4)Z_n(5)$. First, $\mathcal{M}' \gg \mathcal{M}$ implies $R'_{t'_n} \geq R_{t_n}$. Second, we need an estimate of

$$\left| \frac{Z_{t_n}(4)Z_{t_n}(5)}{Z'_{t'_n}(4)Z'_{t'_n}(5)} - 1 \right| \leq 2\left( \ln \frac{Z_{t_n}(4)}{Z'_{t'_n}(4)} + \ln \frac{Z'_{t'_n}(5)}{Z_{t_n}(5)} \right).$$

Let us, for instance, upper bound $\ln Z'_{t'_n}(5)/Z_{t_n}(5)$. Assume $n \geq k \geq \mathrm{Cst}$ and let $\delta_k := Y^-_{t_k}(6) - Y^+_{t_k}(6)$. The proofs of Proposition 3.1(b) and (c) imply [using $Z_{t_n}(5) \geq Z^+_{t_n}(5) \geq n$, $Z_{t_n}(7) \geq (1-\varepsilon)Z_{t_n}(2) \geq (1-\varepsilon)Z_{t_n}(5)/\varepsilon \geq n$ and $n < T^{a,b,m,\varepsilon}_1$] that

$$
\begin{aligned}
\ln Z_{t_n}(5) &- \mathrm{Cst}(x, v_0) \\
&= Y^+_{t_n}(4) + Y^-_{t_n}(6) + \square(Z_{t_n}(5)^{-1}) \\
&= Y^+_{t_n}(4) + Y^+_{t_n}(6) + \delta_k + \square(2/\sqrt{k}) \\
&= Y^+_{t_n}(4) + \widetilde{Y}^+_{t_n}(6) + \delta_k + \square(3/\sqrt{k})
\end{aligned}
$$

and, similarly,

$$
\ln Z'_{t'_n}(5) - \mathrm{Cst}(x, v_0) = Y'^+_{t'_n}(4) + \widetilde{Y}'^+_{t'_n}(6) + \delta_k + \square(2/\sqrt{k})
$$

(recall that $\mathcal{M}$ and $\mathcal{M}'$ are the same until time $t_k$). Therefore,

$$
\tag{39} \ln \frac{Z'_{t'_n}(5)}{Z_{t_n}(5)} \leq Y'^+_{t'_n}(4) - Y^+_{t_n}(4) + \widetilde{Y}'^+_{t'_n}(6) - \widetilde{Y}^+_{t_n}(6) + \frac{4}{\sqrt{k}}.
$$

Now $Y'^+_{t'_n}(4) \leq Y^+_{t_n}(4)$ since $\mathcal{M}' \gg \mathcal{M}$. Let $u_n$ (resp. $u'_n$) be the $n$th visit time to site 7 for $\mathcal{M}$ (resp. $\mathcal{M}'$):

$$
\widetilde{Y}^+_{t_n}(6) = \sum_{k=1}^{Z_{t_n}(7)} \frac{\mathbb{1}_{\{X_{u_k+1}=6\}}}{k},
$$

$$
\widetilde{Y}'^+_{t'_n}(6) = \sum_{k=1}^{Z'_{t'_n}(7)} \frac{\mathbb{1}_{\{X'_{u'_k+1}=6\}}}{k} \leq \sum_{k=1}^{Z'_{t'_n}(7)} \frac{\mathbb{1}_{\{X_{u_k+1}=6\}}}{k}.
$$

In summary, (39) implies, assuming $k \geq \mathrm{Cst}$ and $\varepsilon < \mathrm{Cst}$,

$$
\ln \frac{Z'_{t'_n}(5)}{Z_{t_n}(5)} \leq \sum_{k=Z_{t_n}(7)+1}^{Z'_{t'_n}(7)} \frac{1}{k} \leq \ln \frac{Z'_{t'_n}(7)}{Z_{t_n}(7)} + \frac{1}{Z_{t_n}(7)} + \frac{4}{\sqrt{k}} \leq \ln \frac{Z'_{t'_n}(7)}{Z_{t_n}(7)} + \frac{5}{\sqrt{k}}.
$$

## APPENDIX: GENERAL MARTINGALE RESULTS

Let us recall the following theorem: The first part is due to Doob, while the second part is due to Neveu (see [6], Propositions VII-2.3 and VIII-2.4).

THEOREM A.1. *Let $(M_n)_{n\in\mathbb{N}}$ be a square integrable martingale. Let*

$$
\alpha_n = \mathbb{E}((M_{n+1} - M_n)^2 | \mathcal{F}_n) < \infty, \qquad \langle M \rangle_n = \sum_{k=0}^n \alpha_i.
$$

*Then, for all $r > 1/2$:*

(i) $\{\langle M \rangle_\infty < \infty\} \subset \{\exists\, M_\infty \in \mathbb{R}/M_n \to M_\infty\}$ *a.s.;*

(ii) $\{\langle M \rangle_\infty = \infty\} \subset \{M_n = o(\langle M \rangle_{n-1}^{1/2}(\ln\langle M \rangle_{n-1})^r)\}$ *a.s.*

We make use of the following generalized version of the conditional Borel–Cantelli lemma at various steps of the proof.

LEMMA A.1.    *Let* $\mathbb{G} = (\mathcal{G}_n)_{n\in\mathbb{N}}$ *be a filtration. Let* $(\xi_n)_{n\in\mathbb{N}}$ *[resp.* $(\beta_n)_{n\in\mathbb{N}}$*] be a* $\mathbb{G}$*-adapted sequence of random variables that take values in* $\mathbb{R}_+$ *and are bounded [resp. nondecreasing]. Let* $(\Gamma_n)_{n\in\mathbb{N}^*}$ *be a sequence of* $\mathbb{G}$*-adapted random sets. Let*

$$p_{k-1} = \mathbb{P}(\Gamma_k | \mathcal{G}_{k-1}),$$

$$\Phi_n = \sum_{k=1}^n \xi_{k-1}\mathbb{1}_{\Gamma_k}, \qquad \Phi_n^* = \sum_{k=1}^n \xi_{k-1}p_{k-1},$$

$$\delta_k = \xi_k^2 p_k(1-p_k), \qquad \Delta_n = \sum_{k=1}^n \delta_{k-1}.$$

*Then*:

(i) $\Phi_n \asymp \Phi_n^*$;

(ii) $\{\Delta_\infty < \infty\} \subset \{\Phi_n \equiv \Phi_n^*\}$;

(iii) $\{\sum_{k\in\mathbb{N}} \beta_k\delta_k < \infty\} \subset \{\Phi_n - \Phi_n^* \doteq O(1/\sqrt{\beta_n})\}$.

PROOF.    First observe that

$$M_n = \Phi_n - \Phi_n^*$$

is a martingale and that $\Delta_n = \langle M \rangle_{n-1}$. This implies claims (i) and (ii), using Theorem A.1 [note that $\Delta_n = O(\Phi_n^*)$, since $(\xi_n)_{n\in\mathbb{N}}$ is bounded].

Let us now prove (iii). Let us define the $\mathbb{G}$-adapted random processes

$$\Psi_n = \sum_{k=1}^n \beta_{k-1}^{1/2}\xi_{k-1}\mathbb{1}_{\Gamma_k}, \qquad \Psi_n^* = \sum_{k=1}^n \beta_{k-1}^{1/2}\xi_{k-1}p_{k-1}, \qquad R_n = \Psi_n - \Psi_n^*,$$

with the convention that $R_0 = 0$. Observe that $R_n$ is a martingale and that

$$\langle R \rangle_\infty = \sum_{k=0}^\infty \beta_k\delta_k < \infty,$$

which implies by Lemma A.1 that $R_n$ converges a.s. toward a r.v. $R_\infty \in \mathbb{R}$. Now, observe that

$$M_n = \sum_{k=1}^n \beta_{k-1}^{-1/2}(R_k - R_{k-1})$$

and, therefore,

$$M_\infty - M_n = \sum_{k=n+1}^{\infty} \beta_{k-1}^{-1/2}(R_k - R_{k-1})$$

$$= \sum_{k=n+1}^{\infty} (\beta_{k-1}^{-1/2} - \beta_k^{-1/2})R_k - \beta_n^{-1/2}R_n = O(1/\sqrt{\beta_n}). \qquad \square$$

The following lemma has some similarity to the conditional Borel–Cantelli lemma, although its proof is based on different arguments. Assuming some upper bounds on conditional probabilities for events $\Gamma_n$, depending on the number of times $\Gamma_n$ has arisen, $\Gamma_n$ a.s. holds only finitely often. The result is used in the proof of Lemma 2.10.

LEMMA A.2. *Let* $\mathbb{G} = (\mathcal{G}_n)_{n \in \mathbb{N}}$ *be a filtration, let* $(\gamma_n)_{n \in \mathbb{N}}$ *be a* $\mathbb{G}$-*adapted sequence that takes values in* $\mathbb{R}$ *such that* $\liminf \gamma_n > 0$ *a.s. and let* $(\Gamma_n)_{n \in \mathbb{N}^*}$ *be a sequence of* $\mathbb{G}$-*adapted events. For all* $n \in \mathbb{N}$, *let*

$$\tau_n = \sup\{k \le n \ s.t. \ \Gamma_k \ holds\}.$$

*Then*

$$(40) \quad \left\{\exists a \in \mathbb{R}_+^*, m \in \mathbb{N} \ s.t. \ \forall n \ge m, \mathbb{P}(\Gamma_{n+1}|\mathcal{G}_n) \le \frac{a\tau_n^{\gamma_{\tau_n}}}{n^{1+\gamma_{\tau_n}}}\right\} \subset \left\{\sum_{n=0}^{\infty} \mathbb{1}_{\Gamma_n} < \infty\right\}.$$

PROOF. For all $a, \varepsilon \in \mathbb{R}_+^*$ and $m \in \mathbb{N}$, let $T_{a,\varepsilon,m}$ be the stopping time $T_{a,\varepsilon,m} := \inf\{n \ge m \text{ s.t. } \gamma_n < \varepsilon \text{ or } \mathbb{P}(\Gamma_{n+1}|\mathcal{G}_n) > a\tau_n^{\gamma_{\tau_n}}/n^{1+\gamma_{\tau_n}}\}$. For all $n \in \mathbb{N}$, let $\Lambda_{n+1}^{a,\varepsilon,m} := \Gamma_{n+1} \cap \{T_{a,\varepsilon,m} > n\}$ and $\mathcal{V}_{a,\varepsilon,m} = \{\sum_{n=1}^{\infty} \mathbb{1}_{\Lambda_n^{a,\varepsilon,m}} < \infty\}$.

Let us prove that, for all $a, \varepsilon \in \mathbb{R}_+^*$ and $n \ge m$, $\mathbb{P}(\mathcal{V}_{a,\varepsilon,m}|\mathcal{G}_n) \ge \mathrm{Cst}(a,\varepsilon) > 0$. This enables us to conclude. Indeed, suppose this inequality holds. By a standard martingale theorem, $\mathbb{P}(\mathcal{V}_{a,\varepsilon,m}|\mathcal{G}_n) = \mathbb{E}(\mathbb{1}_{\mathcal{V}_{a,\varepsilon,m}}|\mathcal{G}_n) \to_{n \to \infty} \mathbb{E}(\mathbb{1}_{\mathcal{V}_{a,\varepsilon,m}}|\mathcal{G}_\infty) = \mathbb{1}_{\mathcal{V}_{a,\varepsilon,m}}$ since $\mathcal{V}_{a,\varepsilon,m} \in \mathcal{G}_\infty$. Therefore, $\mathbb{1}_{\mathcal{V}_{a,\varepsilon,m}} \ge \mathrm{Cst}(a,\varepsilon)$ a.s. and $\mathbb{P}(\mathcal{V}_{a,\varepsilon,m}) = 1$ (for all $a, \varepsilon \in \mathbb{R}_+^*$ and $m \in \mathbb{N}$). We deduce that a.s. on $\{\exists a, \varepsilon \in \mathbb{R}_+^*, m \in \mathbb{N} \text{ s.t. } T_{a,\varepsilon,m} = \infty\}$, $\Gamma_n$ only occurs finitely often, which proves the lemma.

Fix $a, \varepsilon \in \mathbb{R}_+^*$, and $m \in \mathbb{N}$. For simplicity, we assume $m \ge 2\sup(1,a)$, and $\gamma_n \le 1$ for all $n \in \mathbb{N}$ [the overestimate of $\mathbb{P}(\Gamma_{n+1}|\mathcal{G}_n)$ remains true if we replace $\gamma_n$ by $\gamma_n \wedge 1$]. Given $n \ge m$, let us estimate

$$\mathbb{P}(\mathcal{V}_{a,\varepsilon,m}|\mathcal{G}_n) \ge \mathbb{P}\left(\bigcap_{k=n+1}^{\infty} (\Lambda_k^{a,\varepsilon,m})^c \Big| \mathcal{G}_n\right) \ge \prod_{k=n}^{\infty} \left(1 - \frac{an^{\gamma_{\tau_n}}}{k^{1+\gamma_{\tau_n}}}\right)$$

$$\ge \exp\left(-2an^{\gamma_{\tau_n}} \sum_{k=n}^{\infty} \frac{1}{k^{1+\gamma_{\tau_n}}}\right) \ge \exp\left(-\frac{2an^{\gamma_{\tau_n}}}{\gamma_{\tau_n}(n-1)^{\gamma_{\tau_n}}}\right)$$

$$\ge \exp\left(-\frac{4a}{\varepsilon}\right) > 0,$$

where we use that $1 - x \geq \exp(-2x)$ for all $x \in [0, 1/2]$ (with $x = an^{\gamma_{\tau_n}}/k^{1+\gamma_{\tau_n}} \leq an^{-1} \leq 1/2$ since $m \geq 2a$) and $(n-1)/n \geq 1/2$ (since $m \geq 2$). This enables us to conclude. $\square$

**Acknowledgments.** I am very grateful to my Ph.D. advisor, M. Benaïm, for helpful discussions and for inspiring the choice of this subject. I have also greatly benefitted from remarks by M. Duflo on estimates, which have been very useful to improve this article. I am indebted to the Associate Editor for his detailed report and for his important advice concerning this paper. I also thank an anonymous referee for very carefully reading the paper.

Faculté des Sciences
Institut de Mathématiques
Rue Emile-Argand 11
Case Postale 2
CH-2007 Neuchâtel
Switzerland
e-mail: pierre.tarres@unine.ch
and
Laboratoire de Statistique et Probabilits
UMR CNRS C5583—Universit Paul Sabatier
Btiment 1R1
118 Route de Narbonne
31062 Toulouse Cedex 4
France
e-mail: tarres@math.ups-tlse.fr